\documentclass{article}
\usepackage{algorithm, algpseudocode, amsmath, amssymb, float, geometry, lmodern, mathtools, subcaption}

\usepackage{color, graphicx, epstopdf}
\usepackage[T1]{fontenc}
\usepackage[utf8]{inputenc}
\usepackage{soul}
\usepackage{authblk}

\def\sb{\textbf{\emph{S}}}

\def\RR{\mathbb{R}}
\def\PP{\mathbb{P}}

\def\TT{\mathbb{T}}
\def\Cc{\mathcal{C}}

\def\Pc{\mathcal{P}}
\def\Sc{\mathcal{S}}
\def\Tc{\mathcal{T}}
\def\Kc{\mathcal{K}}

\def\spvecA#1;{\if;#1;\else #1\cr \expandafter \spvecA \fi}

\newtheorem{remark}{Remark}

\newcommand{\Bezier}{B\'{e}zier\ }
\newcommand{\GIFT}{GIFT~\!}
\newcommand{\vB}{\boldsymbol{B}}
\newcommand{\vC}{\boldsymbol{C}}
\newcommand{\vD}{\boldsymbol{D}}
\newcommand{\veps}{\boldsymbol{\epsilon}}
\newcommand{\vf}{\boldsymbol{f}}
\newcommand{\vI}{\boldsymbol{I}}
\newcommand{\vJ}{\boldsymbol{J}}
\newcommand{\vk}{\boldsymbol{k}}
\newcommand{\vK}{\boldsymbol{K}}
\newcommand{\vP}{\mathcal{P}}
\newcommand{\vu}{\boldsymbol{u}}
\newcommand{\vx}{\boldsymbol{x}}
\newcommand{\vxi}{\boldsymbol{\xi}}
\newlength{\picw}

\begin{document}
\author[1]{Elena Atroshchenko}
\author[2]{Gang Xu \thanks{Corresponding author. \\Email addresses: xugangzju@gmail.com; gxu@hdu.edu.cn (Gang Xu), eatroshchenko@ing.uchile.cl (Elena Atroshchenko), tomar.sk@iitkalumni.org (Satyendra Tomar), stephane.bordas@uni.lu (St\'ephane P.A. Bordas)}}
\author[3]{Satyendra Tomar}
\author[3, 4, 5]{St\'ephane P.A. Bordas}
\affil[1] {University of Chile, Department of Mechanical Engineering, Santiago, 8370448, Chile}
\affil[2] {Hangzhou Dianzi University, Hangzhou 310018, P.R.China}
\affil[3] {Institute of Computational Engineering, University of Luxembourg, Faculty of Sciences Communication and Technology, Luxembourg}
\affil[4] {Institute of Mechanics and Advanced Materials, School of Engineering, Cardiff University, UK}
\affil[5] {Intelligent Systems for Medicine Laboratory
University of Western Australia, Perth, Australia}
\date{}
\title{Weakening the tight coupling between geometry and simulation in isogeometric analysis: from sub- and super- geometric analysis to Geometry Independent Field approximaTion (GIFT)}

\maketitle
\begin{abstract}

This paper presents an approach to generalize the concept of isogeometric analysis (IGA) by allowing different spaces for parameterization of the computational domain and for approximation of the solution field. The method inherits the main advantage of isogeometric analysis, i.e. preserves the original, exact CAD geometry (for example, given by NURBS), but allows pairing it with an approximation space which is more suitable/flexible for analysis, for example, T-splines, LR-splines, (truncated) hierarchical B-splines, and PHT-splines. This generalization offers the advantage of adaptive local refinement without the need to re-parameterize the domain, and therefore without weakening the link with the CAD model. We demonstrate the use of the method with different choices of the geometry and field splines, and show that, despite the failure of the standard patch test, the optimum convergence rate is achieved for non-nested spaces.

\end{abstract}

\section{Introduction}

\subsection{Contribution}
We present an approach which enables the use of separate approximation spaces for the field and geometry in isogeometric finite element methods (IGAFEM). For example, coarse NURBS (non-uniform rational B-splines) approximations can be used for the geometry, and locally adapted PHT splines can be used for the field variables. This endows the method with the flexibility to locally enrich the approximation space for the field variables without modifying the spline space used for parametrization of the geometry of the domain. We verify the approach with various parameterizations of the geometry and approximations of the field variables, and on two and three dimensional geometries. We give a detailed mathematical explanation on the ability of the resulting method to pass the patch test, and to converge with optimal rate. Geometry Independent Field approximaTion (GIFT), as we coined the approach, represents a simple and computationally efficient alternative to IGA, which preserves tight integration with CAD, and allows adaptive local 
refinement of the solution field.

\subsection{Background}
Isogeometric analysis (IGA) was introduced in \cite{Hughes20054135} to establish a direct link between the computer-aided-design (CAD) and analysis. Over the last decade, fuelled by rapid developments in  computer graphics and CAD, and due to a number of advantages offered by spline basis functions over standard (Lagrange) finite element analysis (FEA), such as higher continuity, exact representation of the geometry, simplified integration and their behaviour in dynamics, IGA has created a large amount of interest in computational engineering and science. 

The method has found many applications in various areas such as structural vibrations \cite{Reali}, fluids-structure interaction \cite{Bazilevs2008}, shell analysis \cite{Benson2010276},  and fracture mechanics \cite{de2011x}. Although the finite element version of IGA suffers from significant difficulties associated with generating volumetric parameterisations directly from CAD, coupling IGA with boundary element methods enables to perform computations directly from the CAD description of the boundary of the domain \cite{beer2015advancedbook}. Applications include stress analysis  \cite{Simpson201287,lian2013stress,Scott2013197},  shape optimization \cite{seo2010shape, lian2015implementation, lian2017shape}, fracture mechanics  \cite{peng2014damage,peng2016linear,xuan2016isogeometric,peng2017isogeometric}, geomechanics \cite{beer2015simple, marussig2015fast, beer2016advanced, beer2016isogeometric, GeneralizedIGABEMBeer2016, TrimmedBeer2016, Zechner2016212, beer2017isogeometric}, acoustics \cite{Simpson2014265,khajah2016isogeometric}, and electromagnetics \cite{buffa2010isogeometric}. A  detailed overview of the recent work in the field can be found in \cite{Nguyen201589} and in \cite{lian2012recent}.
Other approaches, which also aim at facilitating the transition between CAD and analysis have been developed since the late 1990's, see e.g., \cite{NME:NME292,Turco}, and include (in order of appearance):
\begin{description}
\item[Implicit/immersed boundary definitions] with the extended finite element method and level sets \cite{belytschko2003structured}, \cite{moes2003computational} and \cite{moumnassi2011finite,moumnassi2014analysis}, which allows an implicit definition of arbitrary solids through constructive geometry and different level set functions. A promising approach, known as the finite cell method was proopsed in 2007 by \cite{parvizian2007finite}.
\item[Subdivision surfaces] \cite{cirak2000subdivision,cirak2011subdivision}, where the idea is to use subdivision surfaces as a common representation for geometric modeling and numerical simulation in a unified framework. Subdivision surfaces is a flexible and efficient tool for modeling surfaces with arbitrary topology, avoiding many problems inherent in traditional geometry modeling approach with tensor-product patches. The corresponding solvers are highly scalable, and  provide an efficient computational tool for numerical simulation and design optimization.
\item[NURBS-enhanced finite element methods] \cite{sevilla2008nurbs, sevilla2011nurbs, sevilla20113d}, where the idea is to construct new finite elements with at least one edge represented as a NURBS or B-spline. This, therefore, enables to use the CAD geometry directly, but to retain the standard finite element method everywhere except in the vicinity of the boundary.
\end{description}

\subsection{NURBS: Advantages and limitations}
NURBS are the most common type of splines used in the CAD industry to describe geometry. NURBS shape functions are defined by means of B-splines (piece-wise polynomial functions defined over a knot vector) and a set of control points with associated weights. The main feature of NURBS is the ability to exactly represent conic sections, e.g. ellipses, parabolas, hyperbolas, etc., which are widely used in computer-aided-design.  NURBS basis functions also posses all the properties of the standard FEM basis functions, such as compact support, linear independence, and partition of unity. Moreover, the high order continuity ($C^{p - k}$, with $p$ the polynomial order and $k$ the knot repetition)  of NURBS basis functions across elements facilitates the solution of partial differential equations (PDEs) of arbitrary high order. However, NURBS are based on tensor-product structure, which does not facilitate local refinement (see Fig.~\ref{global_ref}), NURBS suffer from inability to produce the watertight geometries for  general shapes, which imposes difficulties for mesh generation and requires coupling algorithms for multi-patch geometries \cite{nguyen2014nitsche, coox2016robust, du2015nitsche,xu13cada,xu13jcp,xu17cad}.

\subsection{Motivation}
\begin{figure}[!ht]
\begin{center}
\includegraphics[width=\textwidth]{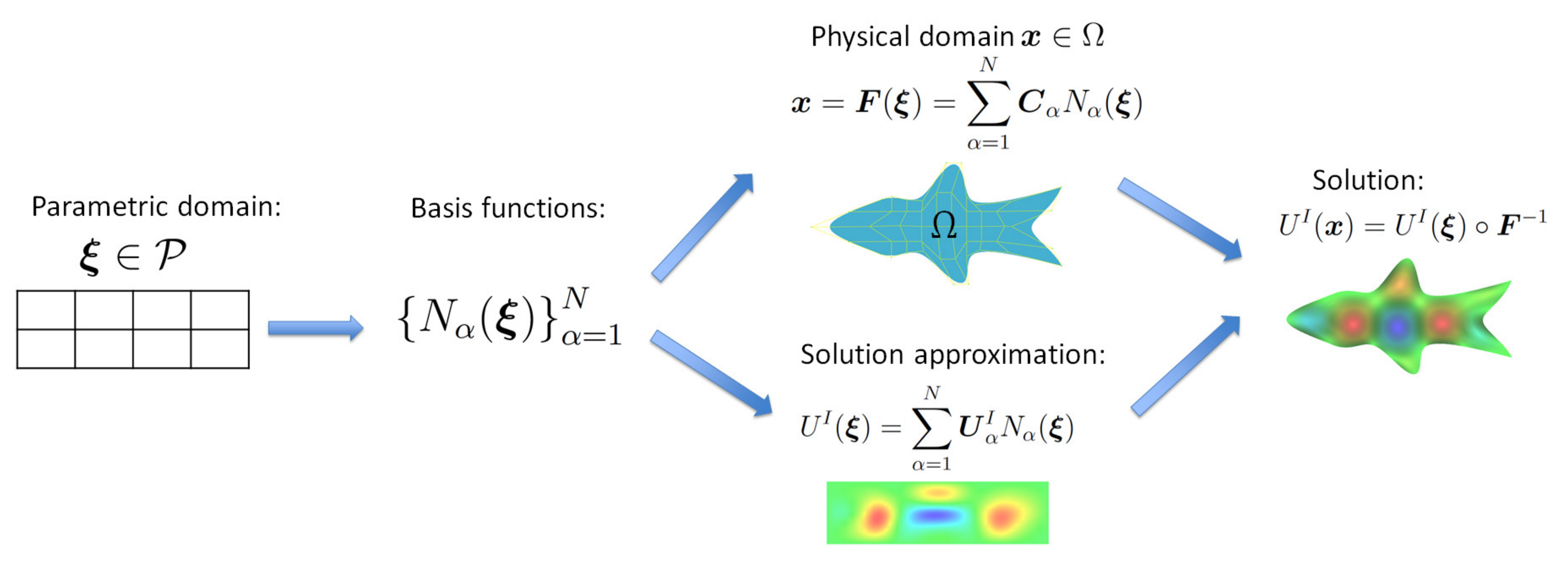}
\caption{Main idea of IGA: the same shape functions are used for geometry parameterization and solution approximation.}\label{main_idea_IGA}
\end{center}
\end{figure}

The main idea of IGA is to use the same shape functions (splines) for both: parameterization of the geometry (computational domain) and approximation of the unknown solution, see Fig.~\ref{main_idea_IGA}. This is a clear advantage in engineering design for example, as any modification of the CAD geometry is directly inherited by the approximation of the unknown fields so that the mesh need not be regenerated at each iteration of the geometrical design process. Because of the nature of the approximations most commonly used in geometrical design (NURBS), which are based on tensor-product constructions, the potential power of IGA is not always realized in practice. Indeed, local refinement is not natural in NURBS approximations. Moreover, situations usually arise when the approximation used for the field variables must be refined locally, whilst keeping the same (coarser) geometrical approximation. This is the case in shape optimization \cite{seo2010shape, lian2015implementation, lian2017shape}, or for problems with  singularities, boundary layers or steep gradients. 

It appears, therefore, that it is useful to develop more general methods where the approximation space used for the geometry and that used for the field variables be decoupled. The investigation of the behaviour of such methods for linear elasticity and Laplace equation is the focus of this paper. 

We base our work on the developments made by the team of Gernot Beer in two important papers for isogeometric boundary element methods: \cite{GeneralizedIGABEMBeer2016,marussig2015fast}, which were further developed in a series of papers targeting geomechanics problems \cite{beer2015simple, beer2016advanced,beer2016isogeometric,GeneralizedIGABEMBeer2016,TrimmedBeer2016,Zechner2016212,beer2017isogeometric}. 

The most important short-comings of IGA, which served as motivation to investigate the possibility to decouple the geometry and the solution fields, can be summarized as follows:
\begin{enumerate}
\item The $h$- and $p$-refinement of the original exact geometry is redundant in IGA. The computational savings on non-refining the geometry were estimated in \cite{GeneralizedIGABEMBeer2016}, and a speed up factor of up to 2 was achieved in the context of BEM. It is expected that for large industrial problems, the savings will be even more significant.
\item NURBS, the de facto standard in CAD community, do not offer local refinement. Therefore, in order to use any other locally-refined splines with the original CAD-model given by NURBS, staying within the iso-parametric concept will require re-parameterization of the original geometry. Note that some locally-refined splines (e.g., T-splines) allow exact representation of the original geometry, but others (e.g., polynomial based splines) do not offer this advantage. Therefore, such operations (particularly, for efficient local refinement) are not only time-consuming, but may also lead to loosing the geometry exactness, and may introduce the need to communicate with the CAD model at every step of the solution refinement process.
\end{enumerate}

\subsection{Alternatives to NURBS}
In many applications, the efficiency of a numerical method can be drastically improved by restricting the refinement to certain  areas (for example, where the solution exhibits high gradients). Therefore, the local refinement issue in IGA is currently an active research topic. A recent survey on locally refinable splines is given in \cite{Li2016}. In the authors' opinion, the most commonly used alternatives which allow local refinement, are THB-splines,  T-splines and PHT-splines (polynomial splines over hierarchical T-meshes), which are briefly discussed below. 

\subsubsection{Hierarchical B-splines and Truncated Hierarchical B-splines (THB-splines)}

Hierarchical B-splines were first introduced in \cite{HierarchicalBSplines1} for surface fitting. They posses such properties as maximum continuity on each refinement level, and linear independence of the basis functions. Later, in \cite{HierarchicalBSplines2}, these shape functions were modified to reduce the local support and assure the partition of unity property. These splines were called  Truncated Hierarchical B- (THB-) splines. The theoretical background, as well as applications to isogeometrical analysis, can be found in \cite{THB0} and \cite{THB1}. 

\subsubsection{T-splines}
\begin{figure}[!ht]
\centering
\begin{subfigure}[b]{0.325\textwidth}
\includegraphics[width = \textwidth]{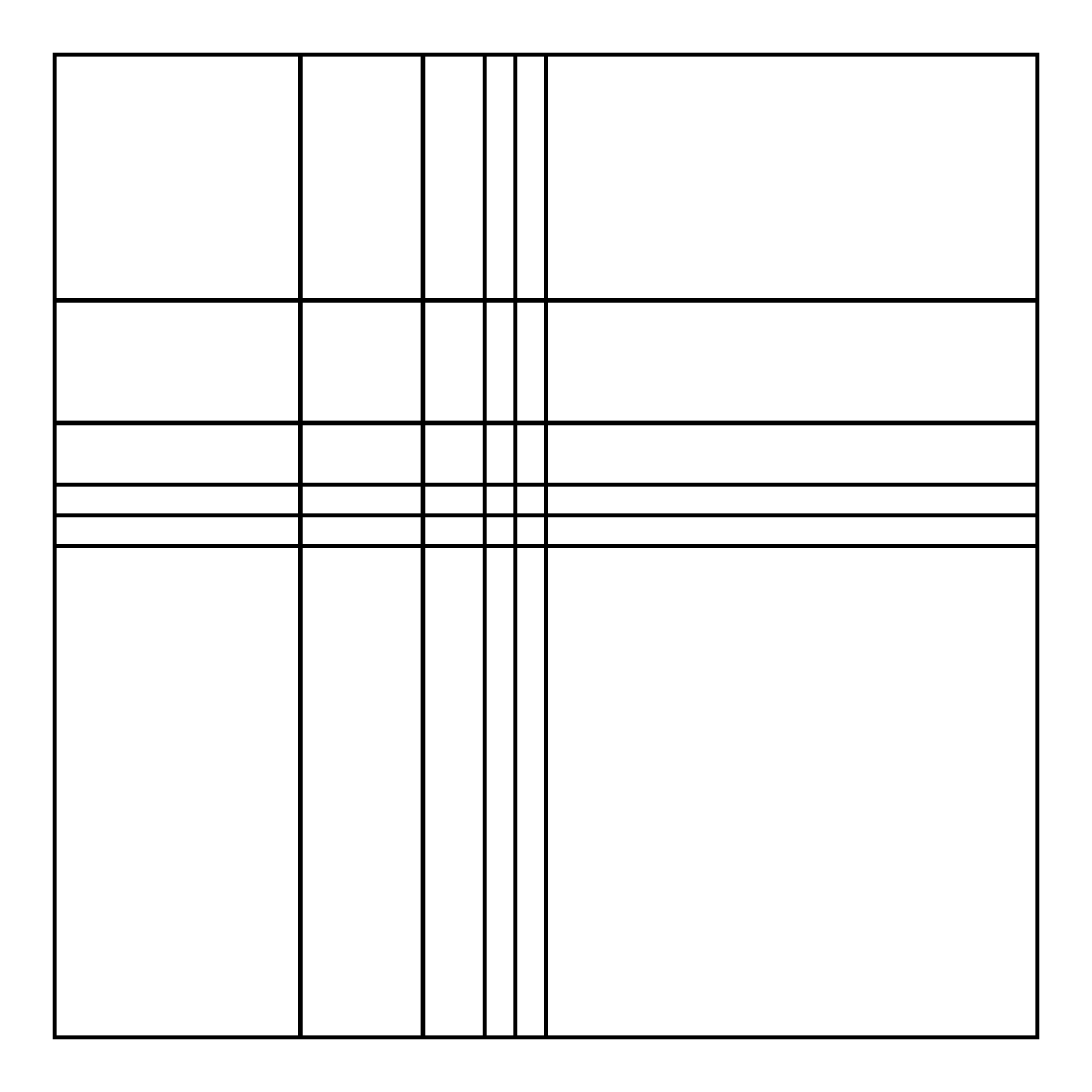}  
\caption{Tensor-product mesh in the parametric space (global refinement).}\label{global_ref}
\end{subfigure}
\quad
\begin{subfigure}[b]{0.325\textwidth}
\includegraphics[width = \textwidth]{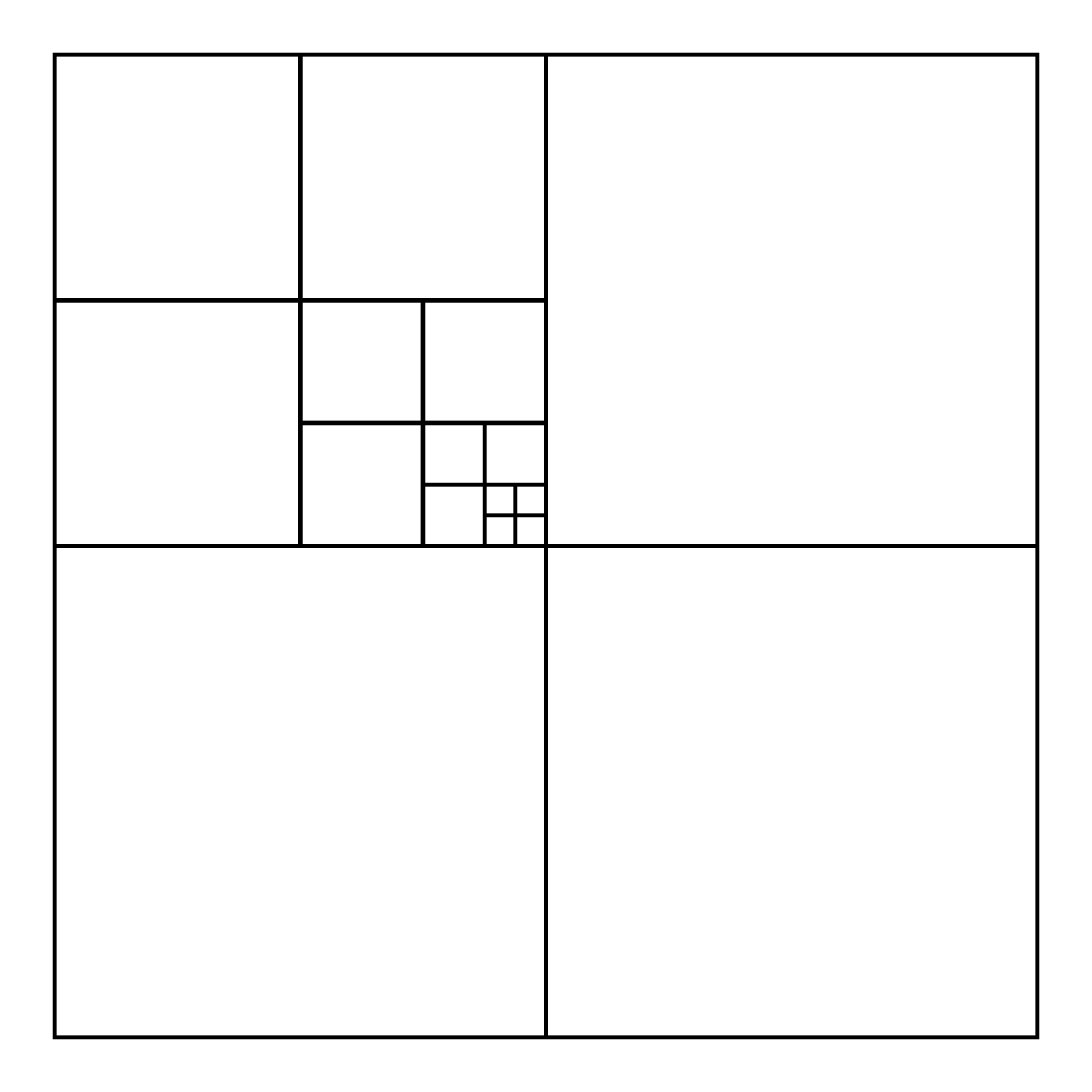} 
\caption{T-mesh in the parametric space (local refinement).}\label{local_ref}
\end{subfigure}
\caption{Refinements in the parametric space.}\label{fig:ref}
\end{figure}

T-splines, introduced by Sederberg et al. in \cite{SederbergZBN-03-TSplines, SederbergCFNZL-04-TSplines}, are defined over so-called T-meshes, see Fig.~\ref{local_ref}. T-splines are piecewise rational functions, which preserve the exactness of NURBS geometry. Moreover, a multi-patch NURBS geometry parameterization can be converted into a single patch T-splines description without gaps and leaks. However, in the most general case, the linear independence of T-splines blending functions is not guaranteed, and this led to the introduction of so-called analysis suitable T-splines \cite{beirao2013analysis,Bressan201517}. The analysis suitable T-splines have been successfully applied to IGA, see e.g. \cite{da2011isogeometric}. In a recent paper \cite{ZhangLi-16-TSplines}, the degree elevation of T-splines have also been studied. However, as described in \cite{scott2012local}, the implementation, in particular the algorithm of knot insertion which preserves analysis-suitability, is complicated.
%
\subsubsection{PHT-splines}
PHT-splines were introduced by Deng et al. in 2006 \cite{Deng200876}. In addition to main properties of B-splines, the main advantage of PHT-splines, which makes them attractive for IGA, is efficient and simple refinement algorithm \cite{NguyenThanh20111892}. However, the trade off is the reduced continuity ($C^{1}$). Nevertheless, this is enough for most applications in solid and structural mechanics. Since PHT-splines are polynomials, the CAD geometry of arbitrary topology may not be represented exactly, which led to the development of rational PHT-splines \cite{Wang20111438}.
\subsubsection{Geometric design and IGA considerations for unstructured quadrilateral meshes}
Recently, Toshniwal et al. \cite{ToshniwalSH-17} proposed a new framework for geometric design and IGA on unstructured quadrilateral meshes $M$. Following the construction of D-patch framework developed in \cite{Reif-97}, and the construction of $C^{1}$ smooth spline elements in the neighborhood of extra-ordinary points (albeit only for PHT-splines) in \cite{NguyenP-16}, the focus of \cite{ToshniwalSH-17} is on the construction of smooth, linearly independent, locally supported spline functions over $M$. To achieve the smoothness at the extra-ordinary points, affine-invariant linear transformations (called \emph{smoothing matrices}) from \cite{Reif-97} are employed. As oppose to \cite{NguyenP-16}, the work in \cite{ToshniwalSH-17} uses the smoothing matrices with non-negative coefficients. Although it yields higher quality surfaces for geometric modeling applications, this comes at the cost of slight (very small) variations in geometries during refinement. The design and analysis spline spaces are then simply defined as the span of these splines. The design space $S_{D}$ is contained within the analysis space $S_{A}$ at $k = 0$ (the coarsest level for analysis), and under the assumption of idempotence of the smoothing matrices, the spaces $S_{A}^{k}$ are nested, i.e. $S_{D} \subseteq S_{A} =: S_{A}^{0} \subseteq S_{A}^{1} \subseteq S_{A}^{2} \ldots \subseteq S_{A}^{k}$. For isogeometric compatibility, simple transformation rules for change of basis from $S_{D}$ to $S_{A}$ are constructed.

\subsection{Generalizing the isogeometric concept}
\begin{figure}[!ht]
\begin{center}
\includegraphics[width=\textwidth]{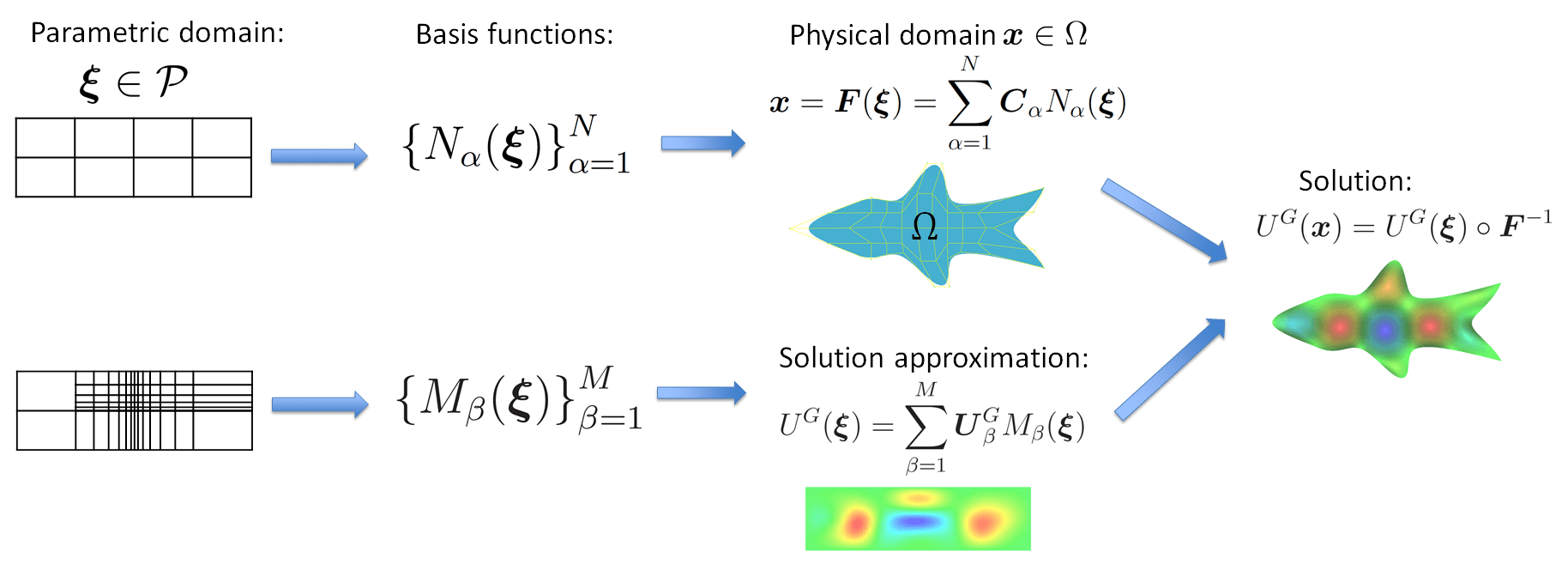}
\caption{Main idea of \GIFT \!: different shape functions are used for geometry parameterization and solution approximation.}\label{main_idea_GIFT}
\end{center}
\end{figure}

In this paper, instead of aiming to construct tailored spline spaces for dealing with local refinement and/or unstructured meshes, we aim to generalize the isogeometric concept, namely, to possibly choose different spaces for design and analysis (where $S_{D}$ may not be contained in $S_{A}$). Thereby, while remaining in the realm of existing technologies for design, we can use suitable spline spaces for analysis.
In order to preserve  the isogeometric concept, with the three choices of splines described above, it is required to convert the original NURBS CAD geometry to a parametrization (exact or approximate) in the same spline space. This contradicts the original objective of the IGA to bring the direct link between the CAD design and analysis. This limitation, caused by tight integration of the geometry parameterization and the approximation of the solution (Fig.~\ref{main_idea_IGA}) in IGA, motivated the authors to formulate a  generalized approach: {\bf Geometry Independent Field approximaTion (\GIFT)}. The main idea of \GIFT is to retain the original CAD geometry whilst flexibly adapting the solution basis to best fit the solution field (Fig.~\ref{main_idea_GIFT}). For example, a NURBS geometry can be used together with a PHT-splines approximation for the solution. The main features of \GIFT can be summarized as follows:
\begin{itemize}
\item Preserve exact CAD geometry (provided in any form, including B-splines or NURBS) at any stage of the solution process;
\item Allow local refinement of the solution by choosing appropriate field approximations, as independently as possible of the geometrical parameterization of the domain, including partition of unity enrichment, and
\item Exploit computational savings by not refining the geometry during the solution refinement process, and by choosing simpler approximations for the solution, i.e. polynomial functions instead of rational functions.
\end{itemize}
This paper aims at showing a proof of concept for \GIFT. We present the general framework of the method, followed by a number of numerical examples. In our implementation of \GIFT, the geometry is given by NURBS because this is the most commonly used form of geometry description in CAD. For solution approximations, we use a variety of approximation schemes, including NURBS functions (which have potentially different degrees and knot vectors compared to those used for the geometry), PHT-splines and B-splines.
\paragraph{Remark} In the engineering community, such different choices of spline spaces for geometry and field give rise to the question of verification of patch-test. The role of the patch test in the convergence of the finite-element based method is a long standing debate. It has been shown in \cite{Stummel-80}, that for the method to converge with optimal rate, the patch test is neither sufficient nor necessary. The present study also investigates the relation of the convergence properties of the method with the standard patch test. The spline spaces for the geometry and the field are chosen such that they may, or may not, pass the patch test.

The organization of the remainder of this article is as follows. We first present a simple introduction to the mathematical formulation of \GIFT \! in Section~\ref{sec:maths}. In Section~\ref{patch_tests}, we present the study of various patch tests. The numerical examples on convergence studies are presented in Section~\ref{num_examples}. These patch tests and numerical examples are carefully designed to study various combinations of bases for geometry and solution approximation. Finally, some conclusions and recommendations for future work are outlined in Section~\ref{conclusions}.
\section{Formulation of \GIFT}\label{sec:maths}
We consider an open domain $\Omega \subset \mathbb{R}^{d}$, $d \geq 2$, with boundary $\Gamma$ consisting of two parts $\Gamma_D$ and $\Gamma_N$, such that: $\Gamma = \overline{\Gamma_D \bigcup \Gamma_N}$, $\Gamma_D \bigcap \Gamma_N = \emptyset$. The domain $\Omega$ is parameterized on a parametric domain $\vP$ by mapping $F$:
\begin{equation}
F: \vP \rightarrow \Omega, \qquad \vx = F (\vxi), \qquad \vx \in \Omega, \quad \vxi \in \vP.
\label{domain_param_1}
\end{equation}
In what follows, we will denote the NURBS basis by $\mathcal{N}_{i,j}$, and the B-splines basis by $\mathcal{B}_{i,j}$, with $i$ denoting the degree in the $\xi$ direction, and $j$ denoting the degree in the $\eta$ direction.
Typically, the geometrical map $F$ is given by a set of basis functions $N_{i_1, i_2, \ldots, i_d}(\vxi)$ and a set of control points $\vC_{i_1, i_2, \ldots, i_d}$ as
\begin{equation}
F (\vxi) = \sum \limits_{i_{1} = 1}^{n_{1}}\sum\limits_{i_{2} = 1}^{n_{2}} \ldots \sum\limits_{i_{d} = 1}^{n_{d}}
\vC_{i_1, i_2, \ldots, i_d} N_{i_1, i_2, \ldots, i_d}(\vxi),
\label{geo_tensor_product}
\end{equation}
where $N_{i_1, i_2, \ldots, i_d}(\vxi)$ can be a tensor product of NURBS, B-splines, T-splines, PHT-splines, etc. For brevity reasons, we introduce two sets of multi-indices $(i_1, i_2, \ldots, i_d)$ of NURBS basis functions by
\begin{subequations}
\label{eq:ind_I}
\begin{align}
\vI = \big\{(i_1, i_2, \ldots, i_d) : i_{1} \in \{1, \ldots, n_{1}\}, \ldots, i_{d} \in \{1, \ldots, n_{d}\} \big\} \\
\vJ = \big\{(i_1, i_2, \ldots, i_d) : i_{1} \in \{1, \ldots, m_{1}\}, \ldots, i_{d} \in \{1, \ldots, m_{d}\} \big\} .
\end{align}
\end{subequations}
Moreover, wherever suitable, for multi-index $(i_1, i_2, \ldots, i_d)$ we will interchangeably use the collapsed notation $\vk$.
Thence, Eq.~(\ref{domain_param_1}) and (\ref{geo_tensor_product}) are written as
\begin{equation}
\vx (\vxi)= \sum \limits_{\vk \in \vI} \vC_{\vk} N_{\vk} (\vxi).
\label{domain_param_2}
\end{equation}
In what follows, we will refer to the set $\{N_{\vk}(\vxi)\}_{\vk \in \vI}$ as the \emph{geometry basis}. For change of variables, we will also need the Jacobian matrix $J(\vxi)$ of the mapping $F$, which is given by
\begin{equation}
J_{ij}(\vxi) = \dfrac{\partial x_{i}}{\partial{\xi}_{j}}(\vxi)
= \sum \limits_{\vk \in \vI} \vC_{\vk i}
\dfrac{\partial N_{\vk}(\vxi)}{\partial \xi_{j}}.
\label{jacobian}
\end{equation}
The departure from classical isogeometric analysis consists in choosing a \emph{solution basis} $\{M_{\vk}(\vxi)\}_{\vk \in \vJ}$, which is different from the geometry basis, and looking for the solution in the form:
\begin{equation}
u (\vxi)= \sum \limits_{\vk \in \vJ} u_{\vk}M_{\vk}(\vxi),
\label{solution_discretization}
\end{equation}
where $u_{\vk}$ are the unknown control variables.
In order to evaluate derivatives of the solution basis function $M_{\vk}(\vxi)$ with respect to variables $\vx$, the standard chain rule is used, which in two dimensions read:
\begin{equation}
\arraycolsep=2.5pt \def\arraystretch{2.5}
\left(
\begin{array}{c}
\dfrac{\partial M_{\vk}(\xi, \eta)}{\partial x} \\
\dfrac{\partial M_{\vk}(\xi, \eta)}{\partial y}
\end{array}
\right)
= \left(
\begin{array}{cc}
\dfrac{\partial \xi}{\partial x} & \dfrac{\partial \eta}{\partial x} \\
\dfrac{\partial \xi}{\partial y} & \dfrac{\partial \eta}{\partial y}
\end{array}
\right)
\left(
\begin{array}{c}
\dfrac{\partial M_{\vk}(\xi, \eta)}{\partial \xi} \\
\dfrac{\partial M_{\vk}(\xi, \eta)}{\partial \eta}
\end{array}
\right).
\label{chain}
\end{equation}
In matrix notations (independent of dimensions), this can be written as follows
\begin{equation}
\nabla_{\vx} M_{\vk} (\vxi)  = J^{-T}(\vxi) \nabla_{\vxi} M_{\vk} (\vxi),
\label{eq:chain_matrix}
\end{equation}
where $J^{-T}(\vxi)$ is the transpose of the inverse of the Jacobian matrix (\ref{jacobian}).
Now let the weak form of the boundary value problem be given by
\begin{equation}
a(u, v) = \ell(v),
\label{weak_form}
\end{equation}
then substituting Eq.~(\ref{domain_param_2}) and (\ref{solution_discretization}) into (\ref{weak_form}) we obtain the linear system of equations
\begin{equation}
\vK \vu = \vf,
\end{equation}
where the stiffness matrix $\vK$ and the force vector $\vf$ are given by
\begin{align}
\vK_{ij} = a(M_{i}(\vx), M_{j}(\vx)), \qquad
\vf_{i} = \ell(M_{i}(\vx)),
\label{m_stiff_v_force}
\end{align}
and vector $\vu$ consists of all unknown control variables $u_{j}$.
For our numerical results, we consider two problems, namely, a) Poisson's equation, and b) the problem of linear elasticity. In Sections~\ref{sec:Poisson}-\ref{sec:elas}, we briefly describe the derivations of the weak form of these problems.

\subsection{Poisson's equation}
\label{sec:Poisson}
For the Poisson's problem, the terms in the weak form are defined as
\begin{equation}
a(u, v) = \int_{\Omega}\nabla u \cdot \nabla v \, d\Omega, \qquad
\ell(v) = \int_{\Omega} f~ v~ d\Omega  + \int_{\partial\Omega_{N}} g~ v~ d\Gamma ,
\label{weak_form_Poisson}
\end{equation}
where $f$ is the source function, and $g$ is the flux prescribed on Neumann part of the boundary. The Eq.~\eqref{m_stiff_v_force} thus results in, for $i, j = 1, \ldots, M$,
\begin{equation}
\begin{split}
\vK_{ij} & = \int_\Omega \nabla_{\vx} M_{i} ({\vx}) \cdot \nabla_{\vx} M_{j} ({\vx})  ~ d \Omega, \\
\vf_{i} & = \int_\Omega f(\vx) ~ M_{i} ({\vx}) ~ d \Omega
+ \int_{\partial \Omega_{N}} g(\vx) M_{i}(\vx) d\Gamma.
\end{split}
\label{eqn:Kij_x}
\end{equation}

\subsection{Linear elasticity}
\label{sec:elas}

Consider the linear elasticity problem in two dimensions, for which the weak form is given by
\begin{equation}
a(u, v) = \int_{\Omega} \veps (u)^{T} \vD \veps (v)~ d\Omega, \qquad
\ell(v) = \int_{\Omega} f~ v~ d\Omega  + \int_{\partial\Omega_{N}} \tilde{t}~ v~ d\Gamma,
\end{equation}
where $u := (u_{x}, u_{y})^{T}$ denotes the vector of unknown displacements, $f$ denotes the body forces, $\tilde{t}$ denotes the tractions prescribed on Neumann boundary, and the matrix of material parameters in plane strain is given by
\begin{equation}
\vD = \dfrac{E}{(1 + \nu)(1 - 2\nu)}\left(\begin{array}{ccc}
                       1 - \nu & \nu & 0\\
                       \nu & 1 - \nu & 0 \\
                       0 & 0 & (1 - 2\nu)/2 
                       \end{array}
                       \right),
\end{equation}
where $E$ is the Young's modulus and $\nu$ is the Poisson's coefficient.
The strain operator, which is defined as
\begin{equation}
\arraycolsep=2.5pt \def\arraystretch{1.5}
\veps(u) = \left(\begin{array}{c}
\partial u_{x}/\partial {x}\\
\partial u_{y}/\partial {y}\\
\partial u_{x}/\partial {y} + \partial u_{y}/\partial {x}
\end{array}
\right)
\label{eq:strain}
\end{equation}
applied to the solution approximation given by Eq.~\eqref{solution_discretization} yields the following stiffness matrix and force vector:
\begin{equation}
\begin{split}
\vK_{ij} & = \int_{\Omega} \vB_{i}^{T} \vD \vB_{j}~ d\Omega,\\
\vf_{i} & = \int_\Omega f({\boldsymbol{x}})  ~M_{i} ({\boldsymbol{x}}) ~ d\Omega
+ \int_{\partial\Omega_{N}} \tilde{t}(\boldsymbol{x}) M_{i}(\boldsymbol{x}) d\Gamma \quad,
\end{split}
\end{equation}
where the strain operator $\vB_{i}$ is defined as
\begin{equation}
\arraycolsep=2.5pt \def\arraystretch{2.5}
\vB_{i} = \left( \begin{array}{cc}
\dfrac{\partial M_{i}(\xi, \eta)}{\partial {x}} & 0\\
0 & \dfrac{\partial M_{i}(\xi, \eta)}{\partial {y}}\\
\dfrac{\partial M_{i}(\xi, \eta)}{\partial {y}} & \dfrac{\partial M_{i}(\xi, \eta)}{\partial {x}}\\
\end{array}\right),
\end{equation}
and derivatives of the solution shape functions with respect to geometry variables are evaluated according to Eq.~(\ref{chain}).
\section{Patch tests}\label{patch_tests}
In this Section, we will study the classical patch test for various combinations of bases for the geometry representation and solution approximation. For all the test cases, we will study two problems, namely, Laplace equation, and linear elasticity problem, which are defined in Section \ref{sec:pt_prob_def}. We present "engineering" analysis of the patch test results (including the values of the field control variables). This can also serve as a guide for engineers to design patch tests in spline based methods.

\subsection{Problem definition}\label{sec:pt_prob_def}
\begin{figure}[!ht]
\begin{center}
\includegraphics[width=0.65\textwidth]{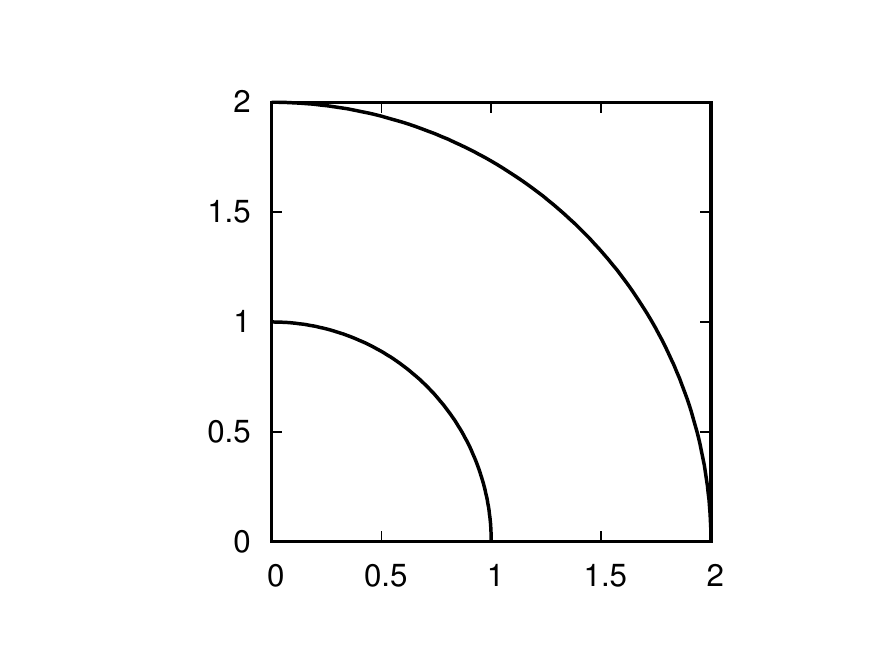}
\caption{Quarter annulus.}
\label{fig:qannu}
\end{center}
\end{figure}
We formulate two problems in a quarter annulus domain $\Omega$, as shown in Fig.~\ref{fig:qannu}. The first problem is for the two-dimensional Laplace equation, where we use the linear solution 
\begin{equation}
u(x,y) = 1 + x + y,
\end{equation}
imposed as the Dirichlet boundary condition on the entire boundary $\partial\Omega$, i.e.
\begin{equation}\label{PT_laplace}
\begin{split}
\Delta{u} &= 0,  \,\,\,\,\text{in}\,\,\,\Omega \\
u|_{\partial\Omega} (x,y)&= 1 + x + y. 
\end{split}
\end{equation}
The second problem is for the linear elasticity, which is given by:
\begin{equation}
\begin{split}
\sigma_{ij,j} & = 0,  \qquad \text{in} \quad \Omega \\
t_{i} &= \sigma_{0}n_{i}, \qquad \text{at} \quad r = 1, 2, \\
u_{2} &= 0, \quad t_{1} = 0, \qquad \text{at} \quad \theta = 0, \\
u_{1} &= 0, \quad t_{2} = 0, \qquad \text{at} \quad \theta = \pi/2, \\
\end{split}
\end{equation}
which admits a linear solution (plane-strain case) given by
\begin{equation}\label{PT_elasticity}
u_{1} (x)= \dfrac{(1 + \nu)(1 - 2\nu)}{E}\sigma_{0}x, \,\,\,\, u_{2} (y)= \dfrac{(1 + \nu)(1 - 2\nu)}{E}\sigma_{0}y.
\end{equation}
\subsection{Geometry parametrization}\label{sec:geo_param}
\begin{figure}[!ht]
\begin{center}
\includegraphics[width=0.65\textwidth]{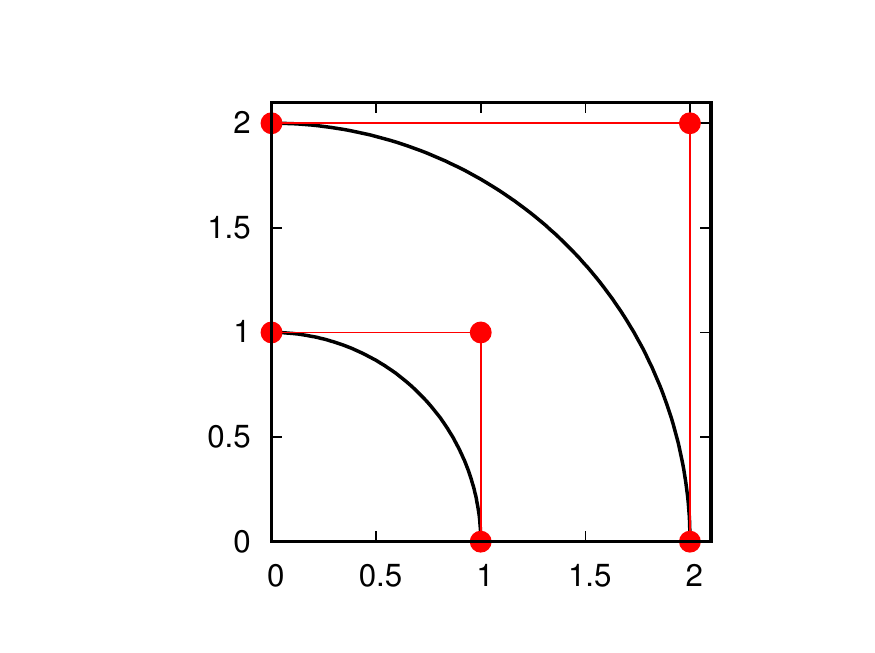}
\caption{Coarsest parametrization of a quarter annulus.}
\label{fig:Q0}
\end{center}
\end{figure}
As shown in Fig.~\ref{fig:Q0}, the coarsest parameterization of the quarter annulus can be given by NURBS of degrees $p_{\xi} = 1$ and $p_{\eta} = 2$, built upon the two knot vectors:
\begin{equation}
\Sigma_{0} = \{0,0,1,1\}\hspace{1cm}\textrm{and}\hspace{1cm}\Pi_{0} = \{0,0,0,1,1,1\},
\label{coarsest_quarter_annulus_knots}
\end{equation}
and the 6 control points $P_{ij}$, $i = 1, 2$, and $j = 1, 2, 3$, given in Table~\ref{tab:param_Q0}.
\begin{table}
\begin{center}
\begin{tabular}{|c|c|c|c|}
\hline 
$(i,j)$ & $P^{Q_{0}, x}_{ij}$ & $P^{Q_{0}, y}_{ij}$ & $w^{Q_{0}}_{ij}$ \\
\hline
(1,1) & 1.0000 & 0.0000 & 1.0000 \\
(1,2) & 1.0000 & 1.0000 & 0.7071 \\
(1,3) & 0.0000 & 1.0000 & 1.0000 \\
(2,1) & 2.0000 & 0.0000 & 1.0000 \\
(2,2) & 2.0000 & 2.0000 & 0.7071 \\
(2,3) & 0.0000 & 2.0000 & 1.0000 \\
\hline 
\end{tabular} 
\caption{Control points and weights for coarsest parameterization $Q_{0}$}
\label{tab:param_Q0}
\end{center}
\end{table}
In what follows we will refer to this parametrization as $Q_{0}$.
In this section, together with $Q_{0}$, we will consider four 4-element parameterizations of the quarter annulus. Three will be called "uniform", denoted by letters $A$, $B$, and $D$, and one will be called "non-uniform", denoted by letter $C$. Note that the parameterization $D$ is obtained from the parameterization $A$ by setting all the weights to $1$. It does not exactly represent the quarter annulus, and it will only be used for solution approximation. The "uniform" parameterizations are obtained by performing the standard operation of knot-insertion on knot vectors $\Sigma_{0}$ and $\Pi_{0}$ to calculate the positions $P_{ij}$ and the weights $w_{ij}$ of the corresponding control points. We add knot value $s$ in $\xi$-direction and knot value $t$ in $\eta$-direction. Then the geometry is given by:
\begin{equation}
\Sigma_{1} = \{0,0,s, 1,1\}\hspace{1cm}\textrm{and}\hspace{1cm}\Pi_{1} = \{0,0,0,t, 1,1,1\}.
\label{knots1}
\end{equation}
The "non-uniform" parametrization is obtained by adding two knot values $s$ and $t$ to the knot vectors $\Sigma_{0}$ and $\Pi_{0}$, analogously to Eq.~\eqref{knots1}, but arbitrarily choosing the weights and the position of two inner control points: $P_{2,2}$ and $P_{2,3}$.
\begin{table}
\begin{center}
\begin{tabular}{|c|c|c|c|c|c|}
\hline 
 & $p_{\xi}$ & $p_{\eta}$ & $ s $ & $ t $ & $P_{ij}$ \& $w_{ij}$ \\
\hline 
$A_{1}$ & 1 & 2 & 2/3 & 1/8 & Table \ref{tab:param_A1} \\
$B_{1}$ & 1 & 2 & 0.17 & 0.81 & Table \ref{tab:param_B1} \\
$C_{1}$ & 1 & 2 & 2/3 & 1/8 & Table \ref{tab:param_C1} \\
$D_{1}$ & 1 & 2 & 2/3 & 1/8 & Table \ref{tab:param_D1} \\
\hline 
\end{tabular} 
\caption{Various parameterizations}
\label{tab:param_all}
\end{center}
\end{table}
The control points $P_{ij}$ and the corresponding weights $w_{ij}$ for parameterizations $A_{1}$, $B_{1}$, $C_{1}$, and $D_{1}$ are listed in Tables~\ref{tab:param_A1}, \ref{tab:param_B1}, \ref{tab:param_C1}, and \ref{tab:param_D1}, respectively. The higher order parameterizations $A_{2}$, $B_{2}$, $C_{2}$, and $D_{2}$ are obtained by standard degree elevation (by one in both directions) of parameterizations $A_{1}$, $B_{1}$, $C_{1}$, and $D_{1}$, respectively. Thereafter, the degrees of the parameterizations $A_{2}$, $B_{2}$, $C_{2}$, and $D_{2}$ will be given by $p_{\xi} = 2$ and $p_{\eta} = 3$. Moreover, for a special case, we will also consider a parameterization $D_{0}$, where $p_{\xi} = 1$ and $p_{\eta} = 1$, and which is obtained by degree reduction of parameterization $D_{1}$.
\begin{table}
\begin{center}
\begin{tabular}{|c|c|c|c|}
\hline 
$(i,j)$ & $P^{A_{1}, x}_{ij}$ & $P^{A_{1}, y}_{ij}$ & $w^{A_{1}}_{ij}$ \\
\hline 
(1,1) & 1.0000 & 0.0000 & 1.0000 \\
(1,2) & 1.0000 & 0.0917 & 0.9634 \\
(1,3) & 0.8320 & 1.0000 & 0.7437 \\
(1,4) & 0.0000 & 1.0000 & 1.0000 \\
(2,1) & 1.6667 & 0.0000 & 1.0000 \\
(2,2) & 1.6667 & 0.1529 & 0.9634 \\
(2,3) & 1.3865 & 1.6667 & 0.7437 \\
(2,4) & 0.0000 & 1.6667 & 1.0000 \\
(3,1) & 2.0000 & 0.0000 & 1.0000 \\
(3,2) & 2.0000 & 0.1835 & 0.9634 \\
(3,3) & 1.6639 & 2.0000 & 0.7437 \\
(3,4) & 0.0000 & 2.0000 & 1.0000 \\
\hline
\end{tabular} 
\caption{Control points and weights for parameterization $A_{1}$}
\label{tab:param_A1}
\end{center}
\end{table}
\begin{table}
\begin{center}
\begin{tabular}{|c|c|c|c|}
\hline
$(i,j)$ & $P^{B_{1}, x}_{ij}$ & $P^{B_{1}, y}_{ij}$ & $w^{B_{1}}_{ij}$ \\
\hline
(1,1) & 1.0000 & 0.0000 & 1.0000 \\
(1,2) & 1.0000 & 0.7509 & 0.7628 \\
(1,3) & 0.1423 & 1.0000 & 0.9444 \\
(1,4) & 0.0000 & 1.0000 & 1.0000 \\
(2,1) & 1.1700 & 0.0000 & 1.0000 \\
(2,2) & 1.1700 & 0.8786 & 0.7628 \\
(2,3) & 0.1665 & 1.1700 & 0.9444 \\
(2,4) & 0.0000 & 1.1700 & 1.0000 \\
(3,1) & 2.0000 & 0.0000 & 1.0000 \\
(3,2) & 2.0000 & 1.5018 & 0.7628 \\
(3,3) & 0.2845 & 2.0000 & 0.9444 \\
(3,4) & 0.0000 & 2.0000 & 1.0000 \\
\hline
\end{tabular} 
\caption{Control points and weights for parameterization $B_{1}$}
\label{tab:param_B1}
\end{center}
\end{table}
\begin{table}
\begin{center}
\begin{tabular}{|c|c|c|c|}
\hline 
$(i,j)$ & $P^{C_{1}, x}_{ij}$ & $P^{C_{1}, y}_{ij}$ & $w^{C_{1}}_{ij}$ \\
\hline 
(1,1) & 1.0000 & 0.0000 & 1.0000 \\
(1,2) & 1.0000 & 0.0917 & 0.9634 \\
(1,3) & 0.8320 & 1.0000 & 0.7437 \\
(1,4) & 0.0000 & 1.0000 & 1.0000 \\
(2,1) & 1.6667 & 0.0000 & 1.0000 \\
\bf{(2,2)} & \bf{1.1000} & \bf{0.3000} & \bf{0.8000} \\
\bf{(2,3)} & \bf{0.7500} & \bf{1.4000} & \bf{0.7500} \\
(2,4) & 0.0000 & 1.6667 & 1.0000 \\
(3,1) & 2.0000 & 0.0000 & 1.0000 \\
(3,2) & 2.0000 & 0.1835 & 0.9634 \\
(3,3) & 1.6639 & 2.0000 & 0.7437 \\
(3,4) & 0.0000 & 2.0000 & 1.0000 \\
\hline
\end{tabular} 
\caption{Control points and weights for parameterization $C_{1}$}
\label{tab:param_C1}
\end{center}
\end{table}
\begin{table}
\begin{center}
\begin{tabular}{|c|c|c|c|}
\hline 
$(i,j)$ & $P^{D_{1}, x}_{ij}$ & $P^{D_{1}, y}_{ij}$ & $w^{D_{1}}_{ij}$ \\
\hline 
(1,1) & 1.0000 & 0.0000 & 1.0000 \\
(1,2) & 1.0000 & 0.0917 & 1.0000 \\
(1,3) & 0.8320 & 1.0000 & 1.0000 \\
(1,4) & 0.0000 & 1.0000 & 1.0000 \\
(2,1) & 1.6667 & 0.0000 & 1.0000 \\
(2,2) & 1.6667 & 0.1529 & 1.0000 \\
(2,3) & 1.3865 & 1.6667 & 1.0000 \\
(2,4) & 0.0000 & 1.6667 & 1.0000 \\
(3,1) & 2.0000 & 0.0000 & 1.0000 \\
(3,2) & 2.0000 & 0.1835 & 1.0000 \\
(3,3) & 1.6639 & 2.0000 & 1.0000 \\
(3,4) & 0.0000 & 2.0000 & 1.0000 \\
\hline
\end{tabular} 
\caption{Control points and weights for parameterization $D_{1}$}
\label{tab:param_D1}
\end{center}
\end{table}
\begin{figure}[!ht]
\begin{center}
\includegraphics[width=0.65\textwidth]{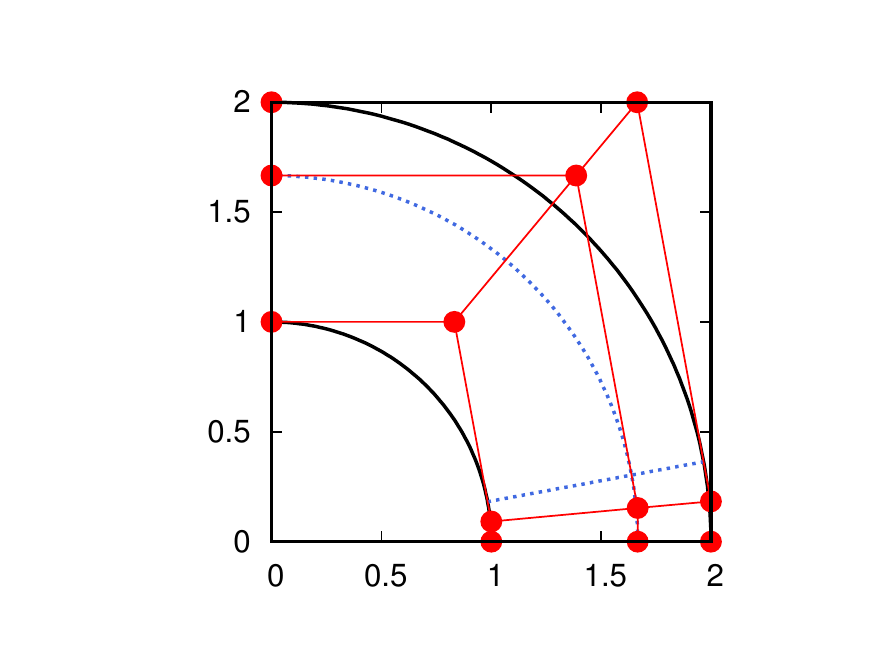}
\caption{Uniform parametrization $A_{1}$: elements and control points.}
\label{fig:A1}
\end{center}
\end{figure}
\begin{figure}[!ht]
\begin{center}
\includegraphics[width=0.65\textwidth]{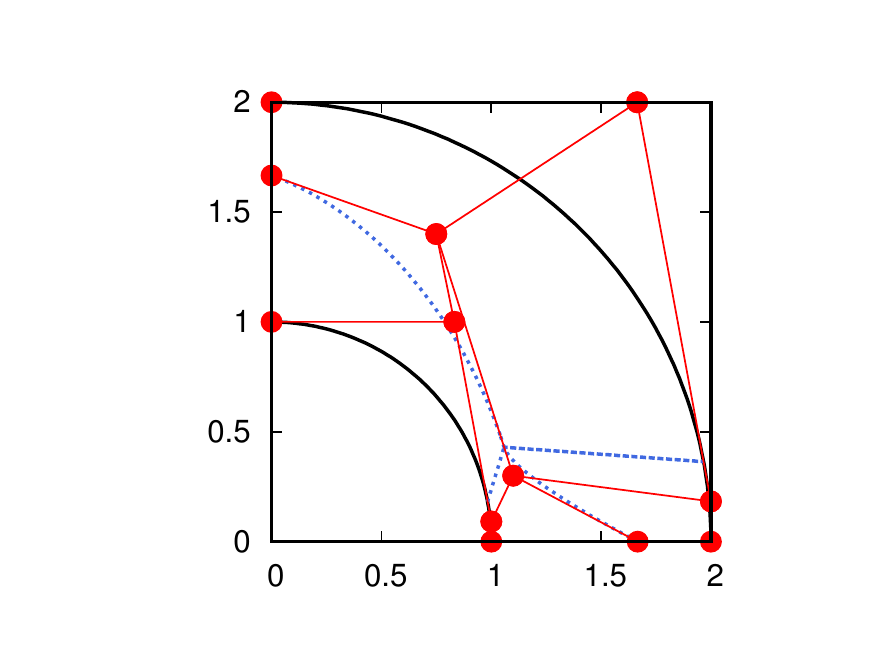}
\caption{Non-uniform parametrization $C_{1}$: elements and control points.}
\label{fig:C1}
\end{center}
\end{figure}
As discussed in the following section, a total of 19 patch tests for each of the problems of Section~\ref{sec:pt_prob_def} (in two-dimensions) are studied. Therefore, for better readability and classification, we introduce the following notation. We will denote a patch test by $T^{\ell}_{G_{i}, S_{j}}$, where the superscript $\ell = 0, 1, \ldots, 5$, denotes the patch test case, $G_{i}$ denotes the basis for geometry parameterization, and $S_{j}$ denotes the basis for solution approximation. The basis $G_{i}$ and $S_{j}$ are chosen from the bases set given in Table~\ref{tab:param_all}, and their degree elevated/reduced versions.

\subsection{Results and discussion of patch tests}
\begin{table}[!ht]
\renewcommand*{\arraystretch}{1.3}
\begin{center}
\begin{tabular}{|c|c|c|}
\hline 
$T^{\ell}_{G_{i}, S_{j}}$ & Laplace Eq.~\eqref{PT_laplace} & Elasticity Eq.~\eqref{PT_elasticity} \\
\hline 
$T^{0}_{Q_{0}, A_{1}}$ & 1.3815e-15 & 3.0871e-14 \\
$T^{0}_{Q_{0}, A_{2}}$ & 5.2147e-15 & 1.7986e-14 \\
$T^{0}_{Q_{0}, C_{1}}$ & 0.0182 & 0.0050 \\
$T^{0}_{Q_{0}, C_{2}}$ & 0.0023 & 0.0012 \\
\hline \hline
$T^{1}_{A_{1}, A_{1}}$ & 1.0023e-15 & 1.1675e-14 \\
$T^{1}_{A_{1}, A_{2}}$ & 4.3958e-14 & 1.2547e-14 \\
$T^{1}_{A_{2}, A_{1}}$ & 1.4059e-15 & 1.5525e-15 \\
\hline \hline
$T^{2}_{B_{1}, A_{1}}$ & 1.4755e-15 & 2.9941e-15 \\
$T^{2}_{B_{1}, A_{2}}$ & 2.1639e-15 & 1.2118e-14 \\
$T^{2}_{B_{2}, A_{1}}$ & 1.0144e-15 & 5.4590e-15 \\
\hline \hline
$T^{3}_{C_{1}, C_{1}}$ & 1.1061e-15 & 1.6439e-14 \\
$T^{3}_{C_{1}, C_{2}}$ & 1.8263e-15 & 2.8737e-15 \\
$T^{3}_{C_{2}, C_{1}}$ & 1.2062e-15 & 5.6517e-14 \\
\hline \hline
$T^{4}_{C_{1}, A_{1}}$ & 0.0203 & 0.0085 \\
$T^{4}_{C_{1}, A_{2}}$ & 0.0016 & 0.0009 \\
$T^{4}_{C_{2}, A_{1}}$ & 0.0203 & 0.0085 \\
\hline \hline
$T^{5}_{A_{1}, D_{1}}$ & 0.0188 & 0.0214 \\
$T^{5}_{A_{1}, D_{2}}$ & 0.0121 & 0.0039 \\
$T^{5}_{A_{1}, D_{0}}$ & 0.5418 & 0.1411 \\
\hline
\end{tabular} 
\caption{Results of various patch tests}
\label{tab:res_tests}
\end{center}
\end{table}
We start the series of patch tests by choosing the coarsest parametrization of the quarter annulus, i.e. $Q_{0}$, paired with different bases for the solution.
The results of the first patch test are denoted by $T^{0}$ in Table~\ref{tab:res_tests}. Note that all the combinations of basis functions fall within the isoparametric or sub-parametric finite element concept. However, the first two combinations pass the test, while the other two fail.
Therefore, in order to understand the connection between the geometry and the solution bases, that leads to either passing or failing the patch test, we introduce the following examples, investigating other choices of the geometry parameterization and solution approximation.

In the next patch test, denoted by $T^{1}$ in Table~\ref{tab:res_tests}, we demonstrate the performance of uniform parameterizations $A_{1}$ and $A_{2}$.
All  three combinations pass the patch test. Note that the test $T^{1}_{A_{1}, A_{1}}$ is a standard IGA patch test within isoparametric concept, which can exactly represent a constant/linear solution on any basis which preserves partition of unity.
To understand the relation between geometry representation and solution approximation, we substitute the geometry parametrization
\begin{equation}\label{PT-ISO-GEO}
\def\arraystretch{1.5}
\begin{split}
\begin{pmatrix}
           x \\
           y 
         \end{pmatrix}
         &= \sum\limits_{ij}\begin{pmatrix}
           P^{A_{1}, x}_{ij} \\
           P^{A_{1}, y}_{ij}
         \end{pmatrix}\times\{\text{Basis of~} A_{1}\} 
\end{split}         
\end{equation}
and the solution approximation
\begin{equation}\label{PT-ISO-GEO-u}
u = \sum \limits_{ij} u_{ij} \times \{\text{Basis of~} A_{1}\}
\end{equation}
into the solution
\begin{equation}\label{PT-sol}
u = 1 + x + y,
\end{equation}
and obtain the solution control variables as
\begin{equation}\label{PT-iso-control-var}
u_{ij} = 1 + P^{A_{1}, x}_{ij} +  P^{A_{1}, y}_{ij}.
\end{equation}
Now recall the following fundamental property of the degree elevation algorithm \cite[Alg.~5.9, P.206]{NURBSbook} (which preserves the exactness of the geometry representation)
\begin{equation}\label{degree_elevation_geo_exactness}
\def\arraystretch{1.5}
\begin{split}
\begin{pmatrix}
           x \\
           y 
         \end{pmatrix}
         &= \sum\limits_{ij}\begin{pmatrix}
           P^{A_{1}, x}_{ij} \\
           P^{A_{1}, y}_{ij}
         \end{pmatrix}\times\{\text{Basis of~} A_{1}\} \\
          &= \sum\limits_{ij}\begin{pmatrix}
           P^{A_{2}, x}_{ij} \\
           P^{A_{2}, y}_{ij}
         \end{pmatrix}\times\{\text{Basis of~} A_{2}\}.
\end{split}
\end{equation}
The test $T^{1}_{A_{1}, A_{2}}$ can be viewed analogous to the FEA sub-parametric approach (the degree of the solution approximation is higher than the degree of the geometry). If the solution is approximated by
\begin{equation}\label{PT-SUPER-sol}
u = \sum \limits_{ij} u_{ij} \times \{\text{Basis of~} A_{2}\}
\end{equation}
then, using \eqref{degree_elevation_geo_exactness}, the corresponding control variables can be found as
\begin{equation}\label{PT-super-control-var}
u_{ij} = 1 + P^{A_{2}, x}_{ij} +  P^{A_{2}, y}_{ij}.
\end{equation}
The test $T^{1}_{A_{2}, A_{1}}$ can be viewed analogous to the FEA super-parametric approach (the degree of the solution approximation is lower than the degree of the geometry). However, the essential difference is that the geometry degree was artificially lifted from $A_{1}$ to $A_{2}$. Since using \eqref{degree_elevation_geo_exactness} in \eqref{PT-sol} yields the solution given by \eqref{PT-iso-control-var}, this means that both the bases $A_{1}$ and $A_{2}$, as well as any two parameterizations with the same property, yield identical solutions for the solution control variables (within the numerical tolerance).

For the second patch test, denoted by $T^{2}$ in Table~\ref{tab:res_tests}, we consider the situation when we have two uniform, but different, parameterizations of the domain, and one of them is used to parametrize the geometry and the second one is used as a basis to approximate the unknown solution.
All the three combinations pass the test. Note that in this test, we used the same basis for the solution approximation as in the test $1$, but a different geometry parametrization $B_{1}$ (instead of $A_{1}$). However, both parameterizations, $A_{1}$ and $B_{1}$, were obtained from the coarsest parametrization $Q_{0}$ of the quarter annulus.
According to the property of the knot insertion algorithm  \cite[Alg.~5.3, P.155]{NURBSbook}, this gives
\begin{equation}\label{knot_insertion_geo_exactness}
\def\arraystretch{1.5}
\begin{split}
\begin{pmatrix}
           x \\
           y 
         \end{pmatrix}
         &= \sum\limits_{ij}\begin{pmatrix}
           P^{A_{1}, x}_{ij} \\
           P^{A_{1}, y}_{ij}
         \end{pmatrix}\times\{\text{Basis of~} A_{1}\} \\
          &= \sum\limits_{ij}\begin{pmatrix}
           P^{B_{1}, x}_{ij} \\
           P^{B_{1}, y}_{ij}
         \end{pmatrix}\times\{\text{Basis of~} B_{1}\}.
\end{split}         
\end{equation}
Due to the property (\ref{knot_insertion_geo_exactness}), geometry parameterizations $A_{1}$ and $B_{1}$ (as well as their higher degree versions $A_{2}$ and $B_{2}$, or any other two parameterizations with the same property) analytically yield the results for the solution control variables in $T^{2}$ to be identical (within the numerical tolerance) to the solution control variables in $T^{1}$, which is given by Eq.~\eqref{PT-iso-control-var} (or for higher degree by Eq.~\eqref{PT-super-control-var}.
Note that the same property holds between bases $Q_{0}$ and $A_{1}$ in the test $T^{0}$, which implies that the solution in $T^{0}_{Q_{0}, A_{1}}$ is expressed by Eq.~\eqref{PT-iso-control-var}, and by adding the property of degree elevation Eq.~\eqref{degree_elevation_geo_exactness}, we can conclude that the solution in $T^{0}_{Q_{0}, A_{2}}$ is given by Eq.~\eqref{PT-super-control-var}.

In the third patch test, denoted by $T^{3}$ in Table~\ref{tab:res_tests}, we investigate the performance of non-uniform parametrization. The results of the third patch test show that all the three combinations pass the test. In fact, in this test only degree elevation algorithm (from $C_{1}$ to $C_{2}$) is used, which makes it fully analogous to $T^{1}$, albeit with a non-uniform parameterization.

In tests $T^{1}$, $T^{2}$ and $T^{3}$, the bases of solution and geometry were related, either both the bases were obtained from the geometry parameterization $Q^{0}$ (as in tests $T^{1}$ and $T^{2}$), or the solution basis was obtained from the geometry basis (as in test $T^{3}$).
The idea of the fourth patch test, denoted by $T^{4}$ in Table~\ref{tab:res_tests}, is to show the combinations of the geometry parameterizations and solution approximations which will fail to produce the linear solution. For this purpose we use the same approximation bases $A_{1}$ and $A_{2}$ for the solution, as in $T^{1}$ and $T^{2}$, but combined with the non-uniform geometry parameterization $C_{1}$ and $C_{2}$.
As it is seen from Table~\ref{tab:res_tests}, this test fails for all the combinations. Due to the fact that $C_{1}$ was not derived from the coarsest parameterization $Q_{0}$ of the quarter annulus, the property analogous to \eqref{knot_insertion_geo_exactness} between $A_{1}$ and $C_{1}$ does not hold.
In the fifth patch test, denoted by $T^{5}$ in Table~\ref{tab:res_tests}, the conditions (\ref{degree_elevation_geo_exactness}) and (\ref{knot_insertion_geo_exactness}) between the geometry parameterization and the solution approximation do not hold, and therefore as expected, all combinations of $T^{5}$ fail. However, it is important for practical applications, and the suggested bases combinations will be used later in the numerical examples. We use the NURBS basis $A_{1}$ for geometry parameterization, and the B-splines basis  $D_{i}$, $i, = \{0, 1, 2\}$ for solution approximation. Two important features of $T^{5}$ are the following:
\begin{itemize}
\item The sub-parametric case is different from what is considered in the previous patch tests, because instead of artificially lifting the degree of the geometry parametrization, we kept the minimum required order of the geometry basis, and reduced the degree of the solution approximation.
\item In all previous patch tests, the basis used for the solution approximation, could also be used to exactly represent the geometry. However, in $T^{5}$, the B-splines basis, which is used for the solution, cannot exactly represent the geometry.
\end{itemize}
{\bf Conclusion:} The {\bf sufficient condition} for the patch test to pass consists in the requirement for the geometry and field bases to be equivalent, up to operations of degree elevation and/or knot insertion.

\section{Numerical examples}\label{num_examples}
In this section we present convergence results of various bases choices for the geometry and the numerical solution. We consider two classes of problems, namely, Laplace and linear elasticity.

\subsection{Laplace equation in two-dimensions}\label{sec:laplace}
For Laplace equation, we choose two geometries, the quarter annulus (see Fig.~\ref{fig:qannu}) which can be exactly represented by $\mathcal{N}_{1, 2}$, and a wedge geometry (see Fig.~\ref{wedge_geo}) which is represented by $\mathcal{N}_{1, 4}$.

\subsubsection{Example 1}\label{example_laplace_1}
\begin{figure}[!ht]
\begin{center}
\includegraphics[width=0.9\textwidth]{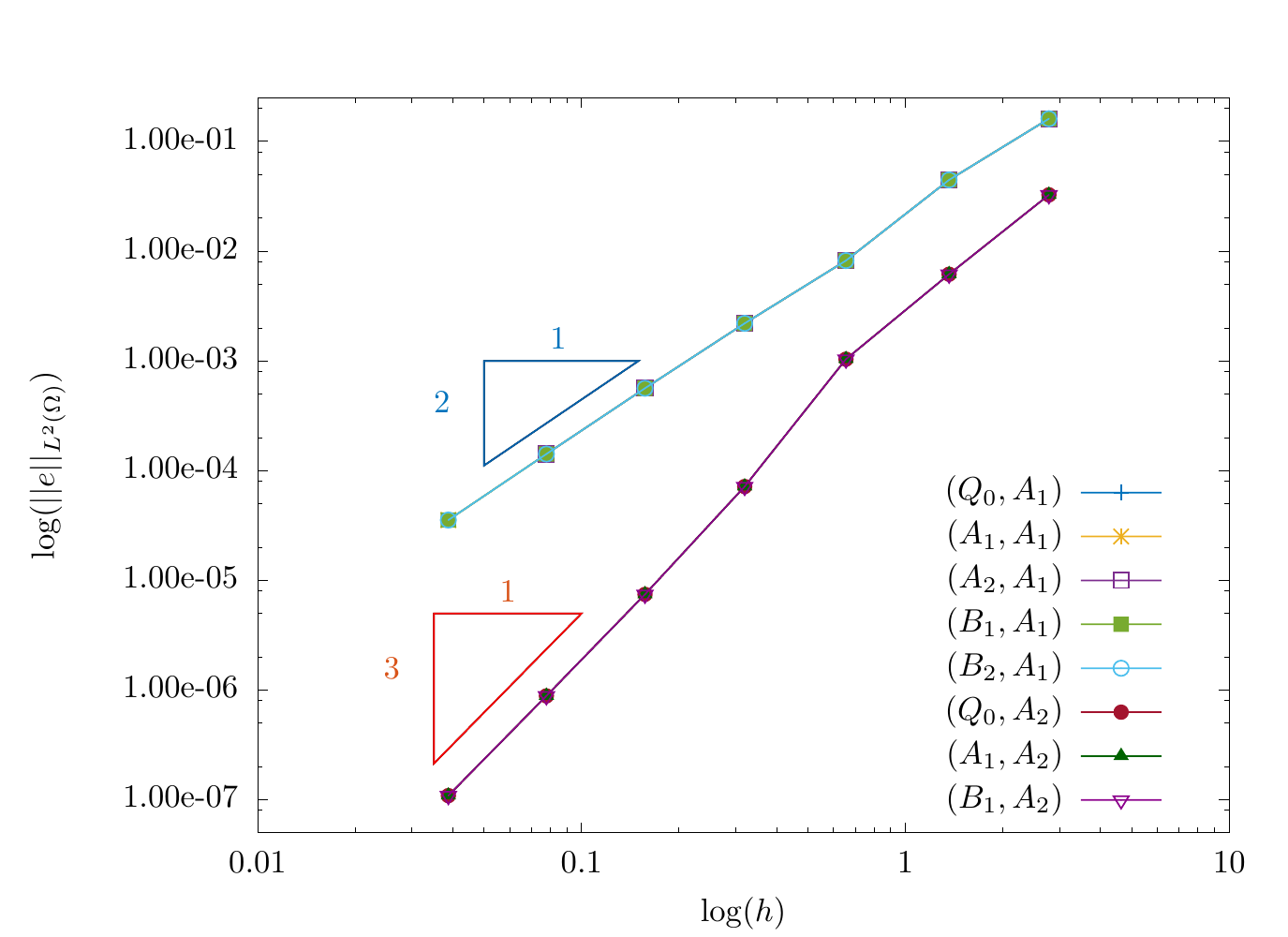}
\caption{Convergence study for Example~\ref{example_laplace_1}. Geometry parameterizations $Q_{0}$, $A_{1}$, $A_{2}$, $B_{1}$ and $B_{2}$ combined with solution bases $A_{1}$ and $A_{2}$.}
\label{example_1_1}
\end{center}
\end{figure}

\begin{figure}[!ht]
\begin{center}
\includegraphics[width=0.9\textwidth]{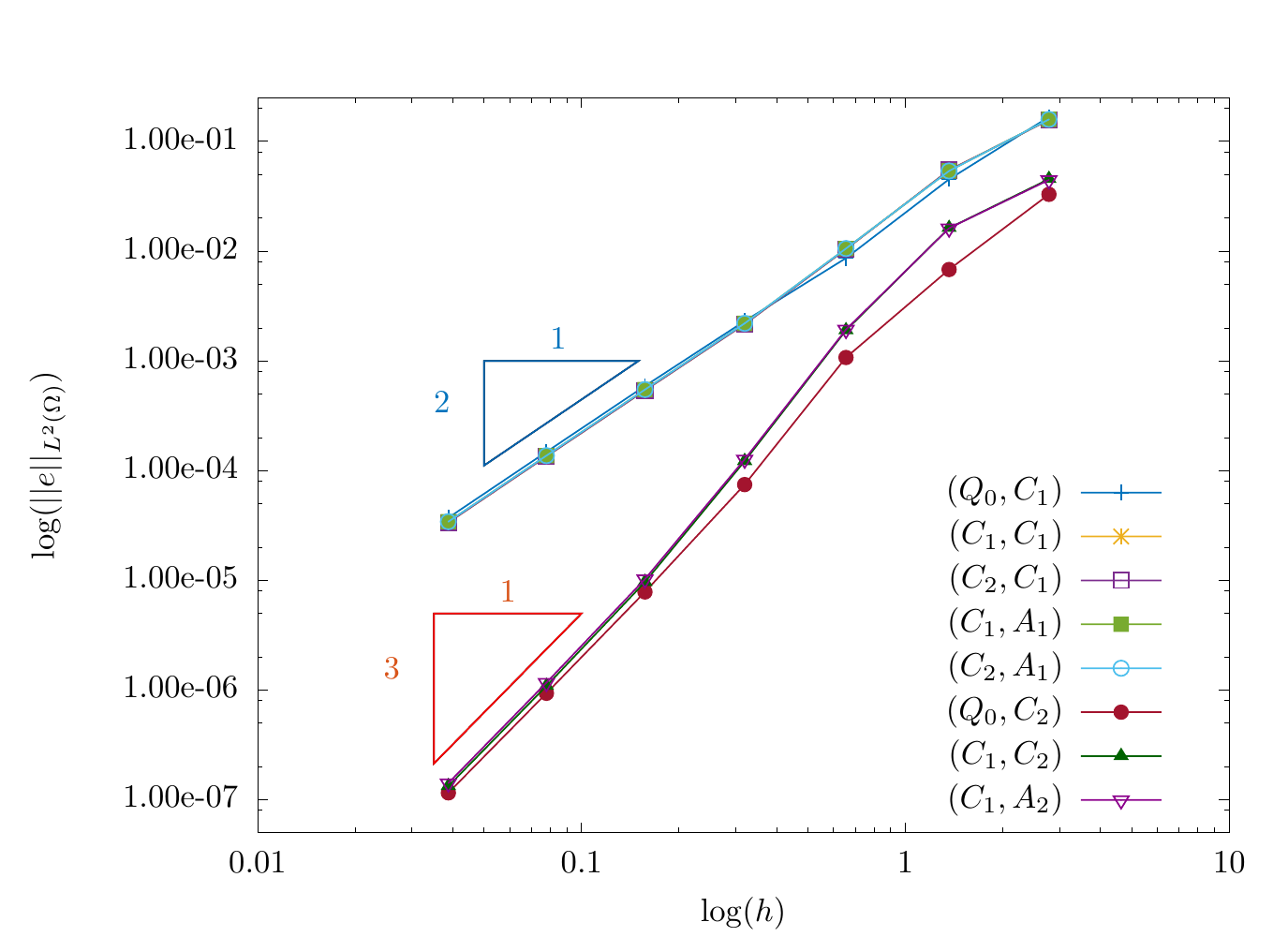}
\caption{Convergence study for Example~\ref{example_laplace_1}. Non-uniform geometry parameterizations $C_{1}$ and $C_{2}$ combined with solution bases $C_{1}$, $C_{2}$, $A_{1}$, and $A_{2}$.}
\label{example_1_2}
\end{center}
\end{figure}

\begin{figure}[!ht]
\begin{center}
\includegraphics[width=0.9\textwidth]{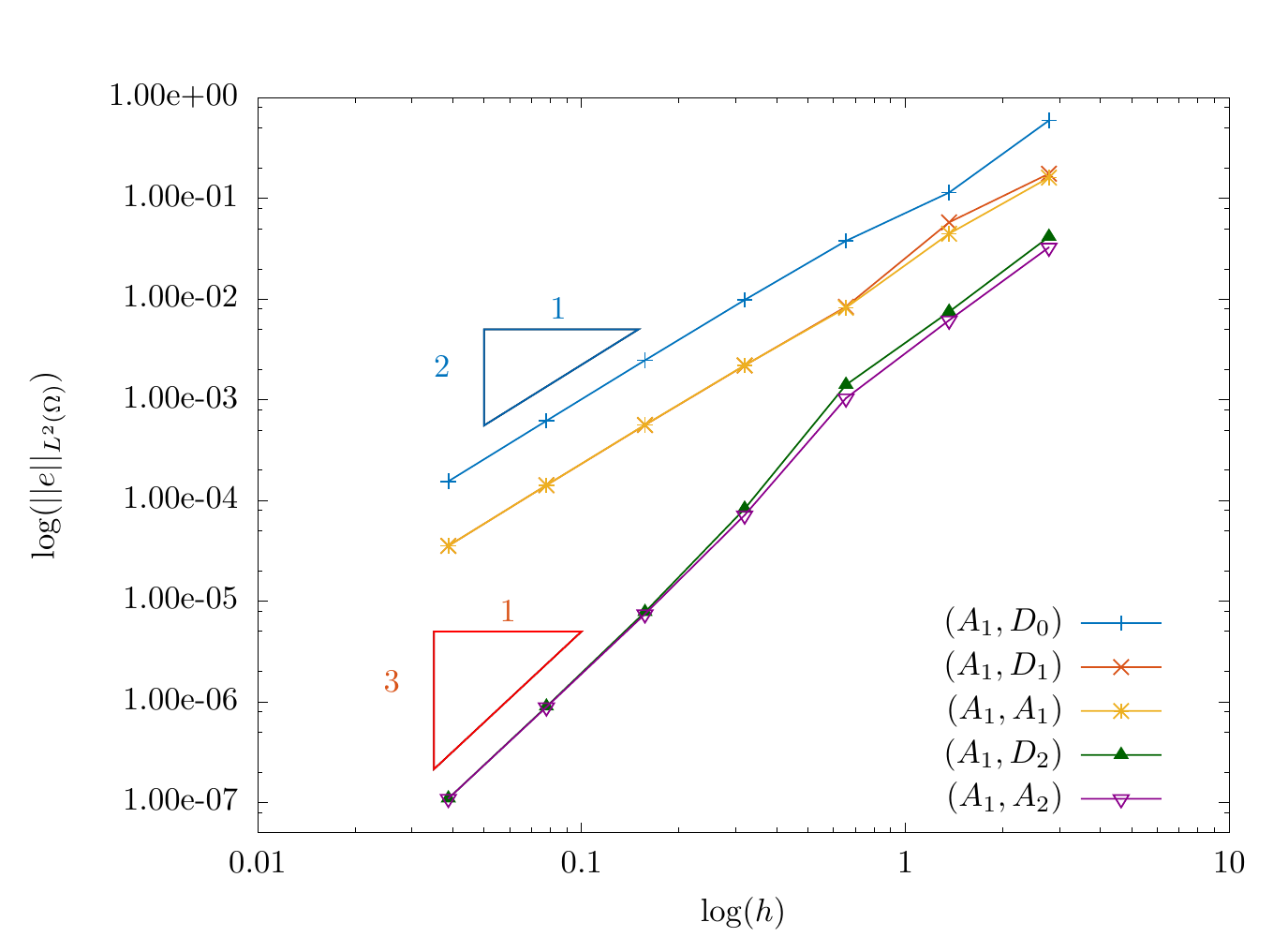}
\caption{Convergence study for Example~\ref{example_laplace_1}. B-splines solution bases: $D_{0}$, $D_{1}$, and $D_{2}$.}
\label{example_1_3}
\end{center}
\end{figure}

As a first numerical example, we consider the problem of Laplace equation in the quarter annulus domain of Fig.~\ref{fig:qannu}. We consider the exact solution
\begin{equation}
u (r,\theta)= r^{-3}\cos(3\theta),
\end{equation}
which is prescribed as Dirichlet boundary condition on all boundaries.
The convergence study cases are based on the choices of geometry and solution bases which are used in the patch tests in Section~\ref{patch_tests} (see Table~\ref{tab:res_tests}).
The convergence rates for all study cases are shown in Figures~\ref{example_1_1}, \ref{example_1_2}, and \ref{example_1_3}. The results are combined as follows.
In Fig.~\ref{example_1_1}, we have collected three choices of geometry parametrization, namely the coarsest $Q_{0}$, and its variations $A_{1}$, $A_{2}$, $B_{1}$, and $B_{2}$. Since $A_{1}$, $A_{2}$, $B_{1}$, and $B_{2}$ are equivalent to $Q_{0}$ up to knot-insertion and degree-elevation, the convergence results are identical, and the slope of the graph depends only on the lowest degree of the solution approximation ($1$ in the case of $A_{1}$, and $2$ in the case of $A_{2}$).
In Fig.~\ref{example_1_2}, we show the results for the non-uniform geometry parameterizations $C_{1}$ and $C_{2}$, and the solution bases $C_{1}$, $C_{2}$, $A_{1}$, and $A_{2}$, in comparison with the coarsest geometry parameterization $Q_{0}$ paired with the solution bases $C_{1}$ and $C_{2}$. Despite the slight difference in the results associated with the choice of the geometry representation, all graphs exhibit the expected convergence rates, which is governed by the lowest degree of the solution approximation.
Finally, in Fig.~\ref{example_1_3}, we included the results for the fifth case (represented by $T^{5}$ in Table~\ref{tab:res_tests}), characterized by the fact that the solution B-splines basis can not represent the geometry exactly. The pairs $A_{1}$-$D_{0}$, $A_{1}$-$D_{1}$, $A_{1}$-$D_{2}$ are shown in comparison with $A_{1}$-$A_{1}$ and $A_{1}$-$A_{2}$. Despite the difference in weights of the basis functions between $A_{1}$-$A_{1}$ and $A_{1}$-$D_{1}$ (as well as between $A_{1}$-$D_{2}$ and $A_{1}$-$A_{2}$) the difference in the results is very minor. This makes us conclude that though the B-splines solution basis fails the patch test, it is nevertheless a suitable basis for the analysis. Moreover, the case $A_{1}$-$D_{0}$, which can be regarded as truly super-parametric in FEA context, also exhibits the expected convergence rate. This surprising result counters the established practice of FEM, where super-parametric approach is not recommended \cite{ZienkiewiczTZ-FEM_V1}.
\subsubsection{Example 2}\label{example_3_Laplace}
\begin{figure}[!ht]
\begin{center}
\includegraphics[width=0.6\textwidth]{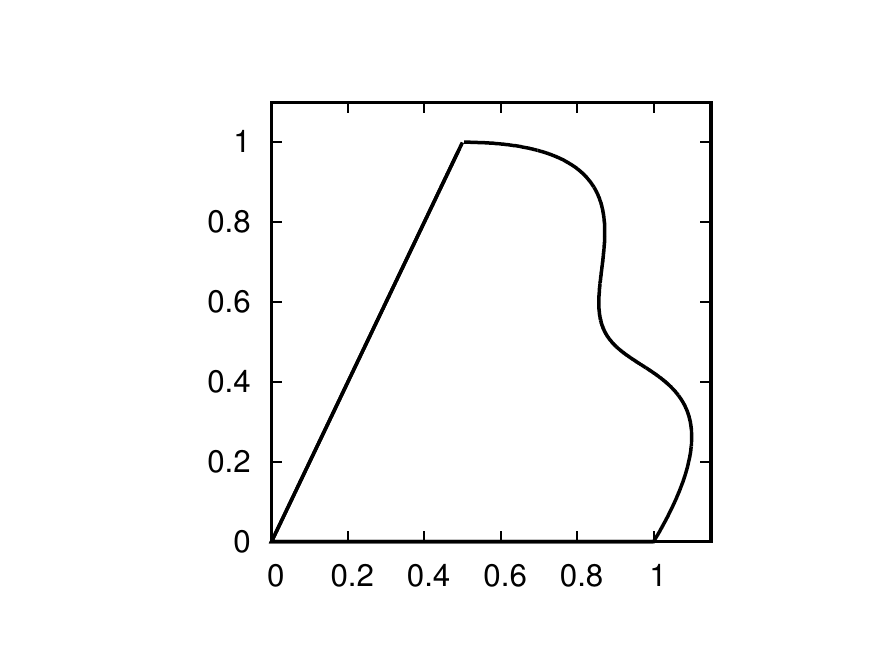}
\caption{Wedge geometry in Example~\ref{example_3_Laplace}.}
\label{wedge_geo}
\end{center}
\end{figure}
\begin{figure}[!ht]
\begin{center}
\includegraphics[width=0.9\textwidth]{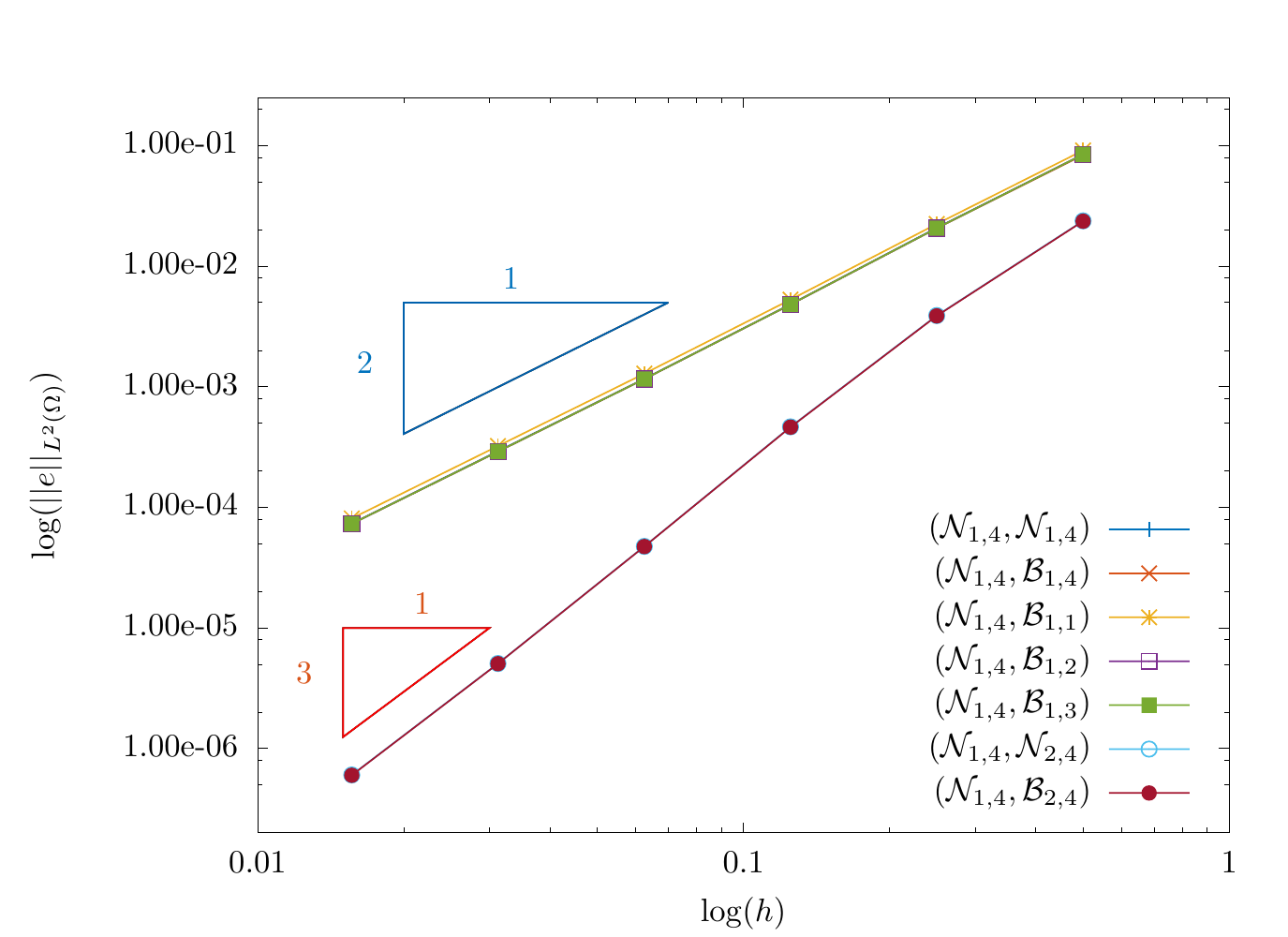}
\caption{Convergence study in example \ref{example_3_Laplace}.}
\label{fig:example_3_Laplace}
\end{center}
\end{figure}

In the second example, we study the Laplace problem in a wedge geometry, see Fig.~\ref{wedge_geo}.
The choice of this geometry is motivated by the observation in the super-parametric case of Example~\ref{example_laplace_1}.
As opposed to the patch test cases where the degree for geometry representation was artificially lifted, the geometry in this case is given by a tensor product of NURBS of degrees $p_{\xi} = 1$ and $p_{\eta} = 4$.
We consider the exact solution
\begin{equation}
u = \log((x + 0.1)^2 + (y + 0.1)^2),
\end{equation}
which is prescribed as Dirichlet boundary condition on all boundaries.
In this example we focus on sub- and super- parametric solution approximations. As the geometry parameterization is fixed (we denote the corresponding NURBS basis by $\mathcal{N}_{1,4}$), to ease with the naming, we use the notation of paired bases as $(\mathcal{N}_{1,4}, S_{k,l})$, where $S_{k,l}$ denotes the basis for solution approximation ($S = \mathcal{N}$ for NURBS and $S = \mathcal{B}$ for B-splines), $k$ denotes the degree in $\xi$ direction, and $l$ denotes the degree in $\eta$ direction. For sub-parametric solution approximations, we only elevate the degree in $\xi$ direction, i.e. $k = 2$, and for super-parametric solution approximations, we only consider degree reduction in $\eta$ direction, i.e. $1 \le l \le 3$. In Fig.~\ref{fig:example_3_Laplace}, we present the results for the following $7$ cases: $(\mathcal{N}_{1,4}, \mathcal{N}_{1,4})$, $(\mathcal{N}_{1,4}, \mathcal{B}_{1,4})$, $(\mathcal{N}_{1,4}, \mathcal{N}_{2,4})$, $(\mathcal{N}_{1,4}, \mathcal{B}_{2,4})$, $(\mathcal{N}_{1,4}, \mathcal{B}_{1,1})$, $(\mathcal{N}_{1,4}, \mathcal{B}_{1,2})$, and $(\mathcal{N}_{1,4}, \mathcal{B}_{1,3})$. 
As we can observe from the numerical studies presented in Fig.~\ref{fig:example_3_Laplace}, together with the exact representation of the geometry, the convergence rate depends only on the lowest degree in the solution approximation, and the results for NURBS and B-splines solution bases are almost identical.
\subsection{Linear elasticity in two-dimensions}\label{sec:lin_elas}
For linear elasticity problem, we choose two geometries, thick walled cylinder (quarter annulus for symmetry boundary conditions) and a plate with a circular hole (see Fig.~\ref{plate}).

\subsubsection{Example 1 (pressurized cylinder)}\label{example_elasticity_1}
\begin{figure}[!ht]
\begin{center}
\includegraphics[width=0.9\textwidth]{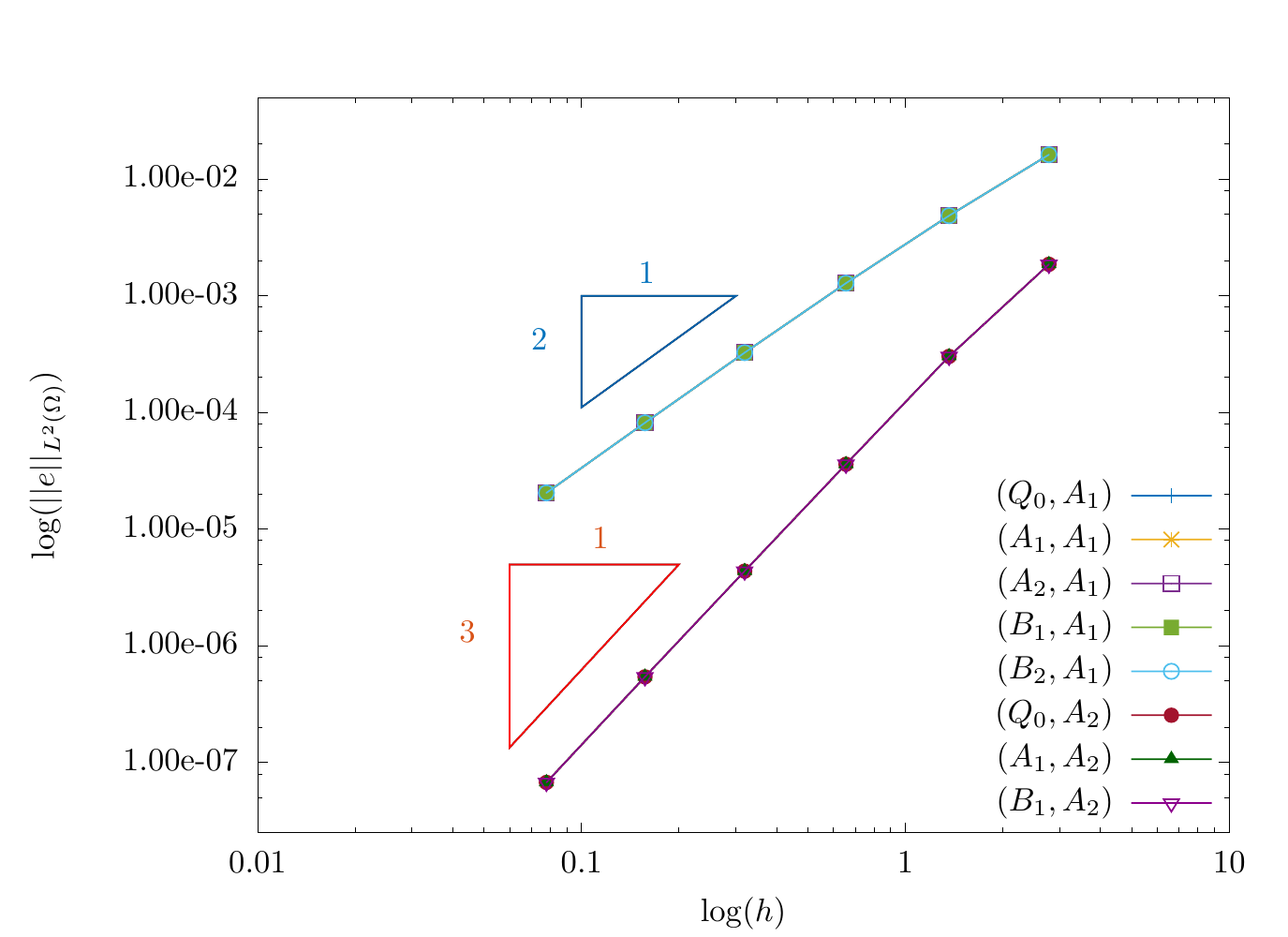}
\caption{Convergence study for Example~\ref{example_elasticity_1}. Geometry parameterizations $Q_{0}$, $A_{1}$, $A_{2}$, $B_{1}$ and $B_{2}$ combined with solution bases $A_{1}$ and $A_{2}$.}
\label{example_elasticity_1_plot_1}
\end{center}
\end{figure}

\begin{figure}[!ht]
\begin{center}
\includegraphics[width=0.9\textwidth]{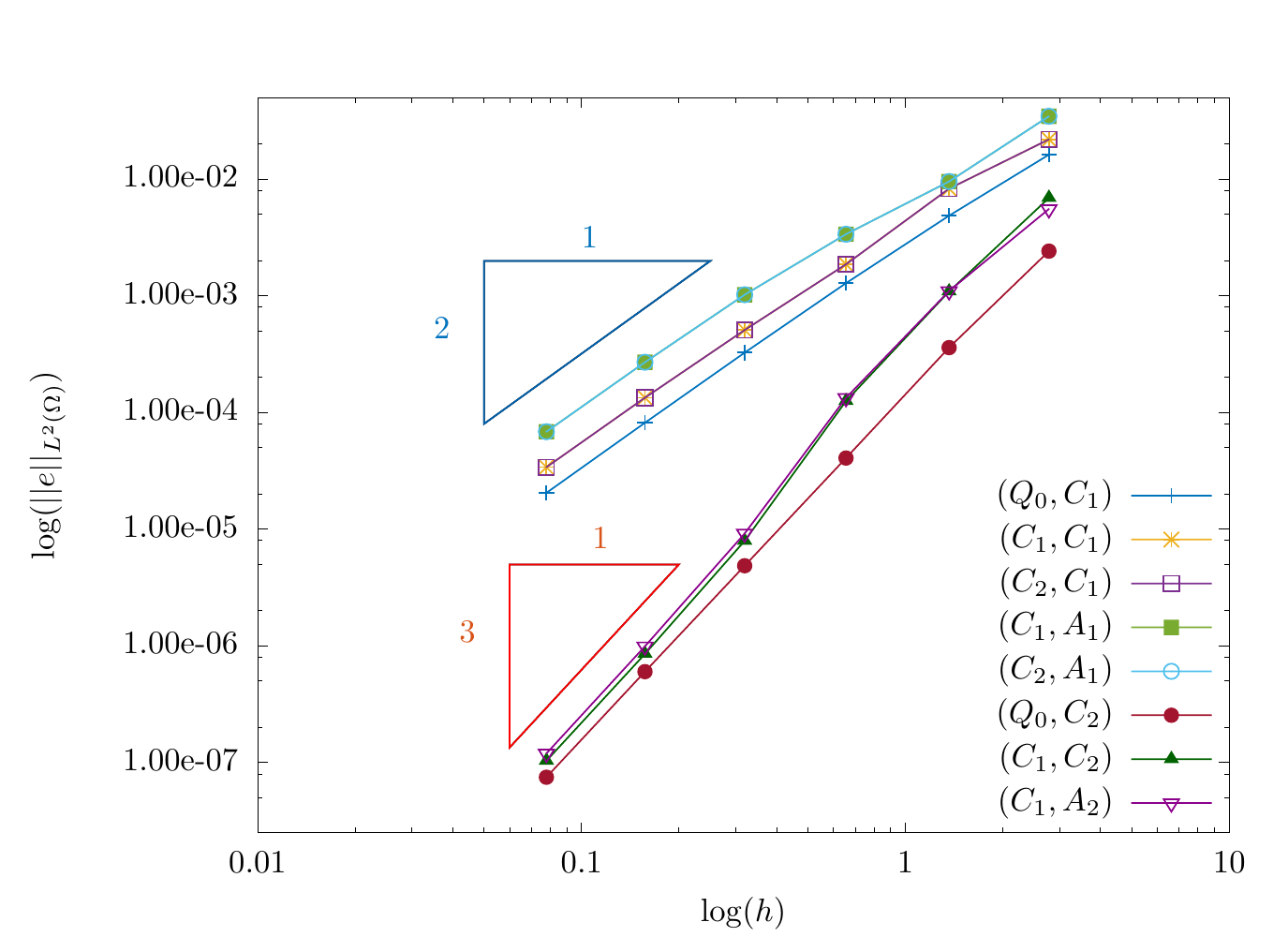}
\caption{Convergence study for Example~\ref{example_elasticity_1}. Non-uniform geometry parameterization $C_{1}$ and $C_{2}$ combined with solution bases $C_{1}$, $C_{2}$ and $A_{1}$, $A_{2}$.}
\label{example_elasticity_1_plot_2}
\end{center}
\end{figure}

\begin{figure}[!ht]
\begin{center}
\includegraphics[width=0.9\textwidth]{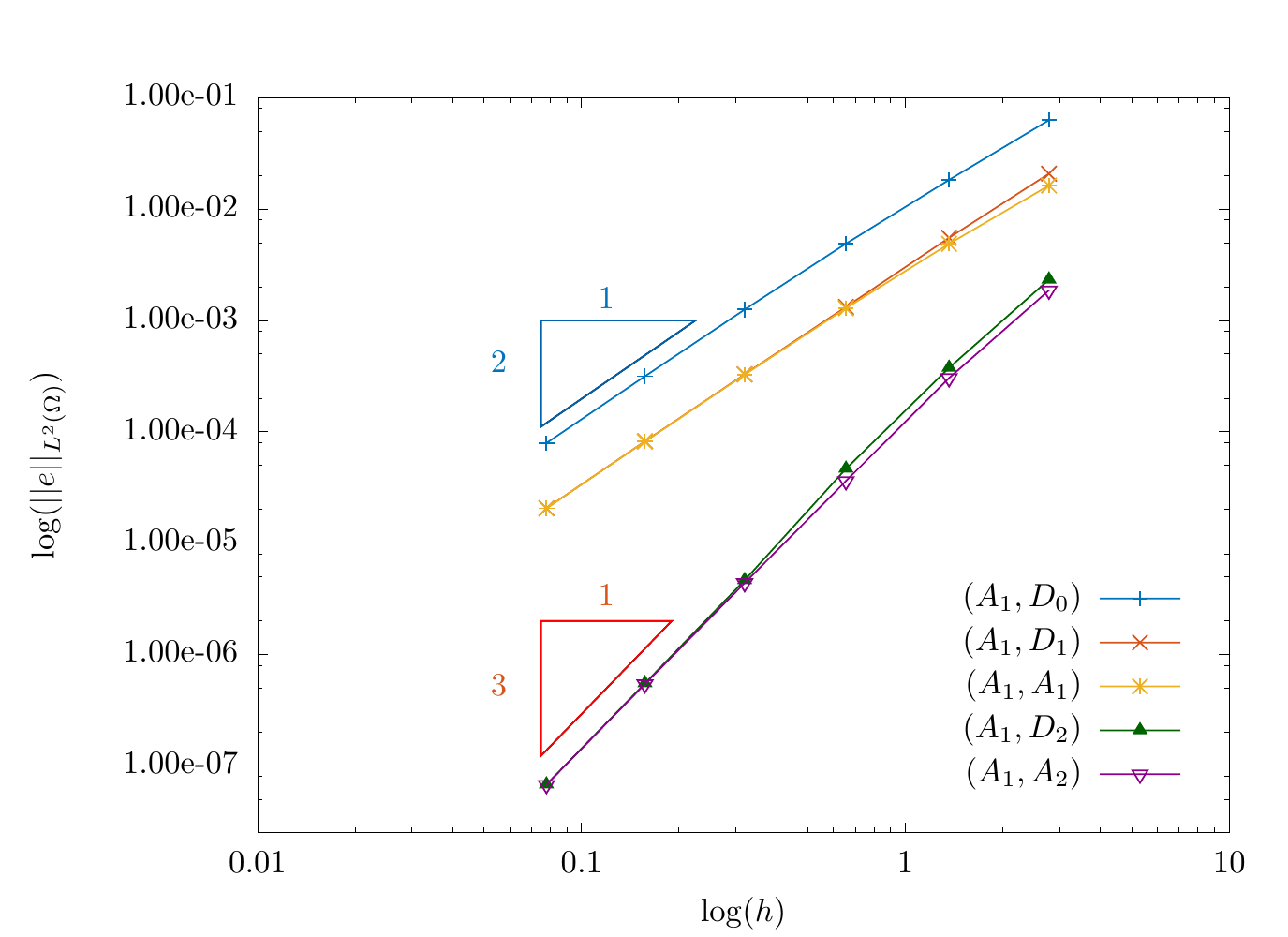}
\caption{Convergence study for Example~\ref{example_elasticity_1}. B-splines solution bases: $D_{0}$, $D_{1}$, $D_{2}$.}
\label{example_elasticity_1_plot_3}
\end{center}
\end{figure}

In the first example of linear elasticity, we consider a problem of a thick walled cylinder under a uniform internal and external pressure, as shown in Fig.~\ref{cylinder}. Instead of the full problem, a quarter annulus is considered with the symmetry boundary conditions along $x = 0$ and $y = 0$, and the following boundary conditions on the inner and outer boundary of the cylinder are imposed:
\begin{figure}[!ht]
\begin{center}
\includegraphics[width=0.5\textwidth]{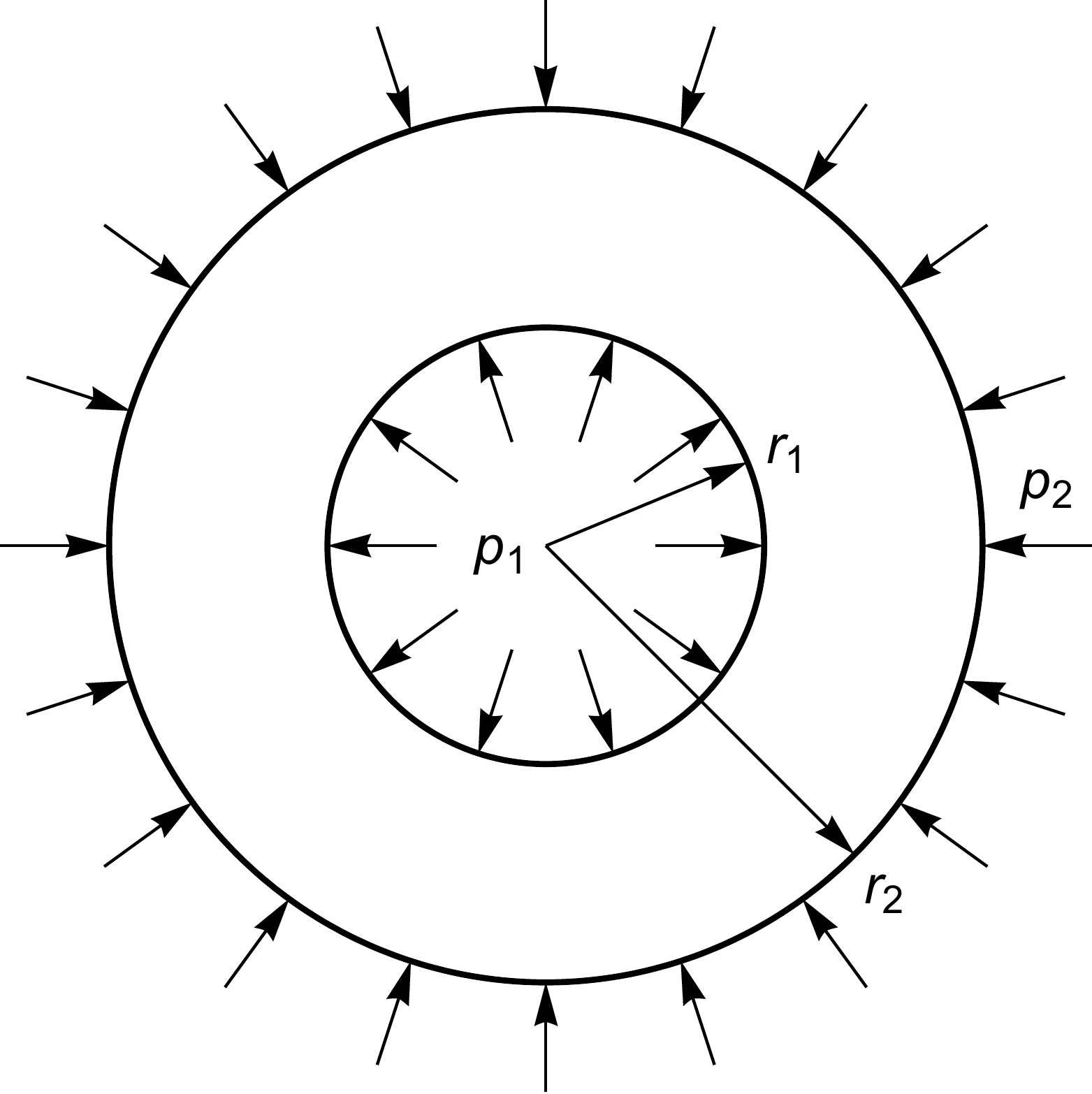}
\caption{Thick-walled pressurized cylinder}
\label{cylinder}
\end{center}
\end{figure}
\begin{equation}
\sigma_{r}(r_{1}) = -p_{1}, \qquad \sigma_{r}(r_{2}) = -p_{2}
\end{equation}
The solution to this problem in polar coordinates, under plane strain assumption, is given by
\begin{align}
\sigma_{r} (r)= s_{1}r^{-2} + s_{2}, \quad
\sigma_{\theta} (r)= -s_{1}r^{-2} + s_{2}, \quad
u_{r} (r)=  s_{3} \left(- s_{1} r^{-1} + s_{2} r\right),
\end{align}
where $s_{1} = \dfrac{r_{1}^2 r_{2}^2 (p_{2} - p_{1})}{r_{2}^2 - r_{1}^2}$, $s_{2} = (1 - 2 \nu) \dfrac{r_{1}^2 p_{1} - r_{2}^2 p_{2}}{r_{2}^2 - r_{1}^2}$, and $s_{3} = \dfrac{1 + \nu}{E}$.
The numerical tests were performed on the same choices of the combination of geometry-solution bases, as in the first example for Laplace equation (see Table~\ref{tab:res_tests}). The results of numerical simulations are organized in three plots: Fig.~\ref{example_elasticity_1_plot_1}, \ref{example_elasticity_1_plot_2} and \ref{example_elasticity_1_plot_3}. The obtained numerical solutions exhibit a similar pattern as that observed in Example~\ref{example_laplace_1}, i.e. the graphs show that the convergence rate of the solution is defined by the lowest degree in the solution basis, independent of the (exact) geometry parameterization.
\subsubsection{Example 2 (plate with a circular hole)}\label{example_elasticity_2}
\begin{figure}[!ht]
\centering
\begin{subfigure}[b]{0.425\textwidth}
\includegraphics[width=\textwidth]{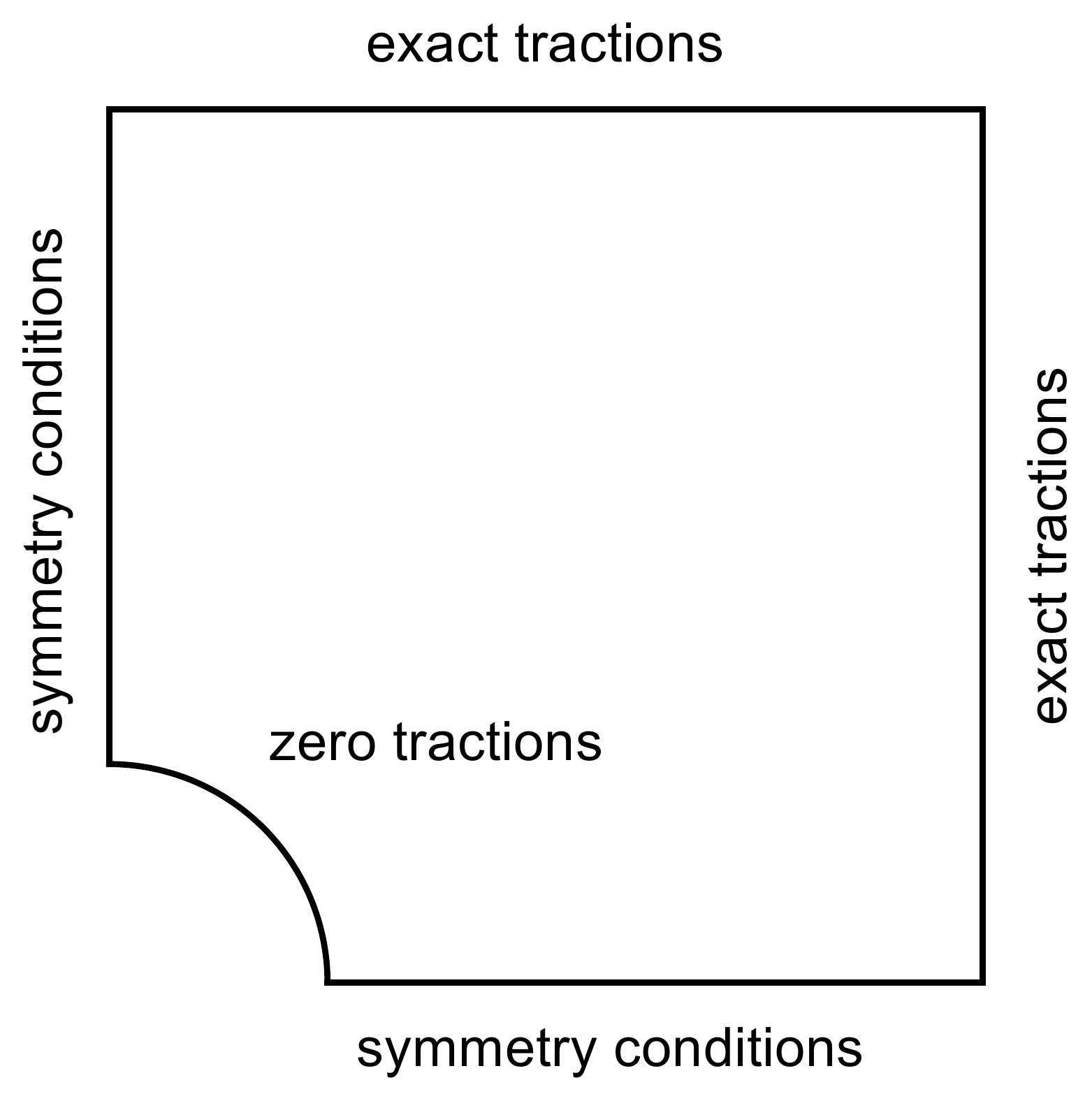}
\caption{Quarter of the plate with a hole, subjected to the boundary conditions given by Eq.~\eqref{plate_solution}.}
\label{plate}
\end{subfigure}
\quad
\begin{subfigure}[b]{0.425\textwidth}
\includegraphics[width=\textwidth]{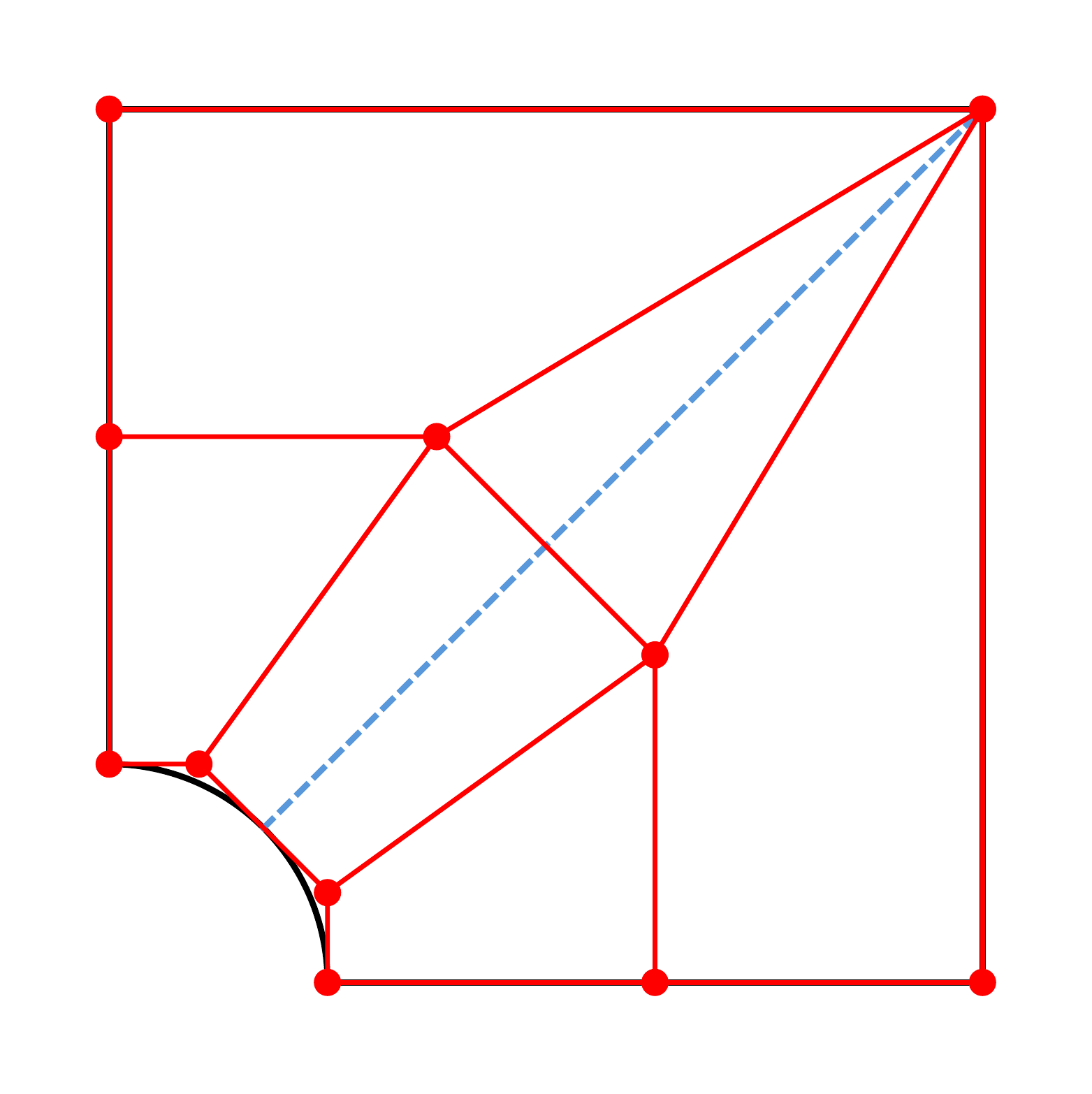}
\caption{Control net and element boundaries in the plate parametrization.}
\label{plate_param}
\end{subfigure}
\caption{Geometry and parameterization for Example~\ref{example_elasticity_2}.}
\end{figure}

\begin{figure}[!ht]
\begin{center}
\includegraphics[width=0.9\textwidth]{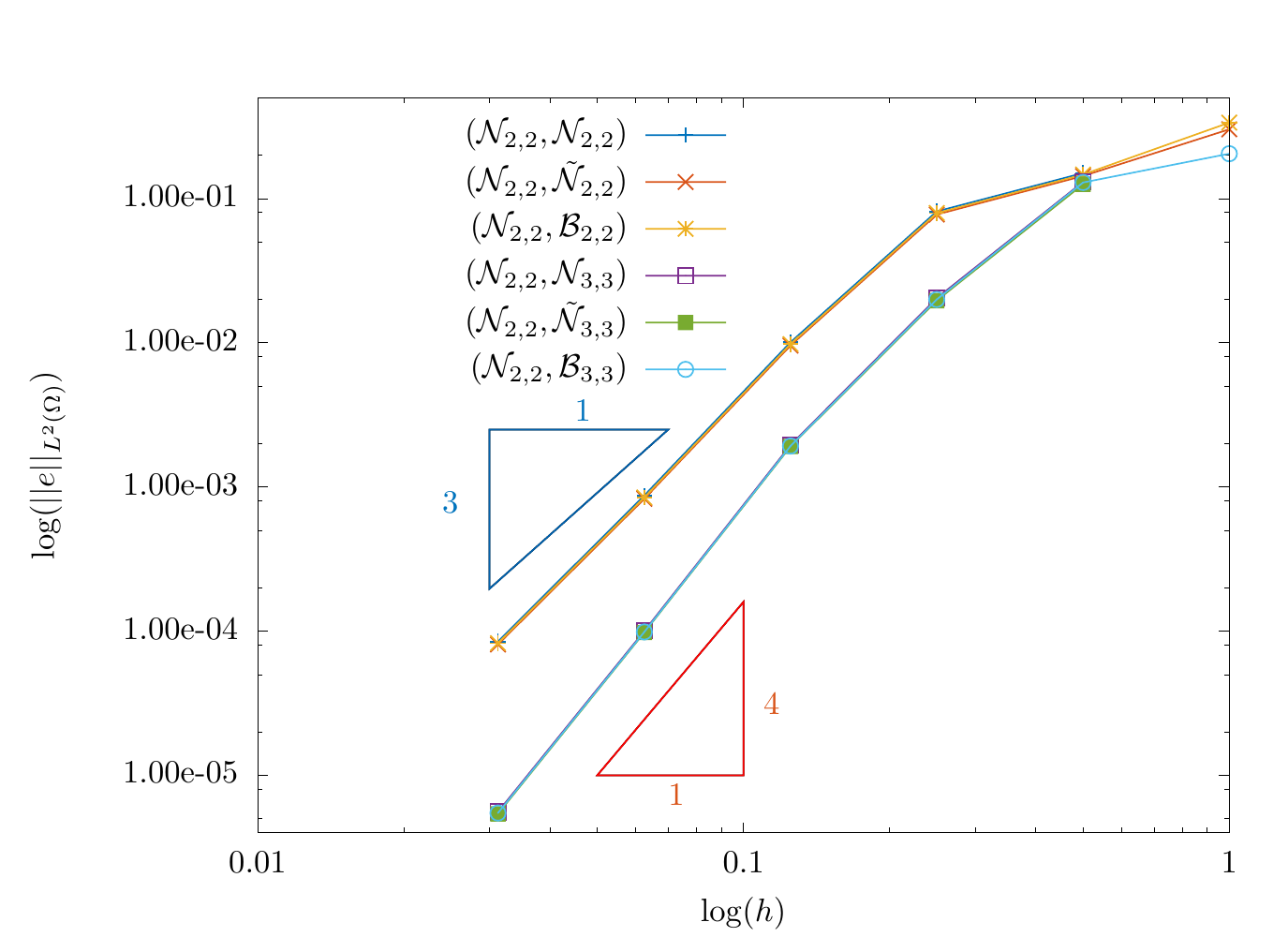}
\caption{Convergence study for Example~\ref{example_elasticity_2}.}
\label{example_elasticity_2_plot_1}
\end{center}
\end{figure}

In the second example of linear elasticity, we consider a typical problem of a plate weakened by a circular hole of radius $a$, and subject to remote tension $T$ in $x$-direction. Only a quarter of a plate (of a finite size $L = 4 a$), as shown in Fig.~\ref{plate}, is modeled with the symmetry boundary conditions along $x = 0$ and $y = 0$, and analytical tractions prescribed on the rest of the boundary according to the solution below:
\begin{equation}
\begin{split}
\sigma_{xx} (r,\theta)&= T - T\dfrac{a^2}{r^2}\left(\dfrac{3}{2}\cos2\theta + \cos 4\theta \right) + T\dfrac{3 a^4}{2 r^4}\cos 4\theta, \\
\sigma_{yy} (r,\theta)&= - T\dfrac{a^2}{r^2}\left(\dfrac{1}{2}\cos2\theta - \cos 4\theta \right) - T\dfrac{3 a^4}{2 r^4}\cos 4\theta, \\
\sigma_{xy} (r,\theta)&= - T\dfrac{a^2}{r^2}\left(\dfrac{1}{2}\sin2\theta + \sin 4\theta \right) + T\dfrac{3 a^4}{2 r^4}\sin 4\theta,\\
u_{x} (r,\theta)&= \alpha_{0} \left(\dfrac{4 r}{a}(1 - \nu)\cos\theta + \dfrac{2 a}{r}(4 (1 - \nu)\cos\theta + \cos 3\theta) - \dfrac{2 a^3}{r^3}\cos 3\theta\right), \\
u_{y}(r,\theta) &= \alpha_{0} \left(\dfrac{r}{a}(-4 \nu)\sin\theta + \dfrac{2 a}{r}(-2(1 - 2\nu)\sin\theta + \sin 3\theta) - \dfrac{2 a^3}{r^3}\sin 3\theta\right),
\end{split}
\label{plate_solution}
\end{equation}
where $\alpha_{0} = (1 + \nu)T a / (4 E)$. In all study cases for this problem, the geometry is parametrized by a basis of second order with the following knot vectors:
\begin{equation}
\Sigma = \{0,0,0,1,1,1\}, \qquad \Pi = \{0,0,0,0.5,1,1,1\}.
\label{knot_vectors_plate}
\end{equation}
The corresponding control points are listed in Table~\ref{cpt_plate}. The parametrization consists of two elements (see Fig.~\ref{plate_param}), and remains unchanged during the solution refinement process (denoted by $\mathcal{N}_{2,2}$).
Following the notation of paired bases introduced in Section~\ref{example_3_Laplace}, we consider the following choices for the solution approximation basis: $(\mathcal{N}_{2,2}, \mathcal{N}_{2,2})$, $(\mathcal{N}_{2,2}, \mathcal{N}_{3,3})$, $(\mathcal{N}_{2,2}, \mathcal{B}_{2,2})$, and $(\mathcal{N}_{2,2}, \mathcal{B}_{3,3})$. Note, that B-Spline basis $\mathcal{B}_{2,2}$ is built on knot vectors (\ref{knot_vectors_plate}).
We also consider another geometry parameterization which is built on the knot vectors (\ref{knot_vectors_plate}), with the boundary control points as listed in Table~\ref{cpt_plate}, but the weights of two inner points being changed to $w_{2,2} = w_{2,3} = 0.9$. This parameterization will be denoted by $\tilde{\mathcal{N}}_{2,2}$. Thereafter, we consider two more choices for the solution approximation basis: $(\mathcal{N}_{2,2}, \tilde{\mathcal{N}}_{2,2})$, and $(\mathcal{N}_{2,2}, \tilde{\mathcal{N}}_{3,3})$.
The results of the these study cases are shown in Fig.~\ref{example_elasticity_2_plot_1}, where it can be seen that the convergence rate for all choices of the solution basis depends only on the order of the solution basis, and for the bases of the same order the results for different basis functions are almost identical.
\begin{table}
\begin{center}
\begin{tabular}{|c|c|c|c|}
\hline 
$(i,j)$ & $P^{x}_{ij}$ & $P^{y}_{ij}$ & $w_{ij}$\\ 
\hline 
(1,1) & 1.0000 & 0.0000 & 1.0000\\
(1,2) & 1.0000 & 0.4142 & 0.8536\\
(1,3) & 0.4142 & 1.0000 & 0.8536\\
(1,4) & 0.0000 & 1.0000 & 1.0000\\
(2,1) & 2.5000 & 0.0000 & 1.0000\\
(2,2) & 2.5000 & 1.5000 & 0.8000\\
(2,3) & 1.5000 & 2.5000 & 0.8000\\
(2,4) & 0.0000 & 2.5000 & 1.0000\\
(3,1) & 4.0000 & 0.0000 & 1.0000\\
(3,2) & 4.0000 & 4.0000 & 1.0000\\
(3,3) & 4.0000 & 4.0000 & 1.0000\\
(3,4) & 0.0000 & 4.0000 & 1.0000\\
\hline 
\end{tabular} 
\caption{Control points in the parametrization of Example~\ref{example_elasticity_2}.}
\label{cpt_plate}
\end{center}
\end{table}
\subsection{Linear elasticity in three-dimensions}\label{sec:lin_elas_3d}
\begin{figure}[!ht]
\begin{center}
\includegraphics[width=0.5\textwidth]{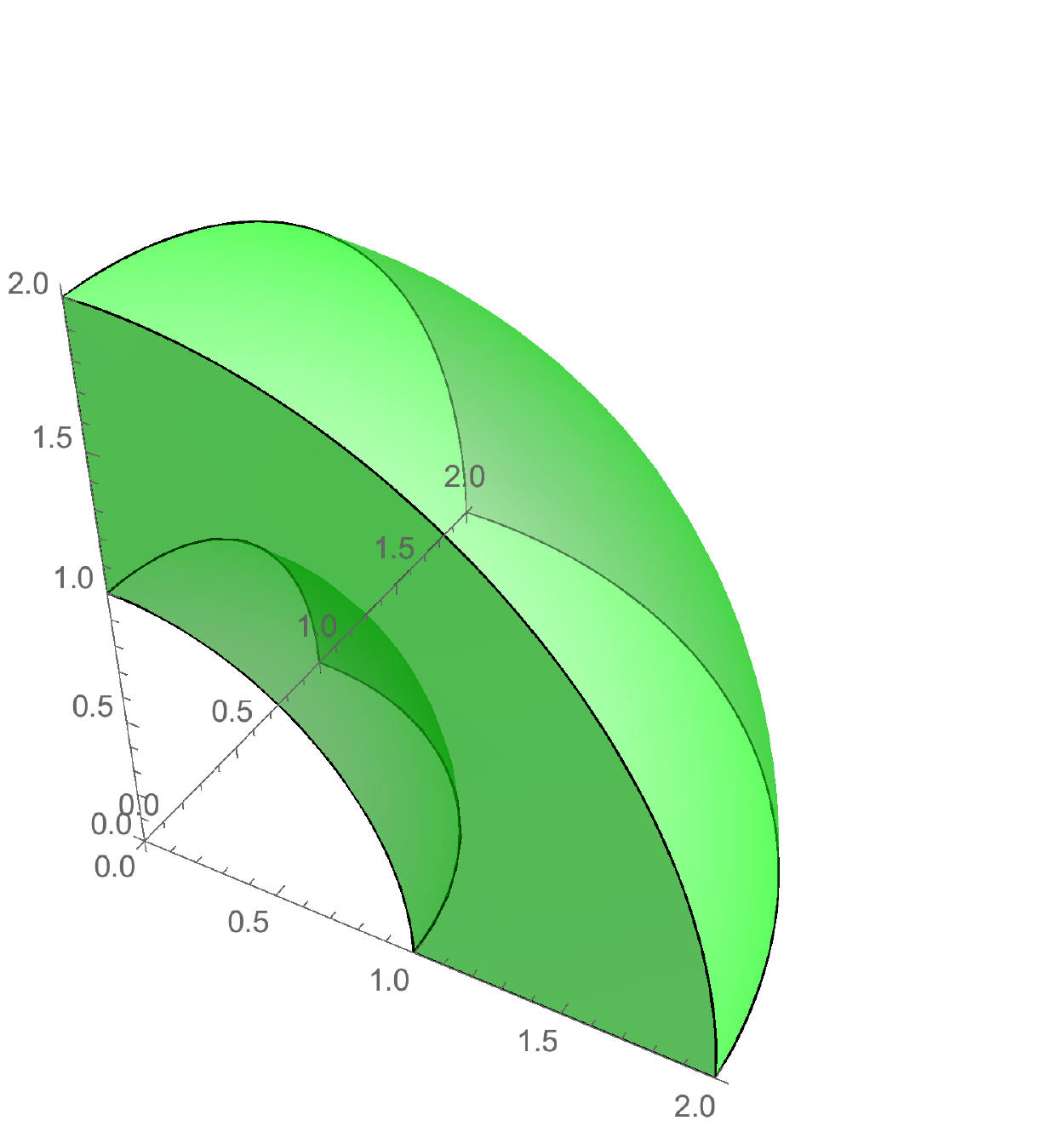}
\caption{One-eighth of a sphere.}\label{sphere_geo}
\end{center}
\end{figure}

\begin{figure}[!ht]
\begin{center}
\includegraphics[width=0.9\textwidth]{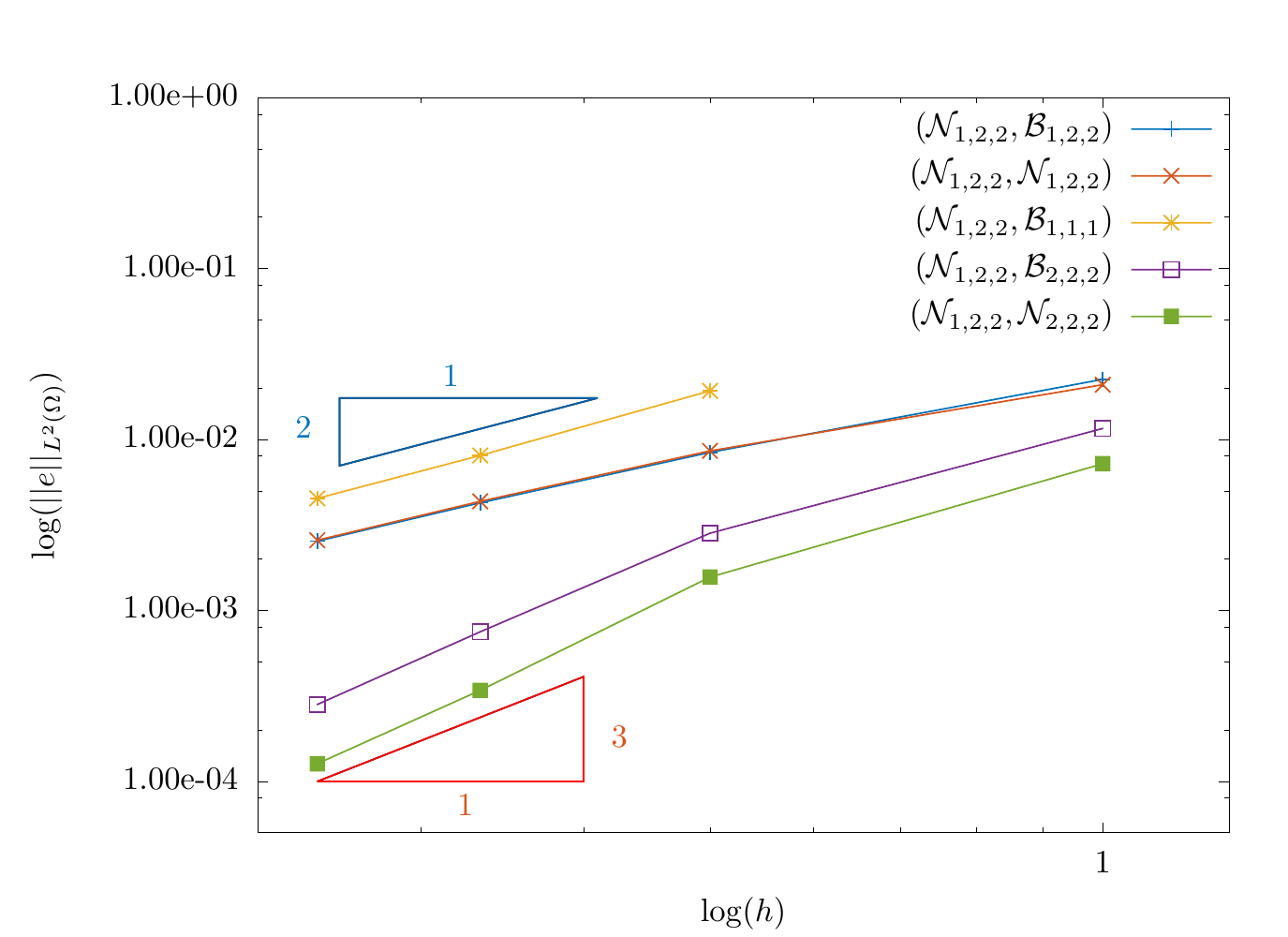}
\caption{Convergence study for the 3D example.}\label{convergence3D}
\end{center}
\end{figure}

For three-dimensional linear elasticity problem, we choose a thick walled sphere (one-eighth of a hollow sphere for symmetry boundary conditions) geometry with inner radius $r_{1}$ and outer radius $r_2$, see Fig. \ref{sphere_geo}. The coarsest NURBS parametrization of this geometry can be given by degrees $p_{\xi} = 1$, $p_{\eta} = 2$ and $p_{\zeta} = 2$, built upon the knot vectors:
\begin{equation}
\Sigma = \{0, 0, 1, 1\}, \quad \Pi = \{0, 0, 0, 1, 1, 1\}, \quad Z = \{0, 0, 0, 1, 1, 1\},
\end{equation}
with the control points listed in Table~\ref{cpt_sphere}. We will refer to this parametrization as $Q_{1}$. Using the notations introduced earlier, the basis for this parameterization will be denoted by $\mathcal{N}_{1, 2, 2}$.
\begin{table}
\begin{center}
\begin{tabular}{|c|c|c|c|c|}
\hline
$(i, j, k)$ & $P^{x}_{ij}$ & $P^{y}_{ij}$ & $P^{z}_{ij}$ & $w_{ij}$\\
\hline
(1,1,1) & 1.0000 & 0.0000 & 0.0000 & 1.0000 \\
(1,2,1) & 1.0000 & 1.0000 & 0.0000 & 0.7071 \\
(1,3,1) & 0.0000 & 1.0000 & 0.0000 & 1.0000 \\
(2,1,1) & 2.0000 & 0.0000 & 0.0000 & 1.0000 \\
(2,2,1) & 2.0000 & 2.0000 & 0.0000 & 0.7071 \\
(2,3,1) & 0.0000 & 2.0000 & 0.0000 & 1.0000 \\
(1,1,2) & 1.0000 & 0.0000 & 1.0000 & 0.7071 \\
(1,2,2) & 1.0000 & 1.0000 & 1.0000 & 0.5000 \\
(1,3,2) & 0.0000 & 1.0000 & 0.0000 & 0.7071 \\
(2,1,2) & 2.0000 & 0.0000 & 2.0000 & 0.7071 \\
(2,2,2) & 2.0000 & 2.0000 & 2.0000 & 0.5000 \\
(2,3,2) & 0.0000 & 2.0000 & 0.0000 & 0.7071 \\
(1,1,3) & 0.0000 & 0.0000 & 1.0000 & 1.0000 \\
(1,2,3) & 0.0000 & 1.0000 & 1.0000 & 0.7071 \\
(1,3,3) & 0.0000 & 1.0000 & 0.0000 & 1.0000 \\
(2,1,3) & 0.0000 & 0.0000 & 2.0000 & 1.0000 \\
(2,2,3) & 0.0000 & 2.0000 & 2.0000 & 0.7071 \\
(2,3,3) & 0.0000 & 2.0000 & 0.0000 & 1.0000 \\
\hline
\end{tabular}
\caption{Control points for the parametrization of one-eighth of a sphere.}
\label{cpt_sphere}
\end{center}
\end{table}
We also consider the following choices of the bases for solution approximation: $\mathcal{N}_{1, 2, 2}$, $\mathcal{N}_{2, 2, 2}$, $\mathcal{B}_{1, 2, 2}$, $\mathcal{B}_{2, 2, 2}$, and $\mathcal{B}_{1, 1, 1}$. Note that, we use uniform refinement of $Q_{1}$ for these bases (similar to the case $A_{1}$ of Section~\ref{sec:geo_param}). And, the basis $\mathcal{B}_{1, 1, 1}$ is built on the knot vectors:
\begin{equation}
\Sigma = \{0,0,1,1\}, \,\,\, \Pi = \{0,0,1,1\}, \,\,\, Z = \{0,0,1,1\}.
\end{equation}
The problem in consideration is a well known thick-walled pressurized sphere, i.e. the hollow sphere subjected to the following boundary conditions (in spherical coordinates $r, \theta, \phi$):
\begin{equation}
\sigma_{rr} = -p_{1} \quad \text{at} \quad r = r_{1}, \qquad
\sigma_{rr} = -p_{2} \quad \text{at} \quad r = r_{2}.
\end{equation}
On the rest of the boundary the symmetry conditions are prescribed. The analytical solution is given by
\begin{equation}
u_{r} (r) = \left(\alpha_{1} r + \alpha_{2} r^{-2}\right) / \alpha_{3},
\end{equation}
where $\alpha_{1} = 2(p_{1} r_1^3 - p_{2} r_2^3) (1 - 2\nu)$, $\alpha_{2} = (p_1 - p_2) (1 + \nu) r_{1}^{3} r_{2}^{3}$, and $\alpha_{3} = 2 E (r_2^3 - r_1^3)$. The numerical results are shown in Fig.~\ref{convergence3D}, where it can again be seen that the convergence rate in all the five cases depends only on the approximation basis for the numerical solution.
\subsection{Numerical solution with PHT splines}\label{sec:laplace_pht}
To further demonstrate the effectiveness of the proposed method, and the use of weak coupling of bases for geometry and simulation, in this section we present numerical results using PHT-splines (to avoid the disruption of this section's readability, the details of PHT-splines based construction is provided in Appendix~\ref{sec:app_PHT}).
We again consider an annulus region in two-dimensions, see Fig.~\ref{fig:qannu}, described by a quadratic $C^1$ NURBS surface with $6 \times 6$ control points, and following knot vectors on the parametric domain
\[\Xi = [0, 0, 0, 0.04, 0.2, 0.36, 1, 1, 1], \quad \mathcal{H} = [0, 0, 0, 0.04, 0.2, 0.36, 1, 1, 1].\]
We consider the Poisson problem, and choose the source function such that the exact solution of the problem has the following form
\begin{equation}
u (r,\theta)= (r - 1) (r - 2) \theta (\theta - \frac{\pi}{2}) \exp \bigl ( -100 (r \cos \theta - 1)^2 \bigr),
\end{equation}
where
\[
r (x,y) = \sqrt{x^2 + y^2}, \quad \theta = \arctan({y/x}).
\]
In Fig.~\ref{fig:Annus2}, we show the results during local refinement operations. From the first row to the fourth row, we show the T-mesh on the parametric domain (left), the numerical solution with $C^1$ PHT-splines representation (middle), and the corresponding exact error color-map (right) on the parametric domain. In Fig.~\ref{fig:AnnusPHTerror}, we present the convergence behavior of four choices of geometry and simulation bases, which are as follows:
\begin{itemize}
\item IGA with cubic NURBS (for the geometry as well as numerical solution). Note that, a quadratic NURBS is sufficient for this geometry, however, to have a fair comparison with the remaining studies, we elevate the degree while maintaining the exact geometry representation.
\item IGA with cubic PHT-splines (for the geometry as well as the numerical solution). Note that, in this case, the computational geometry is only approximate (not exact as in IGA with cubic NURBS).
\item GIFT with cubic B-splines for the numerical solution, and quadratic NURBS for exact geometry representation.
\item GIFT with cubic PHT-splines for the numerical solution, and quadratic NURBS for exact geometry representation.
\end{itemize}
From the convergence plots, we conclude the following:
\begin{enumerate}
\item As observed in earlier studies, IGA with cubic NURBS and GIFT with cubic B-splines exhibit same convergence rate. This is because, with the exact geometry representation, both the solution bases are of the same degree.
\item Owing to the local adaptive refinement of the solution basis, both the cases of PHT-splines (IGA as well as the proposed GIFT) exhibit higher convergence rate than the tensor-product based solution bases (B-splines as well as NURBS).
\item The comparison of IGA with PHT-splines and GIFT with PHT-splines solution highlights an important difference. The advantage of the exact geometry representation in the latter case over an approximate geometry in the former case is very minor. This is due to the fact that the geometry of the computational domain is simple, and can be accurately approximated with the PHT-splines of third degree (geometry approximation error is below the discretization error). However, in realistic industrial problems with complex domains, this advantage will become more pronounced. In comparison to standard IGA with PHT splines, employing the exact coarse NURBS geometry parametrization in GIFT (together with PHT-splines solution) brings two distinct advantages:
\begin{itemize}
\item It eliminates the need to communicate with the original CAD model at each step of the solution refinement process, and the approximation of the boundaries.
\item It also eliminates the need to refine the original coarse geometry, as well as to store and process the refined data, which can lead to significant computational savings for big problems.
\end{itemize} 
\end{enumerate}
\setlength{\picw}{1.9in}
\begin{figure}
\centering
\includegraphics[width=\picw]{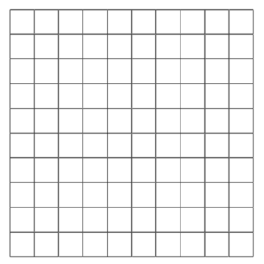}
\quad
\includegraphics[width=\picw]{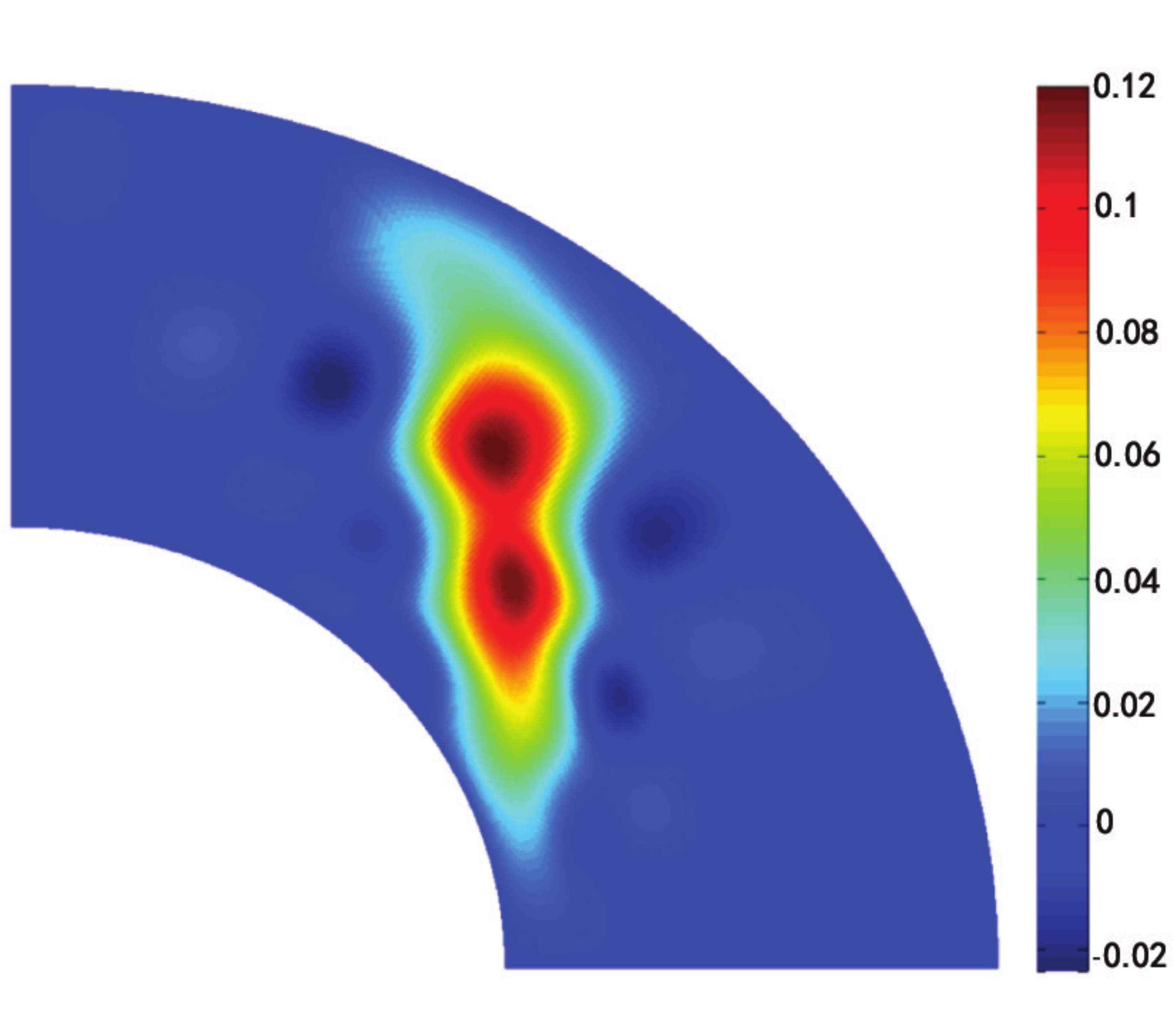}
\quad
\includegraphics[width=\picw]{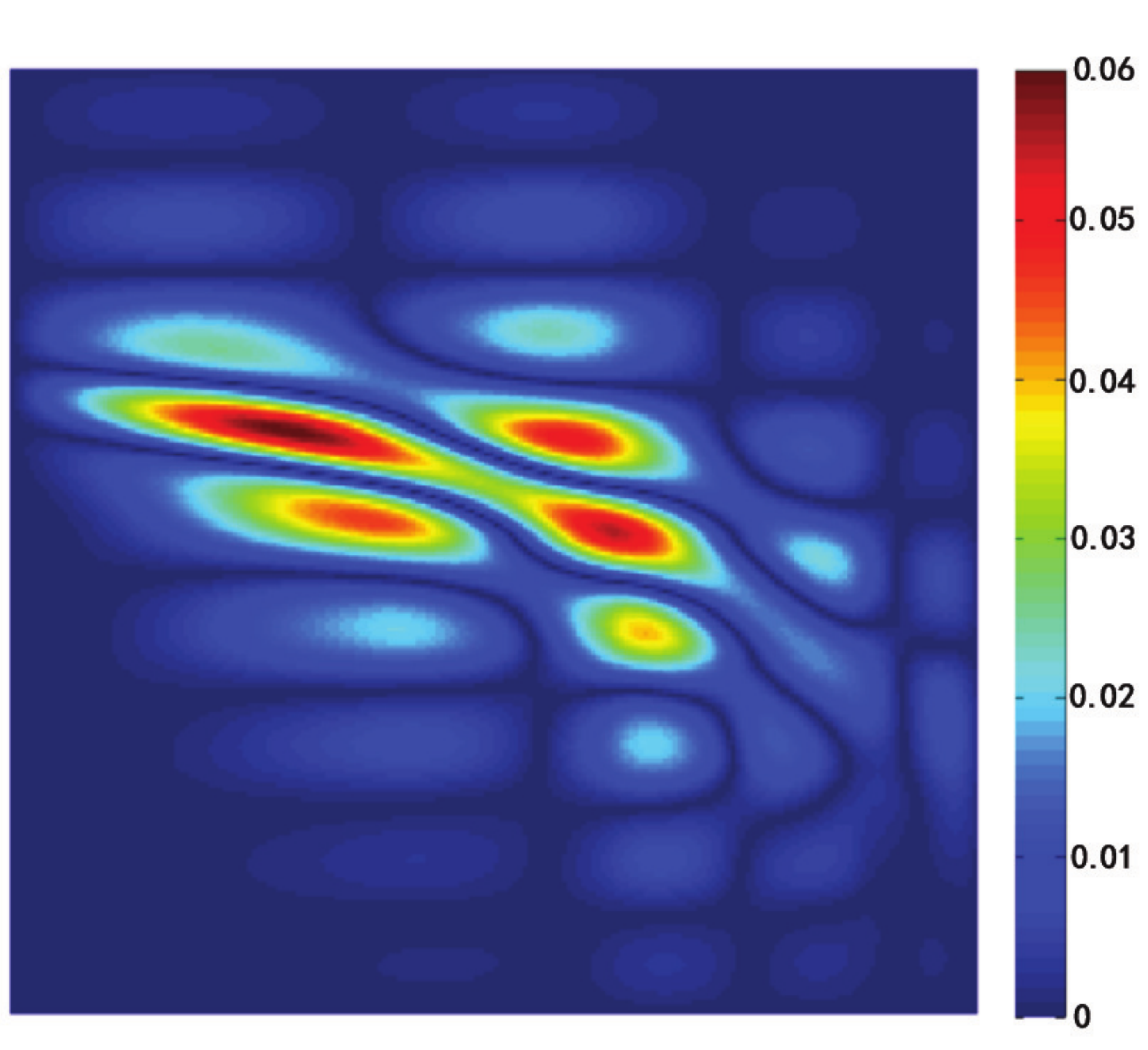}
\\
\includegraphics[width=\picw]{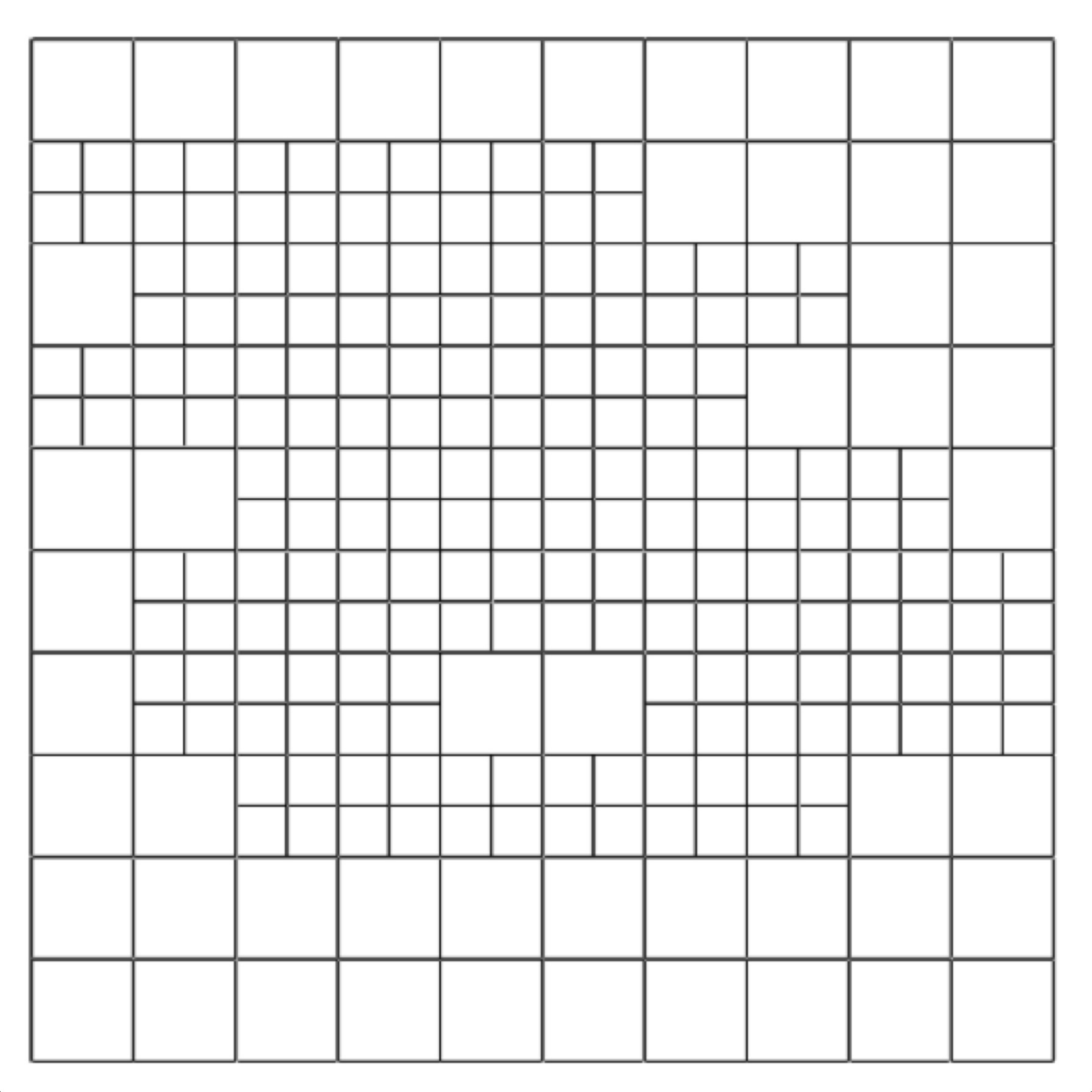}
\quad
\includegraphics[width=\picw]{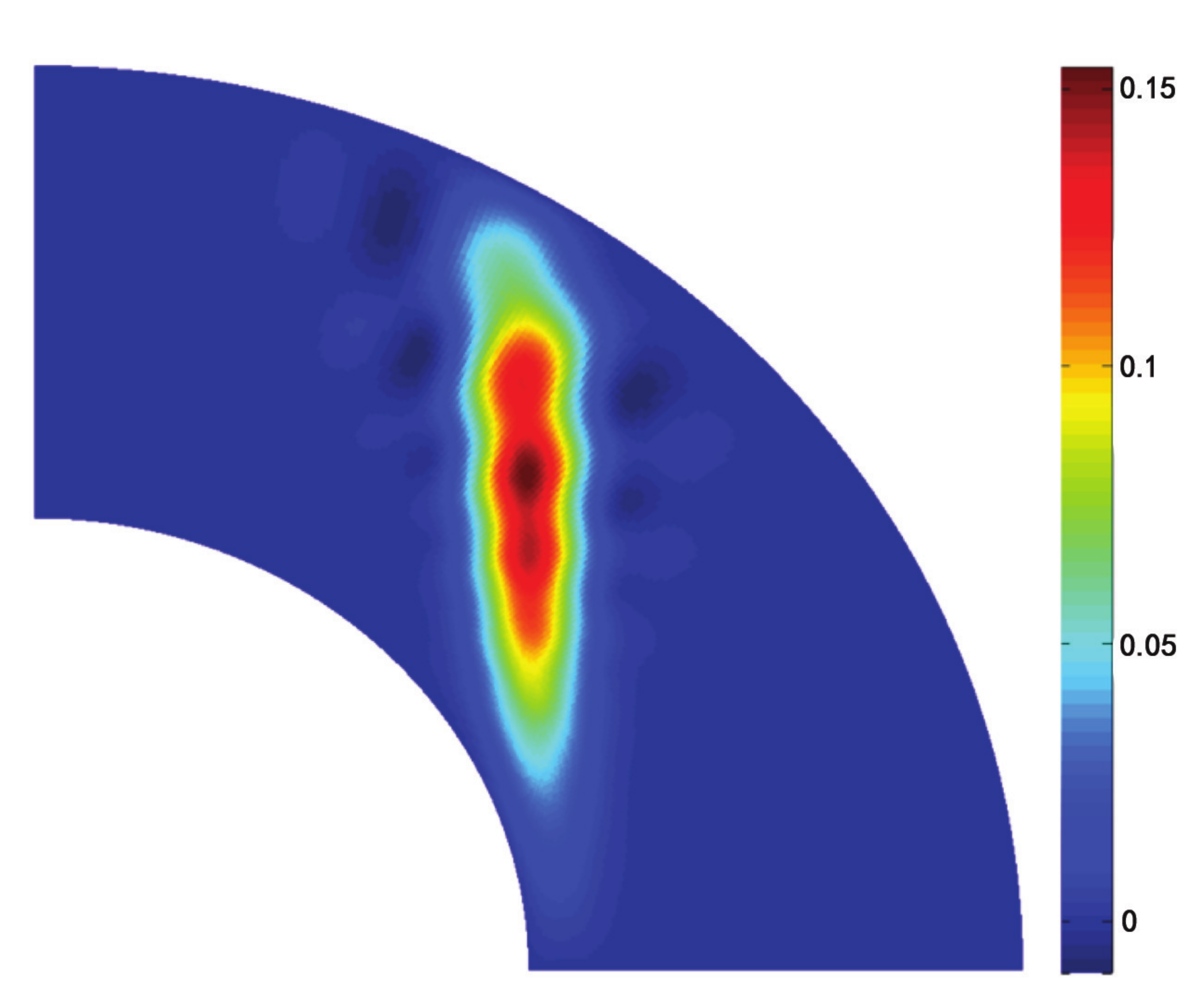}
\quad
\includegraphics[width=\picw]{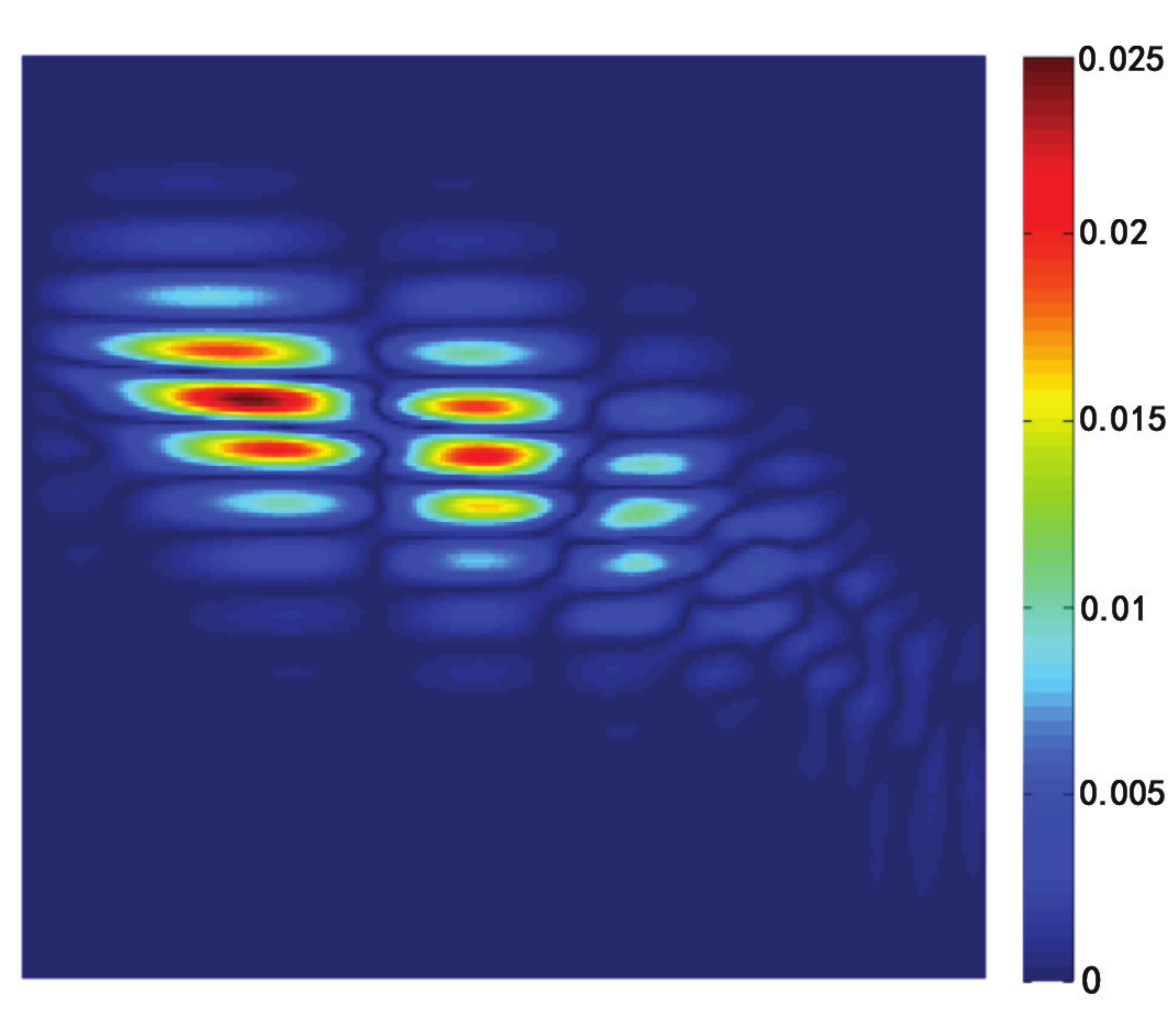}
\\
\includegraphics[width=\picw]{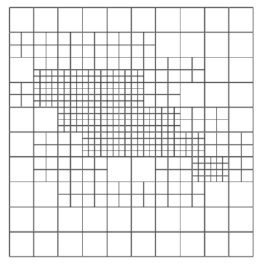}
\quad
\includegraphics[width=\picw]{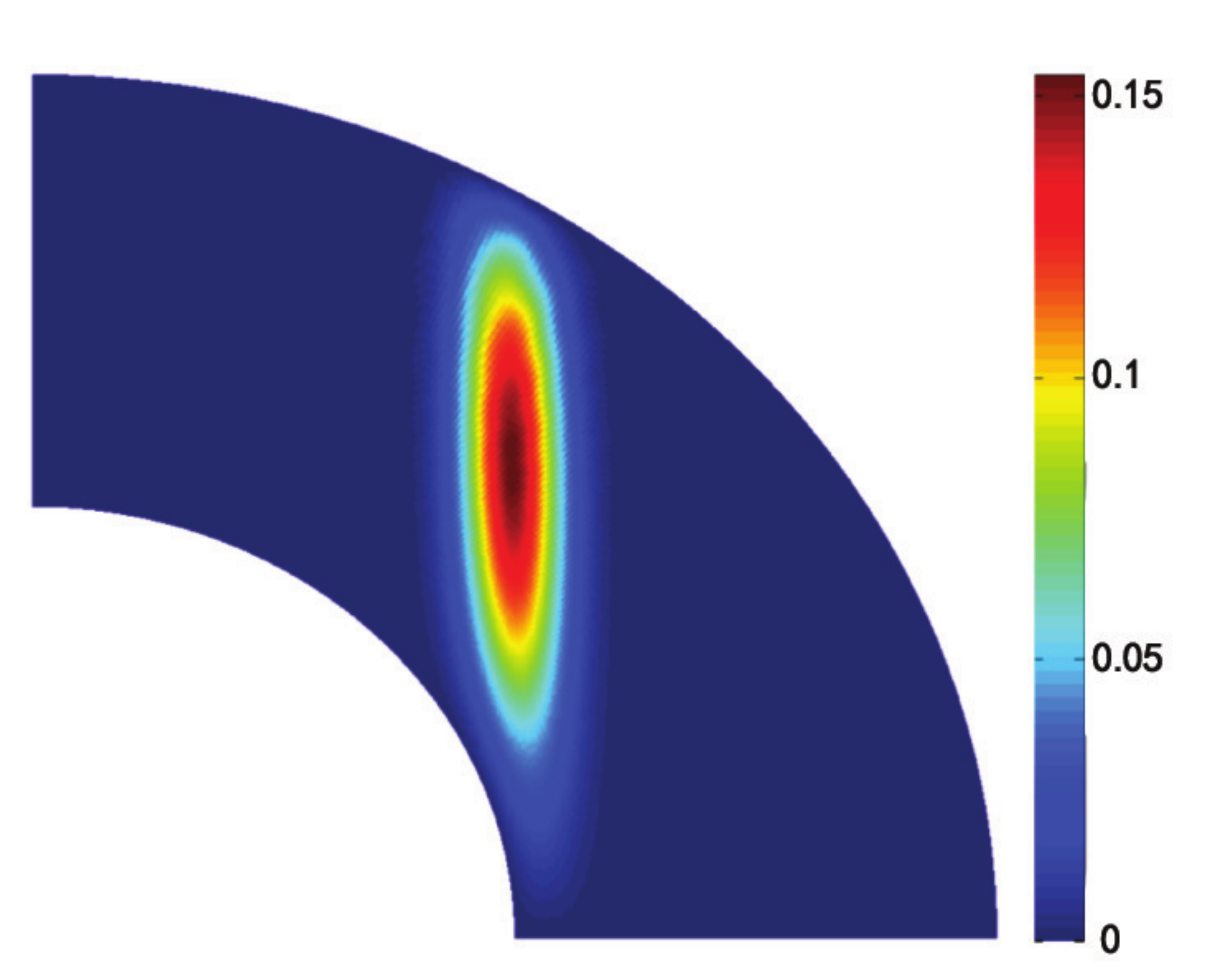}
\quad
\includegraphics[width=\picw]{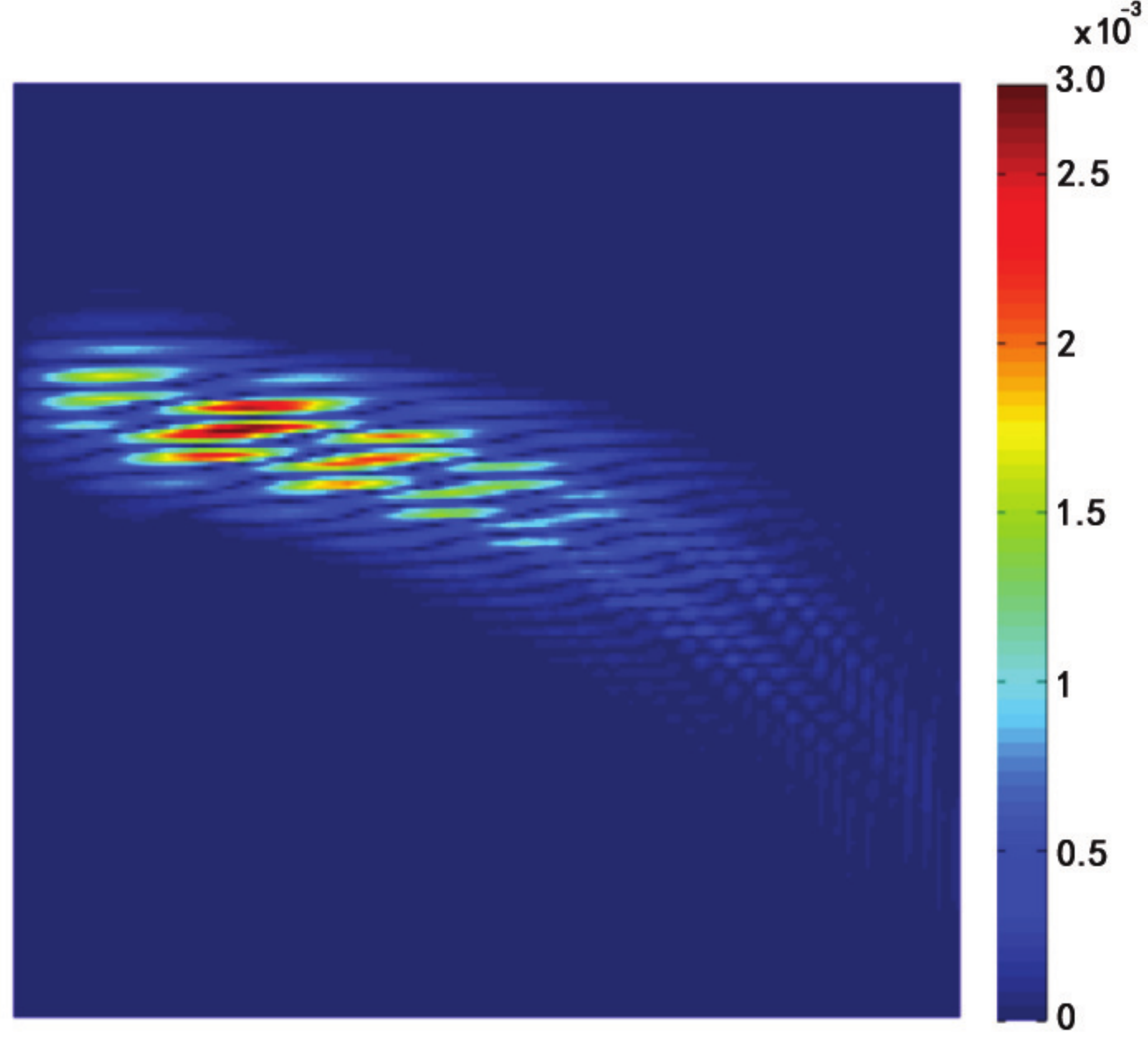}
\\
\includegraphics[width=\picw]{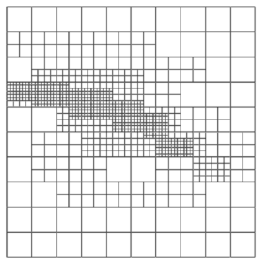}
\quad
\includegraphics[width=\picw]{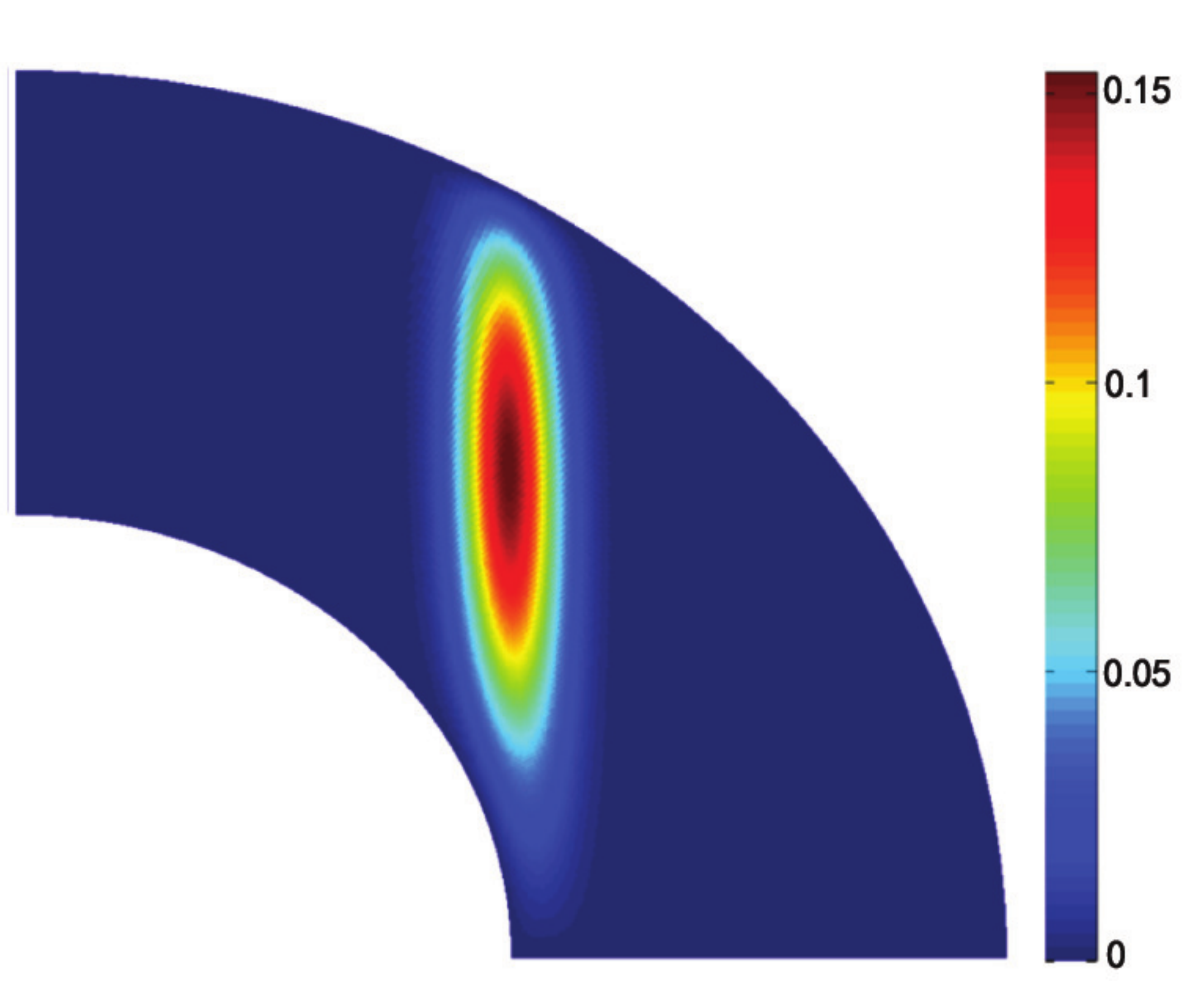}
\quad
\includegraphics[width=\picw]{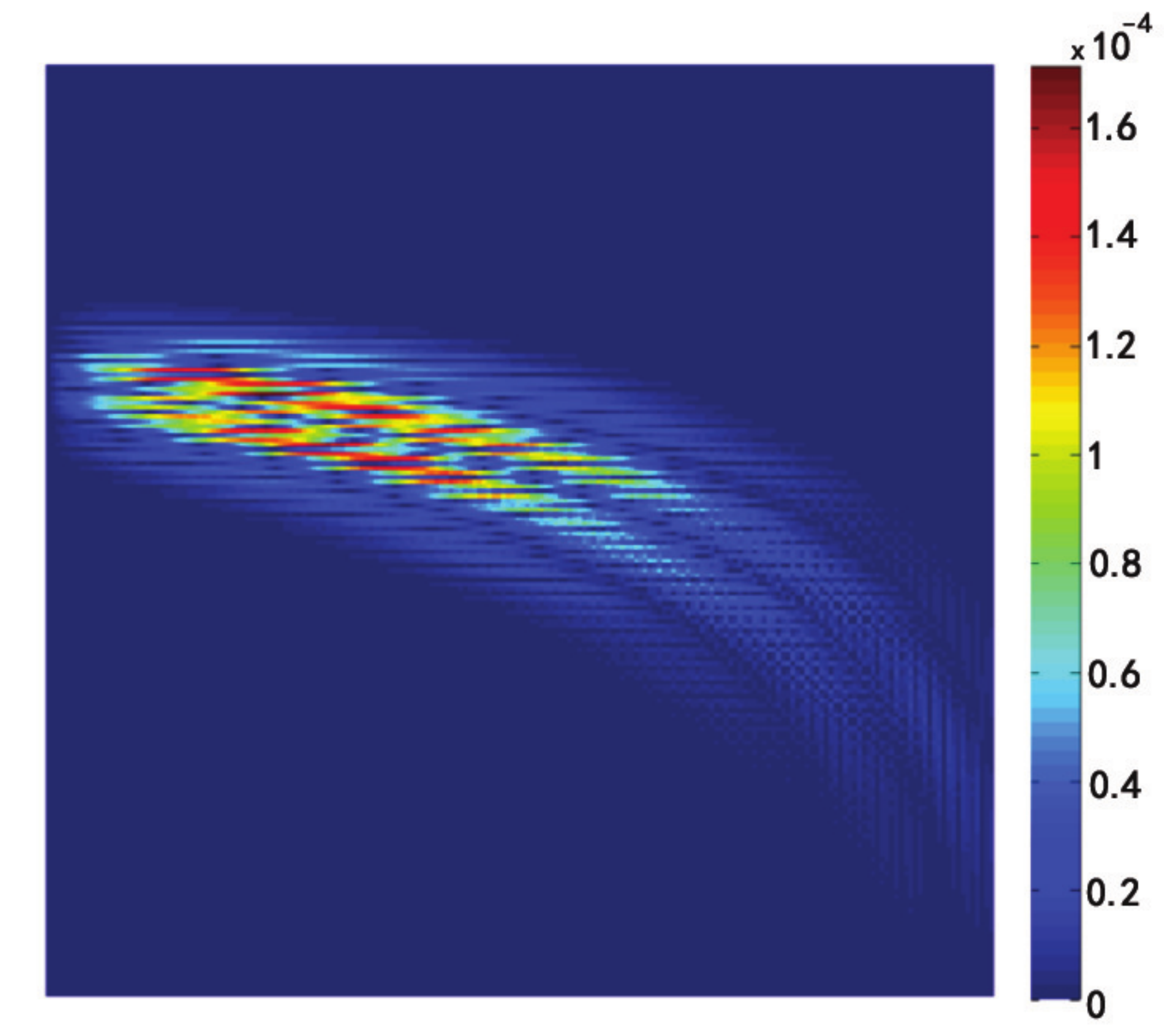}
\\
\caption{GIFT approach for a problem with exact solution having sharp peaks. The geometry representation is based on quadratic $C^{1}$ NURBS, and the numerical solution is based on $C^{1}$ PHT-splines, the latter enabling local refinement. From the first row to the fourth row, we show the T-mesh in the parametric domain (left), the color-map of the numerical solution (middle) and the corresponding error color-map (right) on parametric domain during local refinement operations.}
\label{fig:Annus2}
\end{figure}
\begin{figure}[t]
\centering
\includegraphics[width=3.8in, height=3.6in]{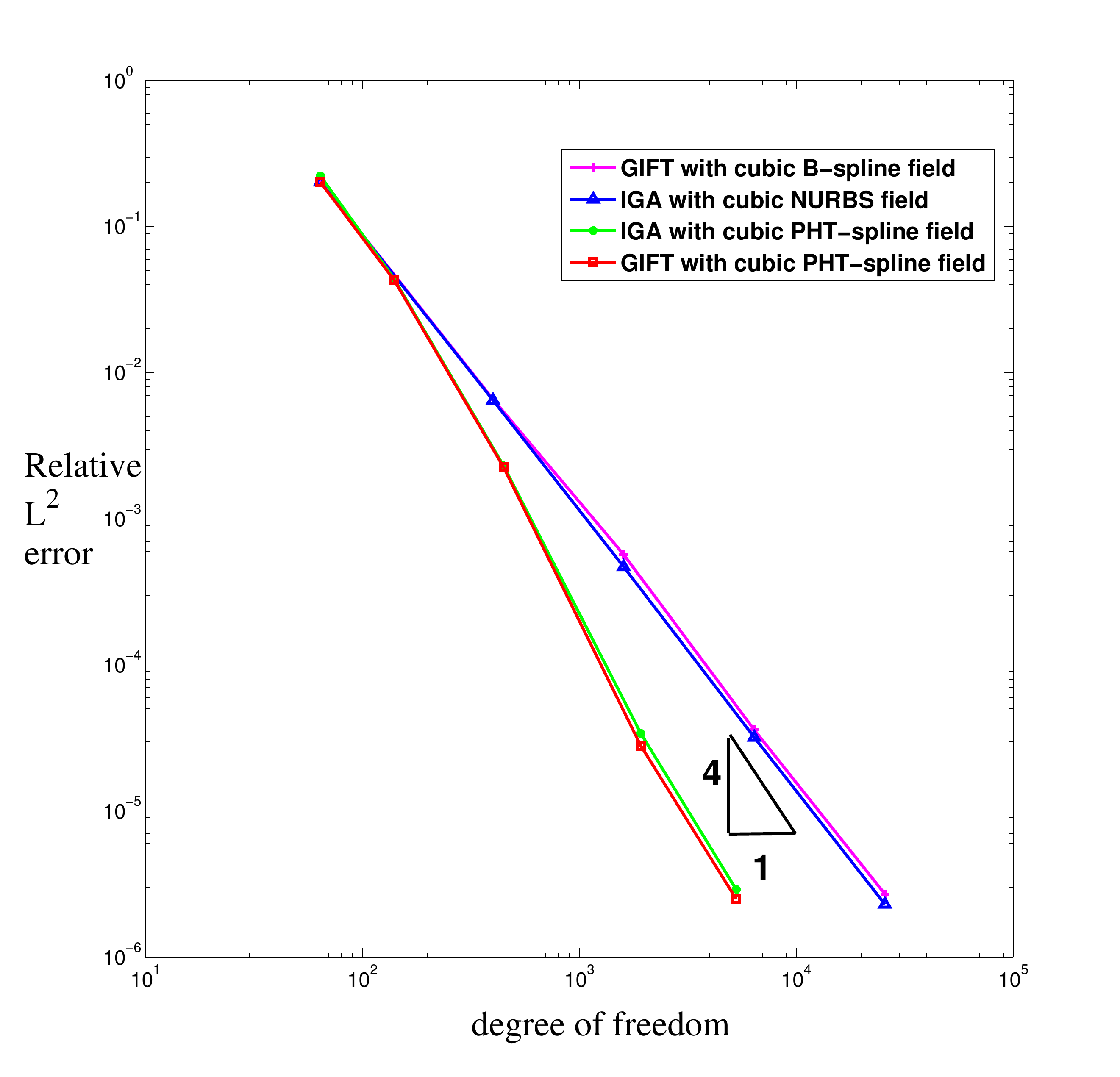}
\caption{Convergence comparison: IGA method with cubic NURBS, IGA method with cubic PHT-splines, GIFT method with cubic B-splines, and GIFT method with cubic PHT-splines (as in Fig.~\ref{fig:Annus2}).}
\label{fig:AnnusPHTerror} 
\end{figure}
\section{Discussion on numerical results}\label{res_disc}

\begin{figure}[!ht]
\begin{center}
\includegraphics[width=0.9\textwidth]{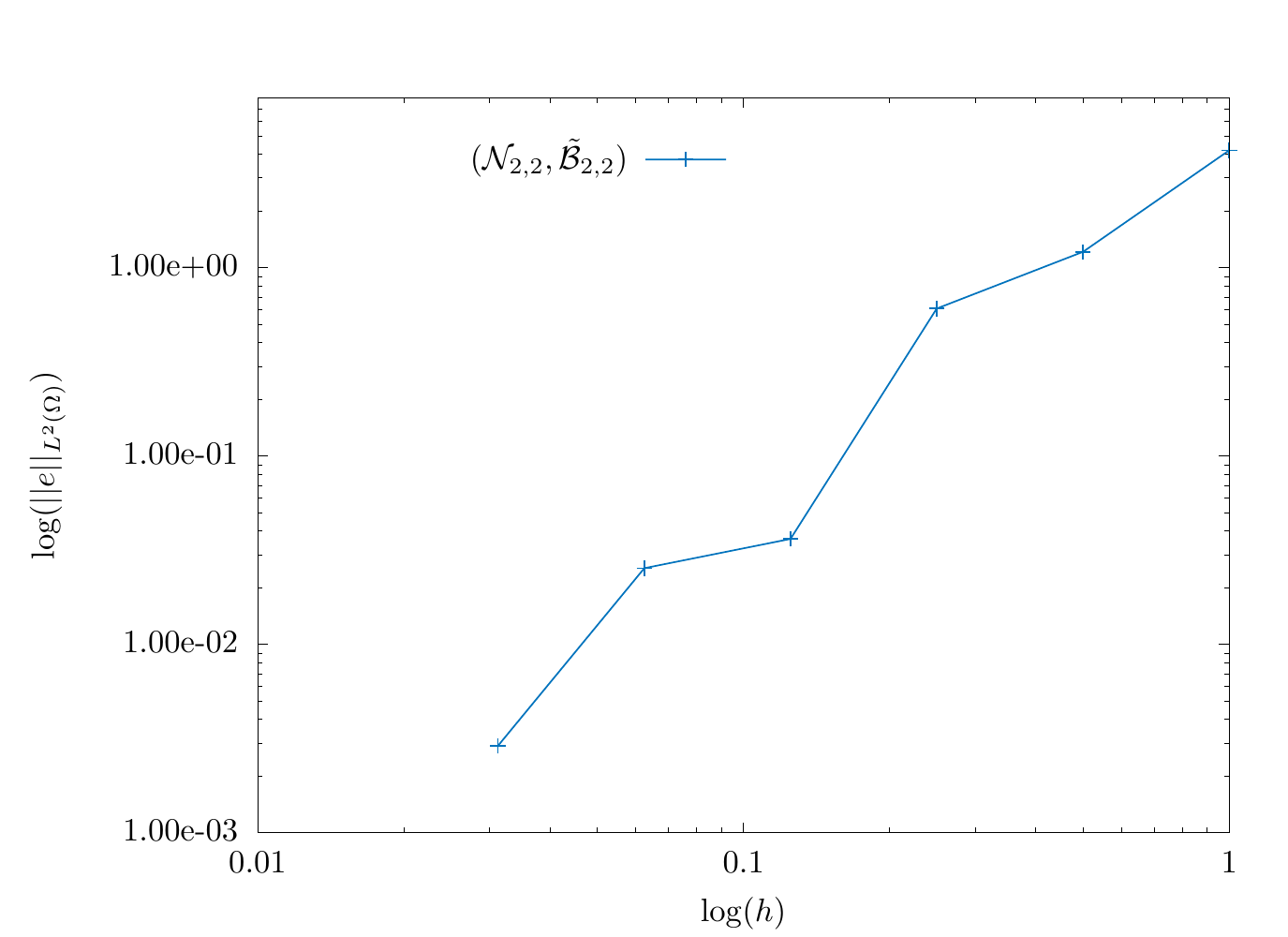}
\caption{Convergence study for Example~\ref{example_elasticity_2} with knot vectors \eqref{knot_vectors_example_divergence_plate}.}
\label{example_elasticity_2_plot_2}
\end{center}
\end{figure}

The numerical results of Section~\ref{patch_tests}, indicate that  a \textbf{sufficient condition} for the combination of the geometry and the solution bases to pass the \textbf{classical patch test} is the requirement that these two bases are equivalent up to the degree elevation \eqref{degree_elevation_geo_exactness} and/or knot insertion operations \eqref{knot_insertion_geo_exactness}, which ensure the preservation of the geometry exactness. %
However, we see from the numerical results in Sections~\ref{sec:laplace}-\ref{sec:laplace_pht}, that despite failing the classical patch test, all  cases presented in Table~\ref{tab:res_tests} exhibit optimal order convergence. This includes the bases derived from the geometry parametrization $A$ and $B$, the non-uniform parametrization $C$, and the B-splines basis $D$ (which cannot represent the geometry exactly). This is not surprising, however. Thanks to \cite[Lemma~3.4-3.5, Theorem~3.1-3.2]{HughesErrorIGA}, with a suitable projector $\Pi_{\mathcal{V}_{h}}$, the optimal order global error estimates hold for every function $v \in H^{l} (\Omega)$. The effect of geometry mapping is included in the definition of the element size $h_{K} = \Vert \nabla F \Vert_{L^{\infty}(Q)} h_{Q}$, where $F$ is the geometry mapping, $Q$ is the element in the parameteric domain, and $h_{Q}$ is the element size in the parameteric domain. Moreover, the constant $C_{\text{shape}}$ appearing in \cite[Lemma~3.4-3.5, Theorem~3.1-3.2]{HughesErrorIGA} depends only on the weight function $w$ and its reciprocal $1/w$ on $\tilde{Q}$, where $\tilde{Q}$ is the support extension of $Q$, and is uniformly bounded with respect to the mesh size. This observation is in alignment with \cite{Stummel-80} that the classical patch test is neither necessary nor sufficient for convergence. Therefore, together with exact geometric representation, a suitable basis (derived-from or related-to geometry parameterization) can be used for optimal order convergence. 
For a quick comparison, in Table~\ref{tab:pt_oc} we recall all the cases studied in this paper.
\begin{table}[!ht]
\renewcommand*{\arraystretch}{1.3}
\begin{center}
\begin{tabular}{|c|c|c|c|c|}
\hline 
Geometry & Solution & Degree & Patch & Optimal \\
parameterization & basis & parity & test & convergence \\
\hline 
$Q_{0}$ & $A_{1}$ & Iso-geometric & \checkmark & \checkmark \\
$Q_{0}$ & $A_{2}$ & Super-geometric & \checkmark & \checkmark \\
$Q_{0}$ & $C_{1}$ & Iso-geometric & $\times$ & \checkmark \\
$Q_{0}$ & $C_{2}$ & Super-geometric & $\times$ & \checkmark \\
\hline \hline
$A_{1}$ & $A_{1}$ & Iso-geometric & \checkmark & \checkmark \\
$A_{1}$ & $A_{2}$ & Super-geometric & \checkmark & \checkmark \\
$A_{2}$ & $A_{1}$ & Sub-geometric & \checkmark & \checkmark \\
\hline \hline
$B_{1}$ & $A_{1}$ & Iso-geometric & \checkmark & \checkmark \\
$B_{1}$ & $A_{2}$ & Super-geometric & \checkmark & \checkmark \\
$B_{2}$ & $A_{1}$ & Sub-geometric & \checkmark & \checkmark \\
\hline \hline
$C_{1}$ & $C_{1}$ & Iso-geometric & \checkmark & \checkmark \\
$C_{1}$ & $C_{2}$ & Super-geometric & \checkmark & \checkmark \\
$C_{2}$ & $C_{1}$ & Sub-geometric & \checkmark & \checkmark \\
\hline \hline
$C_{1}$ & $A_{1}$ & Iso-geometric & $\times$ & \checkmark \\
$C_{1}$ & $A_{2}$ & Super-geometric & $\times$ & \checkmark \\
$C_{2}$ & $A_{1}$ & Sub-geometric & $\times$ & \checkmark \\
\hline \hline
$A_{1}$ & $D_{1}$ & Iso-geometric & $\times$ & \checkmark \\
$A_{1}$ & $D_{2}$ & Super-geometric & $\times$ & \checkmark \\
$A_{1}$ & $D_{0}$ & Sub-geometric & $\times$ & \checkmark \\
\hline
\end{tabular} 
\caption{Summary of patch tests and optimal convergence. See Section~\ref{sec:geo_param} for the notations $Q$, $A$, $B$, $C$, and $D$. For the degree parity, we used the naming convention with respect to the geometry, see Table~\ref{tab:naming}.}
\label{tab:pt_oc}
\end{center}
\end{table}
Nevertheless, it is important to exercise caution while devising a basis for the numerical solution. In all our test cases of Sections~\ref{sec:laplace}-\ref{sec:lin_elas_3d}, the solution basis was constructed based on the same knot vectors as the geometry parametrization. This assures the continuity of the geometry parametrization within solution elements. For Example~\ref{example_elasticity_2}, we now consider the coarsest basis consisting of B-splines of degree $2 \times 2$ defined on the knot vectors (denoted by $\tilde{\mathcal{B}}_{2,2}$):
\begin{equation}
\Sigma = \{0, 0, 0, 1, 1, 1\}, \quad \Pi = \{0, 0, 0, 0.166667, 1, 1, 1\}.
\label{knot_vectors_example_divergence_plate}
\end{equation}
In this case, the singular point $(x = L, ~y = L)$ is inside an element of any mesh for the numerical solution. As it can be seen from Fig.~\ref{example_elasticity_2_plot_2}, the solution does not exhibit the expected convergence rate.
A thorough mathematical derivation of the presented approach is a subject of future research, but it is recommended to avoid those combination of geometry-solution bases where the continuity of the geometry parametrization is violated in the solution elements.

\subsection{Naming convention}\label{sec:naming}

Let $p_{u}$ and $p_{g}$ denote the degrees of the basis functions for the solution, and the geometry, respectively. In the standard FE context, since the primary quantity of interest is the numerical solution, the naming convention is with respect to the parameter $p_{u}$. However, in IGA, the primary quantity is the geometry. Accordingly, for better readability, we propose the naming convention with respect to the parameter $p_{g}$, as presented in Table~\ref{tab:naming}. %
\begin{table}[!ht]
\begin{center}
\begin{tabular}{|c|c|c|c|}\hline
& $p_{u} = p_{g}$ & $p_{u} < p_{g}$ & $p_{u} > p_{g}$ \\ \hline
FEM & Iso-parametric & Super-parametric & Sub-parametric \\
GIFT & Iso-geometric & Sub-geometric & Super-geometric \\
\hline
\end{tabular}
\end{center}
\caption{Naming convention}
\label{tab:naming}
\end{table}
Note that, in standard FEM, the use of $p_{u} < p_{g}$ (super-parameteric) case is not recommended \cite[P.172]{ZienkiewiczTZ-FEM_V1}. However, the results presented in Section~\ref{example_3_Laplace} show that, with exact geometry representation, super-parametric approximations can also deliver optimal orders of convergence.

\section{Conclusions}\label{conclusions}

We presented a method which relaxes the requirement for a tight coupling between the spaces for the representation of the geometry and the approximation of the field variables, but retains both geometrical exactness and the ability to operate directly from CAD files.

The increased independence of the choice of the boundary parameterization and field approximations enables local refinement to capture sharp gradients, without modifying the geometry parameterization. We formulate recommendations regarding the relative choice of discretizations which, as our numerical results indicate, yield optimal convergence.

Important future work includes a detailed mathematical analysis of the method to prove a priori error estimates, and the investigation of this approach for other types of partial differential equations such as wave propagation. A detailed numerical analysis of the suitability of GIFT in a boundary element approach will follow the work of \cite{beer2015simple, marussig2015fast, beer2016advanced, beer2016isogeometric, GeneralizedIGABEMBeer2016, TrimmedBeer2016, Zechner2016212, beer2017isogeometric}.


%
\appendix
\section{PHT-splines}\label{sec:app_PHT}
In this section we present some details on PHT-splines, which were used for the numerical  results in Section~\ref{sec:laplace_pht}.
PHT-splines, proposed by Deng et al. \cite{Deng200876}, are piecewise bicubic polynomials over a hierarchical T-mesh, which inherit the advantageous properties of T-splines. Unlike T-splines, PHT-splines are non-rational polynomial splines, and the refinement algorithm of PHT-splines is local and simple. The blending functions of PHT splines are linearly independent, an important property needed for finite element approximations.
A T-mesh is a rectangular partition of a planar domain with grid lines parallel to the boundary of the domain which allows T-junctions. In T-meshes, the end points of each grid line must lie on two other grid lines, and each cell or facet in the grid must be a rectangle. If a vertex is inside of the domain, it is called an \emph{interior vertex}, otherwise, it is called a \emph{boundary vertex}. There are two types of interior vertices, namely \emph{crossing vertices} (i.e., it possesses valency 4) and \emph{T-vertices} with valency 3, respectively.
Let $\Omega \in \RR ^{m}$ be a rectangular domain with boundary $\partial\Omega$. Denote by $\TT = \cup \Kc$ a hierarchical T-mesh over domain $\Omega$, where $\Kc$ is a cell of the mesh. We further define a spline space
\[
S(p, q, \alpha, \beta, \TT ) 
:= \{s(\xi,\eta) \in \Cc^{\alpha,\beta}(\Omega) | s(\xi,\eta) \in \PP_{pq}, ~\text{for any element} ~\Kc  \in \TT\},
\]
where $\PP_{pq}$ is the space of all the bivariate polynomials with degree $(p, q)$, and the space $\Cc^{\alpha,\beta}(\Omega)$ consists of all continuous bivariate spline functions up to order $\alpha$ in the $\xi$-direction and order $\beta$ in the $\eta$-direction. The dimension formula for the spline space $S(p, q, \alpha, \beta, \TT )$, with $p \geq 2\alpha+1$ and $q \geq 2\beta+1$, has already been provided in \cite{Deng200876}. For the cubic PHT-splines space, the dimension formula can be written as
\[
\emph{Dim} ~\Sc(3, 3, 1, 1, \TT ) = 4(V^b + V ^{+}),
\]
where $V^b$ stands for boundary vertices and $ V^{+}$ stands for interior crossing vertices. From the dimension formula, four basis functions are associated with each basis vertex (boundary vertex or crossing vertex), and they can be built with a hierarchical approach.
For the initial level, i.e., level $0$, denoted as $\TT_0$, the standard bi-cubic tensor-product B-splines are used as basis functions. For simplicity, we set the initial mesh to be a uniform rectangular grid. Suppose that the grid is $[\xi_1, \xi_2, \xi_3, ..., \xi_s] \times [\eta_1, \eta_2, \eta_3, ..., \eta_t]$. Since all the vertices are either crossing vertices or boundary vertices, there are four basis functions to be defined on any vertex $(\xi_i, \eta_j)$. Each basis function at $(\xi_i, \eta_j)$ has support $[\xi_{i-1}, \xi_{i+1}] \times [\eta_{j-1}, \eta_{j+1}]$. These four basis functions are defined to be the B-splines basis functions with knots
\begin{align*}
[\xi_{i-1}, \xi_{i-1}, \xi_i, \xi_i, \xi_{i+1}] \times [\eta_{j-1}, \eta_{j-1}, \eta_j, \eta_j, \eta_{j+1}], \\
[\xi_{i-1}, \xi_i, \xi_i, \xi_{i+1}, \xi_{i+1}] \times [\eta_{j-1}, \eta_{j-1}, \eta_j, \eta_j, \eta_{j+1}], \\
[\xi_{i-1}, \xi_i, \xi_i, \xi_{i+1}, \xi_{i+1}] \times \eta_{j-1}, \eta_j, \eta_j , \eta_{j+1}, \eta_{j+1}], \\
[\xi_{i-1}, \xi_{i-1}, \xi_i, \xi_i, \xi_{i+1}] \times [\eta_{j-1}, \eta_j, \eta_j, \eta_{j+1}, \eta_{j+1}],
\end{align*}
respectively, such that their function values and derivatives vanish outside $[\xi_{i-1}, \xi_{i+1}] \times [\eta_{j-1}, \eta_{j+1}]$.
In the GIFT framework with PHT-splines, the computational domain is in the NURBS space whereas the numerical solution is in the PHT-splines form. Firstly, we construct the parametric domain of the PHT-splines model to represent the numerical solution. We can then get the
initial numerical solution, in which the unknown control variables can be solved by the method presented in Section~\ref{sec:maths}. By using an a posteriori error estimation technique, the supporting cell with large errors in the numerical solution (in the parametric domain) can be marked, and local $h$-refinement is performed only on the numerical solution. Several local refinement steps can be performed until the desired error level is achieved.
\begin{figure}[t]
\centering
\begin{minipage}[t]{3.05in}
\centering
\includegraphics[width=2in]{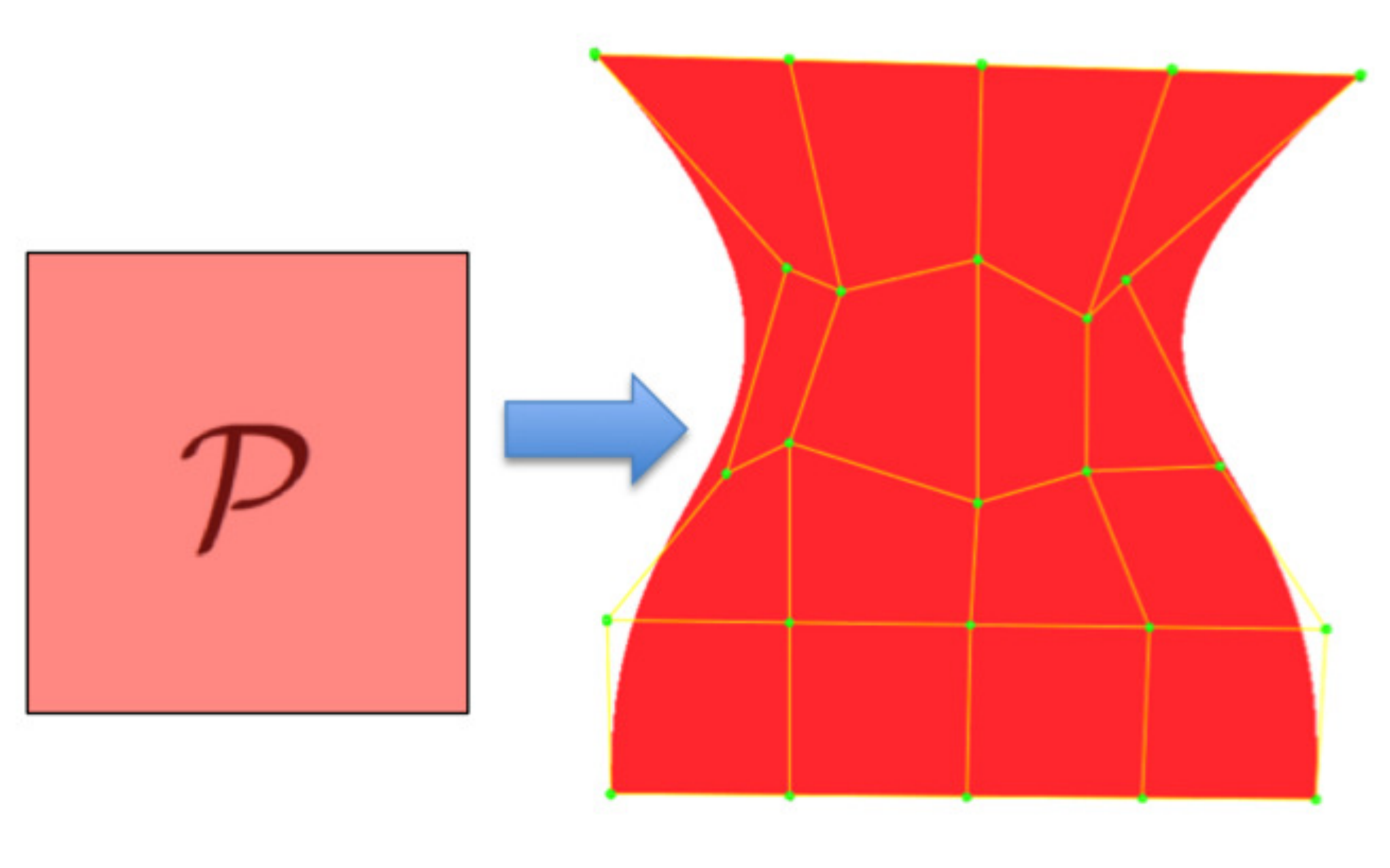} \\ (a)
\end{minipage}
\begin{minipage}[t]{3.05in}
\centering
\includegraphics[width=3.0in]{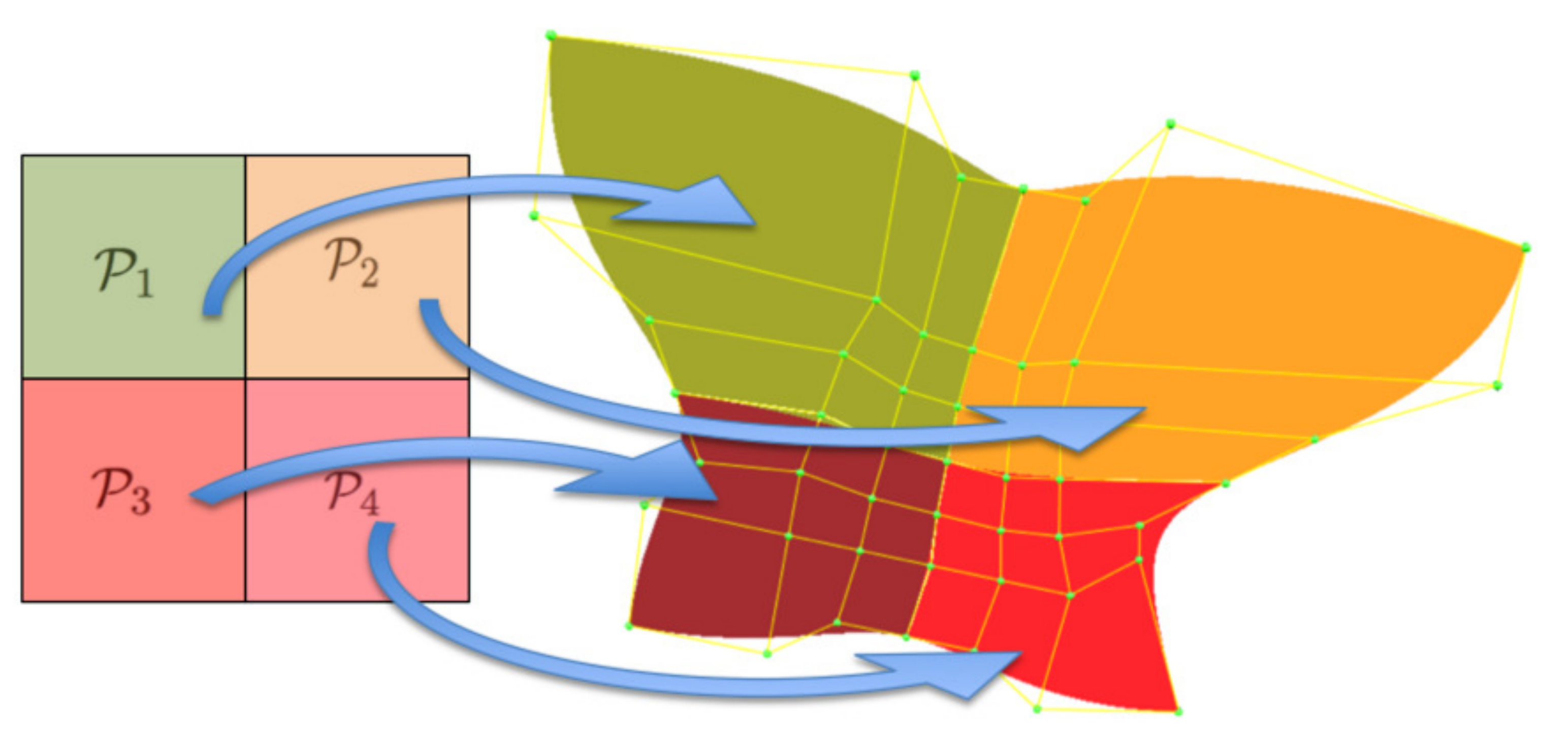} \\ (b)
\end{minipage}    \\
\begin{minipage}[t]{4in}
\centering
\includegraphics[width=4in]{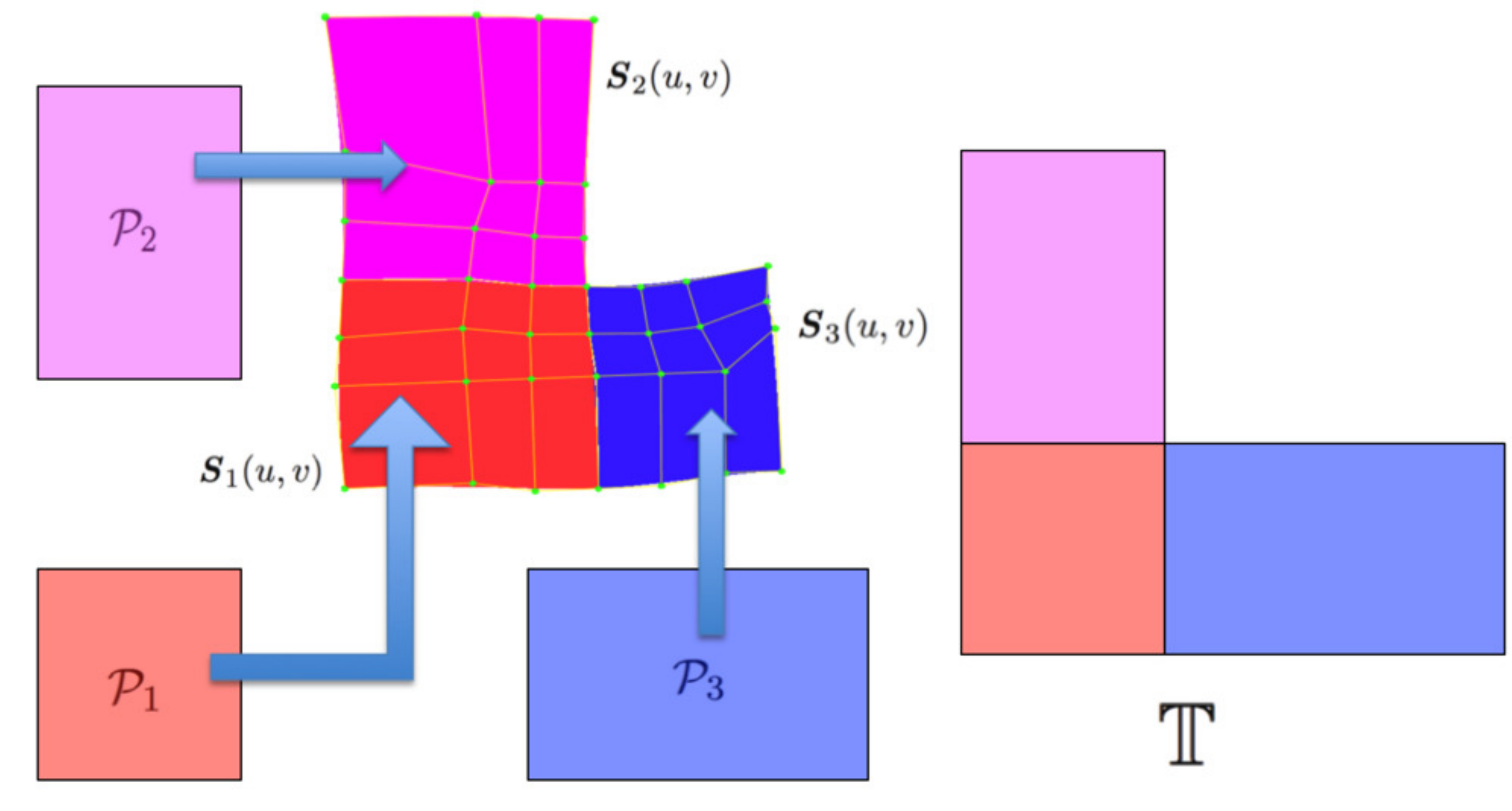}\\ (c)
\end{minipage}  \\
\caption{Parametric T-mesh $\TT$ construction for GIFT with PHT-splines:
 (a) case of a single NURBS patch;
 (b) case of multi-patch in which the parametric domain of each patch forms a quad-mesh;
 (c) other cases in which reparameterization is required.}
\label{fig:PHTpreprocess}
\end{figure}
\subsection{Construction of common parametric domain for NURBS geometry and PHT-splines solution}

As a first step, we need to construct a parametric domain of the numerical solution. This can be seen as a preprocessing stage. Depending on the parametric domain of the NURBS patches under consideration, two kinds of operations are proposed in this step:
\begin{itemize}
\item If the computational domain is made of a single NURBS patch (Fig.~\ref{fig:PHTpreprocess}(a)), or multi-patches in which the parametric domain of each patch forms a quad-mesh as presented in Fig.~\ref{fig:PHTpreprocess}(b), then the parametric domain of the numerical solution is constructed as the partition-mesh formed by the knot lines on the parametric domain of the NURBS patches.
\item For other cases, such as the examples shown in Fig.~\ref{fig:PHTpreprocess}(c), we need to construct the initial parametric mesh $\TT$ according to the topological connection information of the NURBS patches. A reparameterization operation should be performed in this case as described below.
\end{itemize}
Here we present an example to show the reparameterization method. Suppose that the patch $\sb_1(u,v)$ in Fig. \ref{fig:PHTpreprocess}(c) has parametric domain $\Pc_1=[a_1,b_1] \times [c_1,d_1]$, the patch $\sb_2(u,v)$ has parametric domain $\Pc_2=[a_1,b_1] \times [c_2,d_2]$, and the patch $\sb_3(u,v)$ has parametric domain $\Pc_3=[a_3,b_3] \times [c_1,d_1]$. Then according to the topological information of all three patches, the parametric domain for the PHT-splines numerical solution should be constructed as shown in Fig.~\ref{fig:PHTpreprocess}(c). In other words, the parametric domain of $\sb_1(u,v)$ does not change, while the parametric domain of $\sb_2(u,v)$ changes to $[a_1,b_1] \times [d_1, d_1+(d_2-c_2)]$, and the parametric domain of $\sb_3(u,v)$ changes to $[b_1, b_1 + (b_3 - a_3)] \times [c_1, d_1]$.
In order to maintain the geometry of the patches during the transformation of parametric domain, a reparameterization technique should be adopted to obtain the new parametric representation of each patch. Suppose that the initial parametric domain of $\sb(u,v)$ is $[a,b] \times [c,d]$, we can then use the following parameter transformation to achieve a new parameterization $\boldsymbol{F}(\xi,\eta)$ with parametric domain $[e,f] \times [g,h]$
\begin{align}
u(\xi,\eta) & = \frac{1}{f - e}[(f - \xi)a + (\xi - e) b] \\
v(\xi,\eta) & = \frac{1}{h - g}[(h - \eta)c + (\eta - g) d]
\end{align}
\begin{remark}
There is actually no need to derive the explicit parametric representation of the reparameterized surface $\boldsymbol{F}(\xi,\eta)$. The derivative information required in the solving stage, as shown in (\ref{eqn:Kij_x}), can be evaluated from the original parameterization $\sb(u,v)$ through the Jacobian transformation matrix.
\end{remark}
After the planar T-mesh $\TT$ is constructed as the parametric domain of PHT-splines, we can write the initial formula of the numerical solution in PHT-splines form as follows
\begin{equation}\label{eq:PHT}
\Tc (\boldsymbol{\xi}) = \sum_{i=1}^{4n} \, T_{i} \, M_i(\boldsymbol{\xi}),
\end{equation}
where $n$ is the number of basis vertices on $\TT$, and $T_{i}$, $i=1, 2, \cdots, 4n$, are the control variables to be solved.
\subsection{Solution process and local refinement}

After constructing PHT-splines representation of the numerical solution, the control variables in (\ref{eq:PHT}) is obtained by solving the problem as introduced in Section \ref{sec:maths}.
In order to obtain a solution with desired accuracy, refinement operations are often needed to increase the approximation power of the solution space locally, and thereby optimize the computational expenses required for a given accuracy level. Contrary to the classical IGA, the refinement operation in GIFT is only performed on the numerical solution, while the parameterization of the computational domain remains the same during the refinement process. Choosing a PHT-splines approximation offers a natural local refinement.  Local refinement is performed on the sub-patches of the numerical solution, on which the local error indicator is larger than a given marking threshold. The corresponding working flow is summarized in Algorithm~\ref{alg:refine}.

\begin{algorithm}[t]
\caption{Local refinement scheme in GIFT using PHT-splines }
\noindent\textbf{Input}: Planar NURBS parameterization  $\boldsymbol{F} (\xi,\eta) $ of computational domain $\Omega$ \\
\textbf{Output}: PHT-splines numerical solution
\begin{algorithmic}[1]
\State Compute the PHT-splines numerical solution $\Tc (\xi,\eta)$ of model problem \eqref{weak_form} and \eqref{weak_form_Poisson} over the given NURBS parameterization $\boldsymbol{F} (\xi,\eta) $ by GIFT.
\State Calculate the local error indicator $e_\Kc$ patch by patch for the numerical solution (See Section~\ref{sec:res_error}).
\State Mark the parametric cells to be refined by a mean-value marking algorithm (See Section~\ref{sec:mark_elm}).
\State Subdivide the marked cells into four sub-cells on the parametric domain of the  PHT-splines solution.
\State Construct the PHT-splines basis functions over the refined T-mesh of the parametric domain.
\State Compute the new numerical solution $\Tc (\xi,\eta)$ in the refined PHT-splines space.
\State Repeat the above refinement steps until the estimated error is less than a given threshold.
\end{algorithmic}
\label{alg:refine}
\end{algorithm}
\subsubsection{Residual-based error indicator}\label{sec:res_error}

Suppose that $u_{h}$ is the PHT-splines solution of the problem \eqref{weak_form} and \eqref{weak_form_Poisson} by using the GIFT method, and $u$ is the exact solution. Let $e_{h} = u - u_{h}$ be the error of the GIFT approximation $u_{h}$. As the refinement operation in GIFT is only performed on the numerical solution, in order to determine the parametric cell to be refined, it is required to give an error indicator on each cell of the T-mesh in the parametric domain of the numerical solution, rather than an error indicator on the sub-patch with respect to the knot span on the NURBS parameterization $\boldsymbol{F}(\xi, \eta)$.
Suppose that ${\Kc}$ is the cell on the T-mesh $\TT$ of the parametric domain $\Pc$ for the PHT-splines solution $u_{h}$. The residual-based a posteriori error estimate $\|e_h\|^2_{\Pc}$ over the parametric domain $\Pc$ proposed in \cite{xu:cmame,xu13cadb} can be rewritten as follows:
\begin{equation}\label{eq:finalerror}
\|e_h\|^2_{\Pc} \le C \sum_{{\Kc}\in \TT} h_{\Kc}^2 \| f({\boldsymbol x}) + \Delta u_{h} ({\boldsymbol x})\|^2_{L^2({\Kc})},
\end{equation}
where $\boldsymbol{x}={\boldsymbol{F}} (\xi,\eta)=(x(\xi,\eta),y(\xi,\eta))$, $C$ is a positive constant, and $h_{\Kc}$ is the circumference of the sub-patch in the NURBS parameterization $\boldsymbol{F}(\xi,\eta)$ of $\Omega$ with respect to the cell $\Kc$ in the parametric domain $\Pc$.
From (\ref{eq:finalerror}), we can obtain the local error indicator $e_{\Kc}$ on each parametric cell $\Kc$ as follows
\begin{equation}
e_{\Kc}= \sqrt{h_{\Kc}^2\| f({\boldsymbol x}) + \Delta u_{h} ({\boldsymbol x})\|^2_{L^2({\Kc})}}.
\end{equation}
If the parametric cell ${\Kc}$ is written as $[\xi_0,\xi_1] \times [\eta_0,\eta_1]$, we have then
\begin{equation} \label{eq:length}
h_{\Kc}=\int_{\xi_0}^{\xi_1} \|\boldsymbol{F}_{\xi}(\xi,\eta_0)\|_2 ~d \xi +
\int_{\xi_0}^{\xi_1} \|\boldsymbol{F}_{\xi}(\xi,\eta_1)\|_2 ~ d \xi + \int_{\eta_0}^{\eta_1}
\|\boldsymbol{F}_{\eta}(\xi_0,\eta)\|_2 ~d \eta+ \int_{\eta_0}^{\eta_1}  \|\boldsymbol{F}_{\eta}(\xi_1,\eta)\|_2 ~d \eta.
\end{equation}

\subsubsection{Marking strategies}\label{sec:mark_elm}

The local refinement in GIFT requires a marking strategy to decide which elements  should be refined. That is, we should determine a subset $\tilde{\TT}$ of the parametric T-mesh $\TT$ consisting of all those cells $\Kc$ that must be refined when the value of $e_{\Kc}$ is bigger than certain threshold. To determine the set $\tilde{\TT}$, we use the mean-value strategy shown in Algorithm~\ref{alg:mark}.
\begin{algorithm}[t]
\caption{Marking Algorithm. (Mean-value  strategy) }
\noindent\textbf{Input}: Parametric T-mesh $\TT$ of numerical solution, error estimates $e_{\Kc}$ for all parametric cells ${\Kc} \in \TT$.\\
\textbf{Output}: A subset $\tilde{\TT}$ of marked parametric cells to be refined
\begin{algorithmic}[1]
\State Compute $e_{\TT, mean}=\dfrac{\sum_{{\Kc}\in \TT} {e_{{\Kc}}}}{N}$, where $N$ is the number of parametric cells in $\TT$.
\State If $e_{\Kc} \geq e_{\TT, mean}$, mark $\Kc$ for refinement and put it into the set $\tilde{\TT}$.
\end{algorithmic}
\label{alg:mark}
\end{algorithm}
For some physical problems, highly varying distribution of estimated errors may be observed in the following cases (I) very few cells having extremely large estimated errors, (II) some cells having extremely small estimated errors, and (III) other cells having estimated errors that are much smaller than that of (I) while much larger than that of (II). We group the cells into three categories (I), (II) and (III), respectively. In this case, if only the parametric cells in the first group are refined by using the above marking algorithm, it would deteriorate the performance of the local $h$-refinement method. In this paper, a simple modification is proposed. For a given small percentage $\epsilon$,
the $\epsilon\%$ parametric cells are firstly marked with largest estimated error for refinement, and then the mean-value marking approach is further applied to the remaining parametric cells.

\subsubsection{Local refinement}
In order to achieve accurate and efficient simulation results, local refinement should be performed on the marked parametric cells. In this part, the refinement rules of PHT-splines will be introduced.

Suppose that the T-mesh on the parametric domain at level $k$ is denoted by $\TT_k$, and the PHT-splines basis functions on $\TT_k$ are constructed as ${M_{j}^k }$, $j = 1, ..., d_k$. Then the basis functions on $\TT_{k+1}$ can be constructed as follows: some basis functions are from the modifications of the former basis functions on $\TT_{k}$, and others are from the new basis functions associated with the new basis vertices of $\TT_{k+1}$. We represent a PHT-splines basis function by specifying 16 \Bezier coefficients in each cell within the compact support of the basis function. When a cross vertex is added in a cell, the cell can be refined into four subcells. Each subcell supports the original basis function, and also has 16 \Bezier coefficients. Five new vertices are obtained by adding a cross vertex, and some new basis vertices are also introduced. Hence, for the old basis functions, all the \Bezier coefficients associated with the new basis vertices should be reset to zero. The new basis vertices are introduced from two sources, i.e., some are crossing vertices while others are T-vertices from the previous level. The latter become basis vertices as the addition of cross points to the neighboring cells. The new basis functions can be further constructed over their supporting cells as in the initial level. For further details, the readers can refer to \cite{Deng200876}.
Note that in the proposed GIFT framework, the local refinement is only performed on the numerical solution, hence we only need to update the PHT-splines basis functions on the locally-refined T-mesh for the new numerical solution. It is not necessary to derive the updated control variables from the old control variables after local refinement, and the new control variables for the refined numerical solution can be obtained by re-solving the PDE in the new PHT-splines space over the refined T-mesh.

\section*{Acknowledgements}

G. Xu is supported by the National Nature Science Foundation of China under Grant Nos. 61472111, Zhejiang Provincial Natural Science Foundation of China under Grant Nos. LR16F020003, LQ16F020005,  and the Open Project Program of the State Key Lab of CAD\&CG (A1703), Zhejiang University.

S.P.A. Bordas and S. Tomar thank the financial support of the European Research Council Starting Independent Research Grant (ERC Stg grant agreement No. 279578) entitled “Towards real time multiscale simulation of cutting in non-linear materials with applications to surgical simulation and computer guided surgery”, and the support of the Luxembourg National Research Funds INTER/FWO/15/10318764 and INTER/MOBILITY/14/8813215/CBM/Bordas. E. Atroshchenko and S.P.A. Bordas acknowledge the financial support of University of Luxembourg through its Computational Sciences research priority.

%

\end{document}